\documentclass[12pt,reqno]{amsart}
\usepackage{enumerate,color,amsmath,amssymb,amsfonts,bm}
\usepackage{frcursive,fullpage,graphicx}

\newtheorem{Lemma}{Lemma}
\newtheorem{Remark}{Remark}
\newtheorem{Theorem}{Theorem}
\newtheorem{Proposition}{Proposition}
\newtheorem{Corollary}{Corollary}
\newtheorem{Conjecture}{Conjecture}

\def\bbone{{\mathchoice {\rm 1\mskip-4mu l} {\rm 1\mskip-4mu l}
{\rm 1\mskip-4.5mu l} {\rm 1\mskip-5mu l}}}

\begin{document}

\author{Abdelmalek Abdesselam}
\address{Abdelmalek Abdesselam, Department of Mathematics,
P. O. Box 400137,
University of Virginia,
Charlottesville, VA 22904-4137, USA}
\email{malek@virginia.edu}

\dedicatory{Dedicated to the memory of Roland S\'en\'eor (1938--2016)}

\title{A Second-Quantized Kolmogorov-Chentsov Theorem via the Operator Product Expansion}

\begin{abstract}
We establish a direct connection between two fundamental topics: one in probability theory and one in quantum field theory. The first topic is the problem of pointwise multiplication of random Schwartz distributions which has been the object of recent progress thanks to Hairer's theory of regularity structures and the theory of paracontrolled distributions introduced by Gubinelli, Imkeller and Perkowski.
The second topic is Wilson's operator product expansion which is a general property of models of quantum field theory and a cornerstone of the bootstrap approach to conformal field theory. Our main result is a general theorem for the almost sure construction of products of random distributions by mollification and suitable additive as well as multiplicative renormalizations. The hypothesis for this theorem is the operator product expansion with precise bounds for pointwise correlations. We conjecture these bounds to be universal features of quantum field theories with gapped dimension spectrum.
Our theorem can accommodate logarithmic corrections, anomalous scaling dimensions and even lack of translation invariance. However, it only applies to fields with short distance singularities that are milder than white noise. As an application, we provide a detailed treatment of a scalar conformal field theory of mean field type, i.e., the fractional massless free field
also known as the fractional Gaussian field.
\end{abstract}

\maketitle

\tableofcontents 

\section{Introduction}

\subsection{The problem of pointwise multiplication for Schwartz distributions}\label{introtointro}

As is well-known, the pointwise multiplication of Schwartz distributions is, in general, impossible~\cite{SchwartzImp}.
In accordance with the intuition expressed in~\cite[p. 115]{SchwartzBook1}, any deterministic theorem to this effect must involve
a ``compensation principle''. Namely, the regularity of one of the factors must compensate for the other's lack of regularity.
H\"{o}rmander's celebrated theorem~\cite[Theorem 8.2.10]{Hormander} using wave front sets 
is a beautiful implementation of this compensation principle, in a direction-wise manner.
Another instance of this principle is the multiplication theorem in~\cite[\S2.8.1]{BahouriCD}
which uses Bony's paraproducts~\cite{Bony}. Let $\mathcal{C}^{\alpha}(\mathbb{R}^d)$ denote the inhomogeneous
Besov space ${\rm B}_{\infty,\infty}^{\alpha}(\mathbb{R}^d)$. Then the product of smooth functions continuously extends
to $\mathcal{C}^{\alpha}(\mathbb{R}^d)\times\mathcal{C}^{\beta}(\mathbb{R}^d)$,
provided $\alpha+\beta>0$ holds. In this article, however, the main issue being addressed is the problem of multiplication
of {\em random} Schwartz distributions which live in the generalized H\"{o}lder spaces $\mathcal{C}^{\alpha}(\mathbb{R}^d)$
with $\alpha$ {\em negative}. Random distributions are not so nice as to sit in relation to each other in a way that satisfies
a compensation principle which would allow one to multiply them via such a deterministic theorem.
A particularly unfavorable case is when the two factors are the same, i.e.,
one is trying to take the square (or higher powers) of a random distribution.
The probabilistic setting thus brings extra difficulties, but it also comes with a precious advantage: one does not have to multiply
all distributions but only {\em almost all} of them in the sense of the underlying probability measure. One only has to aim
for a pathwise and pointwise multiplication. For a long time, the only known rigorous method for doing this was the Wick product
construction (see, e.g.,~\cite{Segal},~\cite[{\S}V.1]{SimonPphi2},~\cite[\S8.5]{GlimmJbook}
or~\cite{DaPratoT,EJS}).
In the last few years, however, there has been tremendous progress on this problem which goes far beyond the Wick product method.
Notable examples of such advances are the theory of regularity structures~\cite{Hairer} and
the theory of paracontrolled distributions~\cite{GubinelliIP}. The main result of this article is a very general theorem for pointwise and pathwise multiplication of random distributions
which can be seen as a useful complement to these two recent theories. Indeed, our result shows that at the heart of the problem
is one of deepest aspects of quantum field theory (QFT): {\em Wilson's operator product expansion} (OPE).
The latter was discovered in~\cite{Wilson64} and first appeared in published form in~\cite{Brandt}
(see also~\cite{Polyakov,Kadanoff1,Kadanoff2}).
For a presentation of the OPE from a physical perspective and also aimed towards a mathematical
audience, see~\cite[Lecture 3]{Witten}. 

For concreteness, let us consider the simplest instance of the problem of pointwise and pathwise multiplication of random distributions:
the squaring of the fractional massless free field (FMFF)
also known as fractional Gaussian field~\cite{LodhiaSSW}.
Let $C$ denote the continuous bilinear form on the Schwartz space $\mathcal{S}(\mathbb{R}^d)$ of rapidly decaying smooth functions
defined by
\[
C(f,g)=\frac{1}{(2\pi)^d}\int_{\mathbb{R}^d} d\xi\ \frac{\overline{\widehat{f}(\xi)}\widehat{g}(\xi)}{|\xi|^{d-2[\phi]}}
\]
where $[\phi]\in (0,\infty)$ is a parameter called the {\em scaling dimension} of the field $\phi$.
Note that in this article we will write integrals as above, i.e., with the volume element
preceding rather than following the integrand. This is hardly avoidable when integrands take several lines to write, as will be the case
in this article. Also, we simply write $d\xi$ instead of $d^d\xi$ as the dimensionality will be clear from the context.
Our Fourier transform normalization convention is
$\widehat{f}(\xi)=\int_{\mathbb{R}^d}dx\ e^{-i\xi\cdot x}f(x)$.
By the Bochner-Minlos Theorem, there is a unique probability measure $\mathbb{P}$
on the space of temperate distributions $\mathcal{S}'(\mathbb{R}^d)$ such that
$\mathbb{E}\ e^{i\phi(f)}=\exp\left(-\frac{1}{2}C(f,f)\right)$
for all test functions $f\in\mathcal{S}(\mathbb{R}^d)$. We use $\phi$ to denote the random distribution in $\mathcal{S}'(\mathbb{R}^d)$.
Here the (canonical) probability space is $(\Omega,\mathcal{F},\mathbb{P})$
with $\Omega=\mathcal{S}'(\mathbb{R}^d)$ equipped with $\mathcal{F}={\rm Borel}(\mathcal{S}'(\mathbb{R}^d))$.
In contrast to the treatment in~\cite{SimonPphi2,GlimmJbook}, our point of view (which follows~\cite{Fernique}) is to see
$\mathcal{S}'(\mathbb{R}^d)$ as a {\em topological space} when endowed with the strong topology.
Note that, depending on what is most convenient, we will use $\phi(f)$ or $\langle \phi,f\rangle$ for the duality pairing between
a distribution $\phi$ and a test function $f$. We will however avoid writing the formal integral
$\int_{\mathbb{R}^d} dx\ \phi(x) f(x)$.
When discussing the theory of kernels, involving several sets of variables, it is necessary to employ notation of this kind which
names yet shames integration variables with the epithet ``dummy''. 
Instead, we will use $\langle \phi(x),f(x)\rangle_{x}$
together with the subscript notation introduced by Schwartz for function spaces. For example, the distribution $\phi(x)$ is said to belong
to $\mathcal{S}'_x(\mathbb{R}^d)$ in order to emphasize the name of the variable.

In order to define the pointwise square ``$\phi^2(x)$'', the most direct approach is to use a mollifier. Take $\rho_{\rm UV}$ or simply
$\rho$ to be a function in $\mathcal{S}(\mathbb{R}^d)$ such that $\int_{\mathbb{R}^d}dx\ \rho(x)=1$.
Fix $L>1$ and for, $r\in\mathbb{Z}$, define the rescaled function $\rho_r(x)=L^{-rd}\rho(L^{-r}x)$. The convolution
$\phi\ast\rho_{r}$ is well defined in a pointwise
manner and gives us an opportunity to practice the previous notation. Indeed, by definition,
$(\phi\ast\rho_r)(x)=\langle \phi(y), \rho_{r}(x-y)\rangle_{y}$.
The result is a function in $\mathcal{O}_{{\rm M},x}(\mathbb{R}^d)$ (see, e.g.,~\cite[Proposition 7, p. 420]{Horvath}). Here, $\mathcal{O}_{\rm M}$ denotes the space of temperate smooth
functions. It is defined by
\[
\mathcal{O}_{\rm M}(\mathbb{R}^d)=\{f\in C^{\infty}(\mathbb{R}^d)\ |\ 
\forall \alpha\in \mathbb{N}^d,\exists k\in\mathbb{N}, \exists K>0, \forall x\in\mathbb{R}^d,
|\partial^{\alpha}f(x)|\le K\langle x\rangle^k \}\ .
\]
For our notations we use the Bourbaki convention $\mathbb{N}=\{0,1,2,\ldots\}$. We denote the Euclidean norm of a
vector $x\in\mathbb{R}^d$ by $|x|$, and we write $\langle x\rangle=\sqrt{1+|x|^2}$ for its ``inhomogeneous norm''.
For $U$ an open subset of $\mathbb{R}^d$, we of course use $C^{\infty}(U)$ to denote the space
of (real-valued) smooth functions on $U$.
Since $(\phi\ast\rho_r)(x)$ converges to $\phi(x)$ in $\mathcal{S}'_x(\mathbb{R}^d)$ when $r\rightarrow -\infty$, it is natural
to try to define the square of $\phi$
as the distribution $\phi^2$ whose action on a test function $f$ would be given by
$\phi^2(f)=\lim_{r\rightarrow -\infty}\int_{\mathbb{R}^d}dx\ \left[(\phi\ast\rho_r)(x)\right]^2 f(x)$.
While the integral makes perfect sense, unfortunately, the limit usually does not.
Nevertheless, one can define a suitable smeared square $\phi^2(f)$ for any $[\phi]>0$
as a Hida distribution~\cite{AlbeverioL} (see~\cite{Segal} for a related result).
The outcome, however, is in general not a true random variable or function on
$\Omega=\mathcal{S}'(\mathbb{R}^d)$ but a ``second-quantized Schwartz
distribution''~\cite[Ch. 8]{HidaBook}.
It is to $L^2(\Omega,\mathcal{F},\mathbb{P})$ what a Schwartz distribution is to $L^2(\mathbb{R}^d)$.
By a ``second-quantized'' Kolmogorov-Chentsov regularity result, we mean showing that what {\em a priori} is such a generalized
functional of the field $\phi$ is in fact a true random variable. In the simple example under consideration, it is easy to show
that if $[\phi]\in \left(0,\frac{d}{4}\right)$,
then such a second-quantized regularity holds.
Indeed,
if one recenters the distribution-valued random variable $\left[(\phi\ast\rho_r)(x)\right]^2\in
\mathcal{O}_{{\rm M},x}(\mathbb{R}^d)\subset\mathcal{S}'_{x}(\mathbb{R}^d)$, then the desired limit exists.
Namely, the correct (Wick) square $\phi^2$ or rather $:\phi^2:$ is given by
\begin{equation}
:\phi^2:(f)=\lim_{r\rightarrow -\infty}
\int_{\mathbb{R}^d}dx\ \left(\left[(\phi\ast\rho_r)(x)\right]^2-\mathbb{E}\left[(\phi\ast\rho_r)(x)\right]^2\right) f(x)
\label{WickSquare}
\end{equation}
with convergence in every $L^p(\Omega,\mathcal{F},\mathbb{P})$, $p\ge 1$, and almost surely.
This is the simplest case of the Wick product construction.

Let us revisit this simple example in a way that gives the OPE flavor of our general theorem. 
When $[\phi]\in \left(0,\frac{d}{2}\right)$, then one has a {\em pointwise representation} for the covariance $C$, i.e.,
\[
C(f,g)=\int_{\mathbb{R}^{2d}}dx\ dy\ \langle \phi(x)\phi(y)\rangle\ f(x)g(y)
\]
where the pointwise correlation $\langle \phi(x)\phi(y)\rangle$
is defined outside the diagonal by
\[
\langle \phi(x)\phi(y)\rangle=\frac{\varkappa}{|x-y|^{2[\phi]}}
\ \ 
{\rm with}
\ \ 
\varkappa=\pi^{\frac{d}{2}}\times 2^{2[\phi]}\times\frac{\Gamma\left([\phi]\right)}{\Gamma\left(\frac{d}{2}-[\phi]\right)}
\]
as shown, e.g., in~\cite[p. 193]{GelfandS}.
If now one considers a higher order moment, say a fourth order one,
\[
\mathbb{E}\left[\phi(f_1)\phi(f_2)\phi(f_3)\phi(f_4)\right]=C(f_1,f_2)C(f_3,f_4)+C(f_1,f_3)C(f_2,f_4)+C(f_1,f_4)C(f_2,f_3)
\]
for $f_1,\ldots,f_4\in\mathcal{S}(\mathbb{R}^d)$, then one also has 
a pointwise representation
\begin{eqnarray*}
\lefteqn{
\mathbb{E}\left[\phi(f_1)\phi(f_2)\phi(f_3)\phi(f_4)\right]=
} & & \\
 & & \int_{\mathbb{R}^{4d}}dx_1\ dx_2\ dx_3\ dx_4\ \langle \phi(x_1)\phi(x_2)\phi(x_3)\phi(x_4)\rangle
\ f_1(x_1)f_2(x_2)f_3(x_3)f_4(x_4)
\end{eqnarray*}
featuring the pointwise correlation
\begin{eqnarray*}
\lefteqn{
\langle \phi(x_1)\phi(x_2)\phi(x_3)\phi(x_4)\rangle=
} & & \\
 & & 
\frac{\varkappa^2}{|x_1-x_2|^{2[\phi]}|x_3-x_4|^{2[\phi]}}
+\frac{\varkappa^2}{|x_1-x_3|^{2[\phi]}|x_2-x_4|^{2[\phi]}}
+\frac{\varkappa^2}{|x_1-x_4|^{2[\phi]}|x_2-x_3|^{2[\phi]}}\ .
\end{eqnarray*}
Note that, in this discussion, our pointwise correlations are seen as {\em ordinary functions} on the open subset of $\mathbb{R}^{nd}$
where the points $x_1,\ldots,x_n\in\mathbb{R}^d$ are {\em distinct}. The integrals above are also on this open subset.
Define (again at non-coinciding points) the new function
\[
\langle\ :\phi^2:(x_1)\ \phi(x_2)\phi(x_3)\rangle=\frac{2\varkappa^2}{|x_1-x_2|^{2[\phi]}|x_1-x_3|^{2[\phi]}}\ .
\]
Then one has the asymptotic behavior
\[
\langle \phi(x_1)\phi(x_2)\phi(x_3)\phi(x_4)\rangle=\langle \phi(x_1)\phi(x_2)\rangle
\langle \phi(x_3)\phi(x_4)\rangle+\langle\ :\phi^2:(x_2)\ \phi(x_3)\phi(x_4)\rangle+o\left(1\right)
\]
when $x_1\rightarrow x_2$ while the three points $x_2$, $x_3$, and $x_4$ are {\em fixed}.
This is the simplest instance of Wilson's OPE which here would say that, ``inside correlations'', one has
\[
\phi(x_1)\phi(x_2)=\mathcal{C}_{\phi\phi}^{\bbone}(x_1,x_2)\times 1+ \mathcal{C}_{\phi\phi}^{\phi^2}(x_1,x_2) :\phi^2:(x_2)+o(1)
\]
\[
{\rm with}\ \ 
\mathcal{C}_{\phi\phi}^{\bbone}(x_1,x_2)=\frac{\varkappa}{|x_1-x_2|^{2[\phi]}}
\ \ 
{\rm and}
\ \ 
\mathcal{C}_{\phi\phi}^{\phi^2}(x_1,x_2) =1\ .
\]
Our theorem shows how such an OPE, {\em with precise bounds on the remainder},
allows one to establish convergence in $L^p$ and almost surely
for suitably renormalized products
as in (\ref{WickSquare}). We do this in a vast setting which can handle Gaussian and non-Gaussian measures, massive and massless fields,
anomalous scaling dimensions, logarithmic corrections, finite degeneracy in the dimension spectrum,
as well as lack of translation invariance.
Much notation and machinery needs to be introduced before stating our theorem precisely in \S\ref{resultssec}.
This is provided in the following sections.
The general introduction will then resume in \S\ref{introcontd} which includes a discussion of related work.

\subsection{Abstract systems of pointwise correlations}
Let $\mathcal{A}$ be a finite set which we will call an alphabet and will serve to label fields.
Define the big diagonal
\[
{\rm Diag}_n=\{(x_1,\ldots,x_n)\in(\mathbb{R}^d)^n|\ \exists i\neq j, x_i=x_j\}
\]
and the configuration space
${\rm Conf}_n=(\mathbb{R}^d)^n\backslash {\rm Diag}_n$.
An {\em abstract system of pointwise correlations} consists in specifying
for all $n\ge 0$ and $A_1,\ldots,A_n\in\mathcal{A}$ an element of $C^{\infty}({\rm Conf}_n)$
denoted by
$\langle \mathcal{O}_{A_1}(x_1)\cdots \mathcal{O}_{A_n}(x_n)
\rangle$.
This is a purely symbolic notation.
No constituent of the formula has meaning by itself and the whole simply is some smooth function of the tuple
$(x_1,\ldots,x_n)$.
We impose that the degenerate $n=0$ case is taken care of by setting $\langle\emptyset\rangle=1$.
We assume symmetry, namely,
\begin{equation}
\langle\mathcal{O}_{A_{\sigma(1)}}(x_{\sigma(1)})\cdots \mathcal{O}_{A_{\sigma(n)}}(x_{\sigma(n)})\rangle=
\langle\mathcal{O}_{A_1}(x_1)\cdots \mathcal{O}_{A_n}(x_n)\rangle
\label{symmetry}
\end{equation}
for all permutations $\sigma\in\mathfrak{S}_n$.
We assume $\mathcal{A}$ contains a distinguished element $\bbone$ having the following ``forgetful'' property.
For all $n\ge 0$, and $A_1,\ldots,A_n$ in $\mathcal{A}$ we have
\begin{equation}
\langle\mathcal{O}_{\bbone}(z)
\mathcal{O}_{A_1}(x_1)\cdots\mathcal{O}_{A_n}(x_n)\rangle=
\langle\mathcal{O}_{A_1}(x_1)\cdots\mathcal{O}_{A_n}(x_n)\rangle
\label{forgetful}
\end{equation}
for all $(z,x_1,\ldots,x_n)\in {\rm Conf}_{n+1}$.

\subsection{Multilinear multilocal enhancement}\label{multiloc}

For $n\ge 0$, let $\mathcal{V}_n$ be the free module over the algebra $C^{\infty}({\rm Conf}_n)$
with basis $\mathcal{A}^n$.
A basis element indexed by $(A_1,\ldots,A_n)\in \mathcal{A}^n$ will be symbolically denoted
by $\mathcal{O}_{A_1}\otimes\cdots\otimes\mathcal{O}_{A_n}$.
Each $P\in \mathcal{V}_n$ has a unique expression as
\[
P=\sum_{(A_1,\ldots,A_n)\in \mathcal{A}^n} f_{A_1,\ldots,A_n}\ 
\mathcal{O}_{A_1}\otimes\cdots\otimes\mathcal{O}_{A_n}
\]
where the $f_{A_1,\ldots,A_n}$ are in $C^{\infty}({\rm Conf}_n)$.
To any $P\in \mathcal{V}_n$ we associate a function $\langle P\rangle\in C^{\infty}({\rm Conf}_n)$. Namely, it is the function
$(x_1,\ldots,x_n)\mapsto  \langle
P(x_1,\ldots,x_n)\rangle$ 
where, by definition,
\[
\langle
P(x_1,\ldots,x_n)\rangle=\sum_{(A_1,\ldots,A_n)\in \mathcal{A}^n} f_{A_1,\ldots,A_n}(x_1,\ldots,x_n) 
\ \langle\mathcal{O}_{A_1}(x_1)\cdots\mathcal{O}_{A_n}(x_n)
\rangle\ .
\]
The latter is not necessarily symmetric in the arguments $x_i$, $1\le i\le n$.
Note that one can give meaning to $P(x_1,\ldots,x_n)$ as a $C^{\infty}$ function ${\rm Conf}_n\rightarrow\mathbb{R}^{\mathcal{A}^n}$
but we will not use this point of view.
For $P\in\mathcal{V}_m$ and $Q\in\mathcal{V}_n$
we define their concatenation $P\otimes Q\in\mathcal{V}_{m+n}$
by
\[
P\otimes Q=\sum_{(A_1,\ldots,A_{m+n})\in \mathcal{A}^{m+n}} (f_{A_1,\ldots,A_m}\otimes g_{A_{m+1},\ldots,A_{m+n}}) 
\ \mathcal{O}_{A_1}\otimes\cdots\otimes\mathcal{O}_{A_{m+n}}
\]
\begin{flalign*}
{\rm if\ \ \ }& &
P& =\sum_{(A_1,\ldots,A_m)\in \mathcal{A}^m} f_{A_1,\ldots,A_m} 
\ \mathcal{O}_{A_1}\otimes\cdots\otimes\mathcal{O}_{A_m} & &\ \\
{\rm and}& &
Q& =\sum_{(B_1,\ldots,B_n)\in \mathcal{A}^n} g_{B_1,\ldots,B_n} 
\ \mathcal{O}_{B_1}\otimes\cdots\otimes\mathcal{O}_{B_n}  & &\ 
\end{flalign*}
where we used the notation $(f\otimes g)(x_1,\ldots,x_{m+n})=f(x_1,\ldots,x_m)g(x_{m+1},\ldots,x_{m+n})$ for
$f\in C^{\infty}({\rm Conf}_m)$ and $g\in C^{\infty}({\rm Conf}_n)$.
We will use the notation
\begin{equation}
\langle P(x_1,\ldots,x_m)Q(x_{m+1},\ldots,x_{m+n})\rangle=
\langle (P\otimes Q)(x_1,\ldots,x_{m+n}) \rangle
\label{braidingnotation}
\end{equation}
for the evaluation of $\langle P\otimes Q\rangle\in C^{\infty}({\rm Conf}_{m+n})$ on the argument $(x_1,\ldots,x_{m+n})$,
and similarly for higher products $\langle P_1\otimes\cdots\otimes P_N\rangle$.
The latter are unambiguously defined since concatenation is associative.
However, concatenation is not commutative: 
in general $P\otimes Q\neq Q\otimes P$ and also the smooth functions of the argument $(x_1,\ldots,x_{m+n})$
given by $\langle P\otimes Q\rangle$ and $\langle Q\otimes P\rangle$ are different.
Because of the symmetry (\ref{symmetry}), one has instead
$\langle P\otimes Q\rangle=\langle Q\otimes P\rangle \circ\tau$
where $\tau$ is the coordinate permutation map given by
\[
\tau(x_1,\ldots,x_m,x_{m+1},\ldots,x_{m+n})=(x_{m+1},\ldots,x_{m+n},x_1,\ldots,x_m)\ .
\]
Nevertheless, 
if for instance
$P,Q,R$ belong to $\mathcal{V}_2$
then
\[
\langle P(x_1,y_1)Q(x_2,y_2)R(x_3,y_3)\rangle=\langle
R(x_3,y_3)Q(x_2,y_2)P(x_1,y_1)\rangle 
\]
and similarly for other permutations in $\mathfrak{S}_3$.
This is the point of the notation (\ref{braidingnotation}): if $P_1,\ldots,P_N$ are say in $\mathcal{V}_2$,
$\langle P_1(x_1,y_1)\cdots P_N(x_N,y_N)\rangle$
is an unambiguously defined function of the $x$ and $y$ arguments.
There is no need to indicate an order in which the product must be written.
Finally, note that there is a clear notion of subsystem of pointwise correlations corresponding to a subset $\mathcal{A}'$
of the alphabet $\mathcal{A}$ which contains $\bbone$. It is obtained by only keeping data concerning labels in $\mathcal{A}'$.

\subsection{OPE structure}

An OPE structure consists of an abstract system of pointwise correlations together with some extra data.
To $A\in\mathcal{A}$ we associate a number $[A]\in[0,\infty)$ called the scaling dimension of the field labeled by $A$.
We impose $[\bbone]=0$. For $\Delta\in[0,\infty)$, we let $\mathcal{A}(\Delta)=\{A\in\mathcal{A}\ |\ [A]\le \Delta\}$. 
We assume that for each triple $(A,B,C)\in\mathcal{A}^3$
we have an element $\mathcal{C}_{AB}^{C}$ in $C^{\infty}({\rm Conf}_2)$ thus giving rise to a smooth function
$\mathcal{C}_{AB}^{C}(x,y)$ 
of $(x,y)\in{\rm Conf}_2$.
An OPE structure is what is needed in order to formulate physicists's statements that if one fixes a number $\Delta$
then the ``operators'' $\mathcal{O}_A$, $A\in\mathcal{A}$ , satisfy an operator product expansion
\begin{equation}
\mathcal{O}_{A}(x)\mathcal{O}_{B}(y)=\sum_{C\in\mathcal{A}(\Delta)} \mathcal{C}_{AB}^{C}(x,y)\mathcal{O}_{C}(y)
+o(|x-y|^{\Delta-[A]-[B]})
\label{opelittleo}
\end{equation}
``inside correlations''.
A precise statement corresponding to such intuition is that, for all $n\ge 0$ and $D_1,\ldots,D_n$ in
$\mathcal{A}$ and all fixed $(y,z_1,\ldots,z_n)\in{\rm Conf}_{n+1}$,
we have
\[
\langle\mathcal{O}_{A}(x)\mathcal{O}_{B}(y)\mathcal{O}_{D_1}(z_1)\cdots \mathcal{O}_{D_n}(z_n)\rangle=
\]
\[
\sum_{C\in\mathcal{A}(\Delta)} \mathcal{C}_{AB}^{C}(x,y)\langle
\mathcal{O}_{C}(y)\mathcal{O}_{D_1}(z_1)\cdots \mathcal{O}_{D_n}(z_n)
\rangle+o(|x-y|^{\Delta-[A]-[B]})
\]
when taking the limit $x\rightarrow y$.
The very easy example from \S\ref{introtointro}
corresponds to $\mathcal{A}=\{\bbone,\phi,\phi^2\}$ and $\Delta=[\phi^2]=2[\phi]$.
What makes the above statement nontrivial
is that the shape of this asymptotic expansion is independent
of the number $n$ and labels $D_1,\ldots,D_n$ as well as positions $z_1,\ldots,z_n$
for the ``spectator fields'' $\mathcal{O}_{D_1}(z_1),\ldots,\mathcal{O}_{D_n}(z_n)$.

The multilinear multilocal enhancement of \S\ref{multiloc} allows a more elegant rephrasing as
\[
\langle {\rm OPE}(x,y)\ \mathcal{O}_{D_1}(z_1)\cdots \mathcal{O}_{D_n}(z_n)
\rangle=o(|x-y|^{\Delta-[A]-[B]})
\]
where ${\rm OPE}\in \mathcal{V}_2$ is defined by
\begin{equation}
{\rm OPE}=\mathcal{O}_A\otimes\mathcal{O}_{B}-\sum_{C\in\mathcal{A}(\Delta)} \mathcal{C}_{AB}^{C}
\ \mathcal{O}_{\bbone}\otimes\mathcal{O}_C\ .
\label{opelike}
\end{equation}

Note that a subset $\mathcal{A}'\subset\mathcal{A}$ which contains $\bbone$ also defines a sub-OPE structure where one restricts the $[A]$'s and $\mathcal{C}_{AB}^{C}$'s to labels in $\mathcal{A}'$.

\subsection{Probabilistic incarnations}\label{probareal}

Suppose that our system of pointwise correlations is such that
all smooth functions
$\langle\mathcal{O}_{A_1}(x_1)\cdots\mathcal{O}_{A_n}(x_n)
\rangle$
are locally integrable on ${\rm Diag}_n$ and have suitable moderate growth at infinity
so the integrals
\[
\int_{{\rm Conf}_n}\prod_{i=1}^{n}dx_i\ 
f(x_1,\ldots,x_n)\ \langle\mathcal{O}_{A_1}(x_1)\cdots\mathcal{O}_{A_n}(x_n)
\rangle
\]
converge absolutely for all test functions $f\in\mathcal{S}(\mathbb{R}^{nd})$.
Then it makes sense to talk about what we call a {\em probabilistic incarnation}.
It is given by a (not necessarilty complete) probability space $(\Omega,\mathcal{F},\mathbb{P})$
together with $\mathcal{S}'(\mathbb{R}^d)$-valued random variables $\mathcal{O}_A$, $A\in\mathcal{A}$
such that for all $f\in \mathcal{S}(\mathbb{R}^{d})$, the real-valued random variable $\mathcal{O}_{A}(f)$
has moments of all orders and for all $n\ge 0$, $A_1,\ldots,A_n\in\mathcal{A}$ and
all $f_1,\ldots,f_n\in \mathcal{S}(\mathbb{R}^{d})$ one has
\[
\mathbb{E}\left[\mathcal{O}_{A_1}(f_1)\cdots\mathcal{O}_{A_n}(f_n)
\right]=\int_{{\rm Conf}_n}\prod_{i=1}^{n}dx_i\ 
f_1(x_1)\cdots f_n(x_n)\ \langle\mathcal{O}_{A_1}(x_1)\cdots\mathcal{O}_{A_n}(x_n)
\rangle\ .
\]
We used the same notation $\mathcal{O}_A$ on both sides of the equation
but there is no risk of confusion. If $\mathcal{O}_A$ appears in a formula with the expectation symbol $\mathbb{E}$ then
we mean an honest distribution-valued random variable from the probabilistic incarnation.
If not, then we mean a constituent symbol participating in the definition of a pointwise correlation.
Also, when we say $\mathcal{O}_A$ is an $\mathcal{S}'(\mathbb{R}^d)$-valued random variable, we mean that the map $\Omega\rightarrow
\mathcal{S}'(\mathbb{R}^d)$
given by $\omega\mapsto \mathcal{O}_A(\omega)$
is $(\mathcal{F},{\rm Borel}(\mathcal{S}'(\mathbb{R}^d)))$-measurable.
As said earlier, the space of temperate distributions $\mathcal{S}'(\mathbb{R}^d)$ is here
equipped with the strong topology, although the weak-$\ast$
topology would give rise to the same Borel $\sigma$-algebra.

\subsection{Soft hypotheses: kernel semi-regularity}\label{softhyp}

Regarding the functions $\mathcal{C}_{AB}^{C}$ we make the following rather mild hypotheses.
We assume that for all $z\in\mathbb{R}^{d}$ the function $y\mapsto \mathcal{C}_{AB}^{C}(y,z)$ is locally integrable at $y=z$
and has at most polynomial growth at infinity. Namely, we assume local integrability together with
$\forall z\in\mathbb{R}^{d}$, $\exists k\in\mathbb{N}$,
$\exists K>0$, $\exists R>0$ such that for all $y\in\mathbb{R}^{d}\backslash\{z\}$,
\begin{equation}
|y|>R\ \Longrightarrow\ |\mathcal{C}_{AB}^{C}(y,z)|\le K\langle y\rangle^k\ .
\label{softbound}
\end{equation}
As a result, and following Schwartz's subscript notation, $\mathcal{C}_{AB}^{C}(y,z)$ can be seen as a $z$-dependent element
of $\mathcal{S}'_y(\mathbb{R}^d)$ defined by
\[
\langle \mathcal{C}_{AB}^{C}(y,z), f(y)\rangle_y=\int_{\mathbb{R}^{d}\backslash\{z\}}dy\ \mathcal{C}_{AB}^{C}(y,z) f(y)
\]
for all test function $f(y)$ in $\mathcal{S}_y(\mathbb{R}^d)$.
We now make two more assumptions or ``soft hypotheses''.

\noindent{\bf SH1:}
For all $g(z)$ in $\mathcal{S}_z(\mathbb{R}^d)$ and all $f(y)$ in $\mathcal{S}_y(\mathbb{R}^d)$,
the function 
$g(z)\langle \mathcal{C}_{AB}^{C}(y,z), f(y)\rangle_y$
belongs to $\mathcal{S}_z(\mathbb{R}^d)$.
Moreover, the resulting map 
$f(y)\mapsto g(z)\langle \mathcal{C}_{AB}^{C}(y,z), f(y)\rangle_y$
is a continuous map $\mathcal{S}_y(\mathbb{R}^d)\rightarrow \mathcal{S}_z(\mathbb{R}^d)$.

\noindent{\bf SH2:}
For all functions $g(x,z)$ in $\mathcal{S}_{x,z}(\mathbb{R}^{2d})$ and $f(y,w)$ in $\mathcal{S}_{y,w}(\mathbb{R}^{2d})$,
the function defined by
$g(x,z)\langle \mathcal{C}_{AB}^{C}(y,z), f(y,w)\rangle_y$
is in $\mathcal{S}_{x,z,w}(\mathbb{R}^{3d})$
and the map which sends 
$f(y,w)$ to the function $g(x,z)\langle \mathcal{C}_{AB}^{C}(y,z), f(y,w)\rangle_y$
is a continuous map $\mathcal{S}_{y,w}(\mathbb{R}^{2d})\rightarrow \mathcal{S}_{x,z,w}(\mathbb{R}^{3d})$.

Note that in fact (SH1) implies (SH2), as a consequence of the Schwartz theory for Volterra composition~\cite[\S4]{SchwartzVect1}.
Essentially, we are assuming that the kernel $\mathcal{C}_{AB}^{C}(y,z)$ is semi-regular and of moderate growth in $z$.
This is a pedestrian way of saying
$\mathcal{C}_{AB}^{C}\in \mathcal{S}'_{y}(\mathbb{R}^d)\ \widehat{\otimes}\ 
\mathcal{O}_{{\rm M},z}(\mathbb{R}^d)$
in the sense of~\cite[\S4]{SchwartzVect1}.
Note that there is no need~\cite[Theorem 50.1 (f)]{Treves} to specify the tensor product
since all spaces involved are nuclear~\cite[Ch. 2, Theorem 10, p. 55]{Grothendieck}. 
It is easy to check (SH1) and (SH2) by hand in the case of the FMFF
from \S\ref{introtointro}, as will be done in \S\ref{softhardsec}.
The functional analytic difficulties in the more general framework of this section are due to {\em non-translation invariance}
(see, e.g.,~\cite{Dolezal,Meidan,Zemanian}
where similar issues arise).

\subsection{Hard hypotheses: factorized nearest neighbor bounds}\label{hardhyp}

From now on we assume that $[A]\in [0,\frac{d}{2})$, for all $A\in\mathcal{A}$.
We sharpen the bound (\ref{softbound}) on the OPE ``structure constants'' $\mathcal{C}_{AB}^{C}$
by requiring $\forall A,B,C\in\mathcal{A}$, $\forall \epsilon>0$, $\exists k\in\mathbb{N}$, $\exists K>0$, $\forall (x,y)\in{\rm Conf}_2$,
\begin{equation}
|\mathcal{C}_{AB}^{C}(x,y)|\le \frac{K}{|x-y|^{[A]+[B]-[C]+\epsilon}} \langle x\rangle^k\langle y\rangle^k\ .
\label{hardCbd}
\end{equation}
An element in $\mathcal{V}_2$ is called {\em OPE-like} if it is of the form
(\ref{opelike}) for some $A,B\in\mathcal{A}$ and $\Delta\in\{[C]\ |\ C\in\mathcal{A}\}$.
An element in $\mathcal{V}_2$ is called {\em CZ-like} if it is of the form
$\mathcal{O}_{A}\otimes\mathcal{O}_{\bbone}-
\mathcal{O}_{\bbone}\otimes\mathcal{O}_{A}$
for some $A\in\mathcal{A}$. ``CZ'' stands for Calder\'on-Zygmund.

For all $m,n,p\ge 0$, for all OPE-like elements ${\rm OPE}_1,\ldots,{\rm OPE}_m$, 
given by
\[
{\rm OPE}_i=\mathcal{O}_{A_i}\otimes\mathcal{O}_{B_i}-\sum_{C_i\in\mathcal{A}(\Delta_i)}
\mathcal{C}_{A_i B_i}^{C_i}\ \mathcal{O}_{\bbone}
\otimes\mathcal{O}_{C_i}
\]
and all CZ-like elements
${\rm CZ}_{m+1},\ldots,{\rm CZ}_{m+n}$ given by
${\rm CZ}_i=\mathcal{O}_{B_i}\otimes\mathcal{O}_{\bbone}-
\mathcal{O}_{\bbone}\otimes\mathcal{O}_{B_i}$
with arbitrary $B_{m+1},\ldots,B_{m+n+p}\in\mathcal{A}$,
we require:
$\exists \eta>0$, $\exists \gamma>0$, $\forall \epsilon>0$, 
$\exists k\in\mathbb{N}$, $\exists K>0$,
\[
\prod_{i=1}^{m+n}\bbone\left\{
|y_i-x_i|\le\eta\min_{j\neq i}|x_i-x_j|
\right\}\times
\left|
\left\langle
\prod_{i=1}^{m} {\rm OPE}_i(y_i,x_i)
\prod_{i=m+1}^{m+n} {\rm CZ}_i(y_i,x_i)
\prod_{i=m+n+1}^{m+n+p} \mathcal{O}_{B_i}(x_i)
\right\rangle
\right|\le
\]
\[
K\prod_{i=1}^{m+n+p}\langle x_i\rangle^k\times\prod_{i=1}^{m+n}\langle y_i\rangle^k
\times
\prod_{i=1}^{m}\left\{|y_i-x_i|^{\Delta_i+\gamma-[A_i]-[B_i]}\times\left(\min_{j\neq i}|x_i-x_j|\right)^{-\Delta_i-\gamma-\epsilon}\right\}
\]
\begin{equation}
\times
\prod_{i=m+1}^{m+n}\left\{|y_i-x_i|^{\gamma}\times\left(\min_{j\neq i}|x_i-x_j|\right)^{-[B_i]-\gamma-\epsilon}\right\}
\times
\prod_{i=m+n+1}^{m+n+p}\left(\min_{j\neq i}|x_i-x_j|\right)^{-[B_i]-\epsilon}
\label{efnnb}
\end{equation}
for all collections of $2m+2n+p$ distinct points in $\mathbb{R}^d$.
When, for fixed $i$, we write $\min_{j\neq i}$,
we mean the minimum over $j$ such that $1\le j\le m+n+p$ and $j\neq i$.

We call (\ref{efnnb}) the {\em enhanced factorized nearest neighbor bound} (EFNNB).
It includes, as the $m=n=0$ special case, the {\em basic factorized nearest neighbor bound} (BFNNB) for pointwise correlations:
$\forall \epsilon>0$, $\exists k\in\mathbb{N}$, $\exists K>0$,
\begin{equation}
\left|
\left\langle
\mathcal{O}_{B_1}(x_1)\cdots\mathcal{O}_{B_p}(x_p)
\right\rangle
\right|\le
K\prod_{i=1}^{p}\langle x_i\rangle^k\times\prod_{i=1}^{p}\left(\min_{j\neq i}|x_i-x_j|\right)^{-[B_i]-\epsilon}
\label{bfnnb}
\end{equation}
for all $(x_1,\ldots,x_p)\in{\rm Conf}_p$. The $\epsilon$ is only needed in order to account for eventual logarithmic corrections.
In a conformal field theory (CFT) as discussed in \cite{Abdesselam3DCP},
one can take $\epsilon=0$. The $k$ allows more generality for our main theorem and is natural in the
setting of probability theory on spaces of temperate distributions,
but it is not needed in usual QFT models.
In order to work in the space $\mathcal{D}'(\mathbb{R}^d)$ of general Schwartz distributions, the previous bounds can be readily adapted
by dropping $k$ altogether and letting the constant $K$ depend on the radius of a large ball in which all the points must be confined.
 
For further reference, we introduce the following terminology regarding the EFNNB.
The points $x_i$ or rather their {\em labels} $i$, $1\le i\le m+n+p$, are called {\em effective} because they can be someone else's nearest
neighbor and they all participate in the computation of the minimums over distance.
By contrast, the points $y_i$ or more precisely their labels $i$, $1\le i\le m+n$, are called {\em virtual}.
They only communicate
with their $x_i$ which serves as a local point of reference, the center of their ``mini-universe'' or ``solar system''. 

\subsection{Worst-case scenario planning}\label{worstsec}

For the bound on $\mathcal{C}_{AB}^{C}$,
the first quantifier is ``$\forall$'' regarding the format, i.e., the choice of triple $(A,B,C)$.
If one anticipates needing the bound for several but {\em finitely many} formats, it is advantageous to have the same $K,k$ for all formats.
Namely, the order of quantifiers can be modified to $\forall\epsilon>0$, $\exists k\in\mathbb{N}$, $\exists K>0$, $\forall$ format, etc.
One just needs to take the largest $k$ and $K$.

A similar property holds for the EFNNB. In this case a format consists in a choice
of $m,n,p$, the $A_i$, $B_i$, $\Delta_i$ featuring in the OPE-like elements, and the choice of $B_i$'s for $m+1\le i\le m+n+p$.
Let us denote such a format by $\mathsf{F}$.
First pick the smallest $\eta$ and $\gamma$, namely, set $\eta=\min_{\mathsf{F}}\eta_{\mathsf{F}}$ and $\gamma=
\min_{\mathsf{F}}\gamma_{\mathsf{F}}$.
Then if one chooses $\epsilon>0$, the indicator function on the left-hand side of (\ref{efnnb})
with the new $\eta$ is bounded by the similar one for $\mathsf{F}$. The correlator defined by the format $\mathsf{F}$
on the left-hand side is bounded as before and the old majorant with $\gamma_{\mathsf{F}}$ is converted into the new one
with the uniform $\gamma$,
at the cost of creating an additional factor
\[
\prod_{i=1}^{m+n}\left\{
\frac{|y_i-x_i|}{\min_{j\neq i}|x_i-x_j|}
\right\}^{\gamma_{\mathsf{F}}-\gamma}
\le \eta^{(m+n)(\gamma_{\mathsf{F}}-\gamma)}\ .
\]
Finally, take $k=\max_{\mathsf{F}} k_{\mathsf{F}}$ and
$K=\max_{\mathsf{F}} [K_{\mathsf{F}} \eta^{(m+n)(\gamma_{\mathsf{F}}-\gamma)}]$
while keeping in mind that $m,n$ depend on the format $\mathsf{F}$.
As a result, the bound (\ref{efnnb}) also works with the new order of quantifiers
$\exists \eta>0$, $\exists \gamma>0$, $\forall \epsilon>0$, 
$\exists k\in\mathbb{N}$, $\exists K>0$, $\forall\mathsf{F}$.
By the same reasoning, one can restrict the common $\gamma$
to be less than a specified positive number if needed, e.g., for local integrability reasons.

\begin{Remark}
In \S\ref{probatocomb} and \S\ref{mainestimatesec} we will be rather implicit on how we do the worst-case scenario planning.
We suggest to the reader to go over these sections a first time without worrying about how small $\gamma$ and $\epsilon$ have to be
and then come back to this section in order to see that the use of the same $k$ and  $K$ is legitimate. The reason for this is that
the number of formats involved in the proofs is finite, because $\mathcal{A}$ is finite and the correlation in Proposition \ref{combiprop}
has fixed order.
\end{Remark}
\subsection{Non-degeneracy condition}

In order to use the $\mathcal{O}_{A}\times\mathcal{O}_{B}$ OPE for the definition of a field $\mathcal{O}_{C}$,
we need to be able to peel it off by a formula one may write intuitively as
\[
\mathcal{O}_{C}\sim \frac{1}{\mathcal{C}_{AB}^{C}}\left[\mathcal{O}_{A}\times\mathcal{O}_{B}
-\sum_{D\neq C} \mathcal{C}_{AB}^{D}\mathcal{O}_{D}\right]
\]
and this requires a non-vanishing or non-degeneracy condition on the OPE coefficient $\mathcal{C}_{AB}^{C}$.

We say that the triple $(A,B,C)$ is {\em non-degenerate} if
$\forall \epsilon>0$, $\exists k\in\mathbb{N}$, $\exists K>0$, $\forall (x,y)\in{\rm Conf}_2$,
\begin{equation}
\mathcal{C}_{AB}^{C}(x,y)\ge \frac{1}{K|x-y|^{[A]+[B]-[C]-\epsilon}}\times \langle x\rangle^{-k} \langle y\rangle^{-k}\ . 
\label{nondeg}
\end{equation}

\subsection{Statement of results}\label{resultssec}

We assume that we have an OPE structure with finite alphabet $\mathcal{A}$ and such that $[A]\in \left[0,\frac{d}{2}\right)$
for all $A\in\mathcal{A}$. We assume the hypotheses from \S\ref{softhyp} and \S\ref{hardhyp} hold for this OPE structure.
We suppose that we have a subset $\mathcal{B}\subset \mathcal{A}$ together with a probabilistic incarnation $(\mathcal{O}_B)_{B\in\mathcal{B}}$
on some probability space $(\Omega,\mathcal{F},\mathbb{P})$
for the subsystem of pointwise correlations corresponding to $\mathcal{B}$.
We assume that $C_{\ast}$
is an element of $\mathcal{A}\backslash\mathcal{B}$ such that
$\mathcal{A}([C_{\ast}])\backslash\{C_{\ast}\}\subset\mathcal{B}$.
We also assume that we have two elements $A$, $B$ in $\mathcal{B}$ such that the triple $(A,B,C_{\ast})$ is non-degenerate.
We let $\rho$ be a mollifier as in \S\ref{introtointro}, except that we now add the more restrictive hypotheses
that $\rho$ is compactly supported with ${\rm supp}\ \rho\in \bar{B}(0,1)$ (the closed Euclidean ball of radius one around the origin),
and is pointwise nonnegative.
We again use the notation $\rho_r(x)=L^{-rd}\rho(L^{-r}x)$ for the rescaled mollifier and use this to define the function
\[
Z_{r}(x)=\left\{
\int_{{\rm Conf}_2}dy\ dz\ \rho_{r}(x-y)\rho_r(x-z)\ \mathcal{C}_{AB}^{C_{\ast}}(y,z)
\right\}^{-1}
\]
in $\mathcal{O}_{{\rm M},x}(\mathbb{R}^d)$.
We then introduce the random element of $\mathcal{O}_{{\rm M},x}(\mathbb{R}^d)$
given by
\begin{equation}
M_{r}(x)=Z_r(x)\left[
\mathcal{O}_{A,r}(x)\mathcal{O}_{B,r}(x)-\sum_{C\in\mathcal{A}(\Delta)\backslash\{C_{\ast}\}} \widetilde{\mathcal{O}}_{C,r}(x)
\right]
\label{candidate}
\end{equation}
\begin{flalign*}
&{\rm where}\ \ & &
\mathcal{O}_{A,r}(x) = (\mathcal{O}_A\ast\rho_r)(x)=\langle \mathcal{O}_{A}(y),\rho_r(x-y) \rangle_y & & \ \\
&{\rm and\ similarly}\ \ & & 
\mathcal{O}_{B,r}(x) = (\mathcal{O}_B\ast\rho_r)(x)=\langle \mathcal{O}_{B}(z),\rho_r(x-z) \rangle_z & & \ \\
&{\rm while}\ \ & &
\widetilde{\mathcal{O}}_{C,r}(x) = \langle \mathcal{O}_{C}(z),g_r(x,z)\rangle_{z} & & \ \\
&{\rm with}\ \ & &
g_{r}(x,z) = \rho_r(x-z)\times\int_{\mathbb{R}^d\backslash\{z\}}dy\ \rho_r(x-y)\ \mathcal{C}_{AB}^{C}(y,z)\ . & & \ 
\end{flalign*}
Note that the dependence on the sample $\omega\in\Omega$ has been suppressed from the notation
and that $\mathcal{O}_A$, $\mathcal{O}_B$, $\mathcal{O}_C$
designate the distribution-valued random variables provided by the probabilistic incarnation for $\mathcal{B}$.
We view $M_r(x)$ as the random Schwartz distribution whose action on a test function
$f\in\mathcal{S}(\mathbb{R}^d)$ is of course given by
\[
M_r(f)=\int_{\mathbb{R}^d}dx\ M_{r}(x)f(x)=\langle  M_r(x), f(x)\rangle_x\ .
\]
It is not trivial to show that $M_r(f)$ is indeed well defined, $\mathcal{F}$-measurable,
and in every $L^p(\Omega,\mathcal{F},\mathbb{P})$,
$p\ge 1$.
This will be done in \S\ref{probatocomb}. The main result in this article is as follows.

\begin{Theorem}\label{maintheorem}\ 

\begin{enumerate}
\item
For any test function $f$, and when taking $r\rightarrow-\infty$, the random variable $M_{r}(f)$ converges
in every $L^p(\Omega,\mathcal{F},\mathbb{P})$, $p\ge 1$,
and $\mathbb{P}$-almost surely to a random variable which we will denote by $\mathcal{O}_{C_{\ast}}(f)$.
\item
The limit is independent from the choice of mollifier $\rho$.
\item
There exists a Borel-measurable map
\[
\mathcal{P}:\prod_{C\in\mathcal{B}}\mathcal{S}'(\mathbb{R}^d)\rightarrow \mathcal{S}'(\mathbb{R}^d)
\]
such that for all $f\in\mathcal{S}(\mathbb{R}^d)$,
$\mathcal{O}_{C_{\ast}}(f)=\left[\mathcal{P}\left((\mathcal{O}_C)_{C\in\mathcal{B}}\right)\right](f)$,
$\mathbb{P}$-almost surely.
\item
If one extends the probabilistic incarnation to $\mathcal{B}\cup\{C_{\ast}\}$
by adding the $\mathcal{S}'(\mathbb{R}^d)$-valued random
variable $\mathcal{P}\left((\mathcal{O}_C)_{C\in\mathcal{B}}\right)$, then the result is a probabilistic incarnation
of the system of pointwise correlations corresponding to the new set of labels $\mathcal{B}\cup\{C_{\ast}\}$.
\end{enumerate}
\end{Theorem}

Clearly, this can be iterated. By growing $\mathcal{B}$ and also using the trivial construction of 
derivatives in the sense of Schwartz distributions of already existing random fields, one can construct all the
composite operators of the fractional massless free field from \S\ref{introtointro} in this way,
provided one remains under the Calder\'on-Zygmund $\frac{d}{2}$
threshold for scaling dimensions. This will be explained in detail in \S\ref{examplesec}
which provides an example where all the hypotheses of our theorem are satisfied.
Also note that if one constructs say $\phi^4$ as $\phi^2\times \phi^2$ or as $\phi\times \phi^3$, the result is the same.
Namely, our renormalized product construction is {\em associative}, as an easy consequence of Theorem \ref{maintheorem}, Part (4).
Indeed, by a the same $L^2$ distance computation as in \S\ref{finishproof} below, a priori different probabilistic incarnations of the same symbolic field $\mathcal{O}_{C\ast}$ must be equal $\mathbb{P}$-almost surely.

\subsection{Wider context and structure of the article}\label{introcontd}

\subsubsection{The OPE in mathematics and physics}
There exist several versions of the OPE. The one used here is the pointwise OPE in position or $x$-space.
One can also express the OPE in a smeared sense, i.e., asymptotics such as (\ref{opelittleo})
are to be interpreted in the sense of distributions involving suitable test functions~\cite[\S2]{WilsonZ}.
Finally, there is the OPE in Fourier or momentum space when modifing two momenta by adding $\xi$ and $-\xi$ respectively
and then taking the limit $|\xi|\rightarrow\infty$ (see, e.g.,~\cite[\S10.4]{ZinnJustin}).
The latter was historically important since it provided the theoretical counterpart (see, e.g.,~\cite[Ch. 14]{Collins})
of the SLAC experiments which confirmed the quark picture~\cite{Breidenbach}.
Such deep inelastic scattering experiments
were an important milestone, since they helped guide researchers towards the discovery of
asymptotic freedom and the elaboration of the Standard Model of particle physics.
As to the pointwise OPE, it is a cornerstone of the conformal bootstrap (see, e.g.,~\cite{DiFrancescoMS} for the 2D situation and~\cite{Rychkov,SimmonsD} for higher dimensions).
Note that the pointwise OPE has led to important developments in mathematics too. The notions of vertex operator algebras
(see, e.g.,~\cite{FrenkelB}), as well as chiral and factorization algebras~\cite{BeilinsonD} can be seen as ways of capturing
the mathematical structure of the pointwise OPE. These mathematical applications pertain to algebraic geometry and
representation theory or, more specifically, the area known as the geometric Langlands program (see, e.g.,~\cite{Frenkel}).
This article shows that the pointwise OPE is also important for {\em probability theory}.
Our construction of the renormalized product $\mathcal{O}_{C_{\ast}}$ is the adaptation to the probability
context of what is known as a ``point-splitting procedure'' in the QFT context~\cite{Dirac,Heisenberg,HeisenbergE,Valatin}.
The theorem from \S\ref{resultssec} relies on a functional-analytic part done in \S\ref{probatocomb} using tools
from what one may call Schwartz-Grothendieck-Fernique theory in view of the foundational
works~\cite{SchwartzBook1,SchwartzBook2,SchwartzVect1,Grothendieck,Fernique}.
It also relies on Proposition \ref{combiprop} which is a combinatorial estimate
in the pure tradition of the \'Ecole Polytechnique school of constructive QFT
founded by Roland S\'en\'eor. This estimate allows one to go from the pointwise OPE to a smeared OPE.

\subsubsection{Where our bounds and methods came from}\label{motivsec}
Let us now explain the origin of the BFNNB and EFNNB. They both were outcomes of the study
of the $p$-adic or hierarchical fractional $\phi^4$ model in~\cite{AbdesselamCG1}.
The author obtained the BFNNB with $\mathcal{A}=\{\bbone,\phi,\phi^2\}$ for the non-Gaussian
hierarchical scaling limit constructed in~\cite{AbdesselamCG1}, together with the local integrability proof in
\S\ref{warmupsec} (see~\cite[Theorem 3]{AbdesselamSlides}). The results from~\cite{AbdesselamCG1} give an explicit bi-infinite series over scales representation
for the pointwise $\phi$ and $\phi^2$ correlations in terms of tree-like
compositions of certain maps, 
among which the most important is the one corresponding to degree two vertices forming
linear chains in the tree. This map, $\dot{V}\mapsto RG_{\rm dv}[\vec{V}_{\ast},\dot{V}]$
in the notation
of~\cite{AbdesselamCG1}, is the renormalization group (RG) acting in the space $\mathcal{E}_{\rm pt}$
of ``non-integrated'' or point-like operator perturbations of the infrared fixed point $\vec{V}_{\ast}$.
Its differential at $\dot{V}=0$, for a rescaling by a factor of $L$, has the eigenvalue $L^{-[A]}$
for the operator $\mathcal{O}_{A}$ instead of $L^{d-[A]}$, and is analogous
to $e^{-(\log L)D}$ where $D$ is the generator of dilations in CFT radial quantization (see, e.g.,~\cite[\S6]{SimmonsD}). The (fusion) tree is determined by geometry, i.e., the relative
positions of the points $x_1,\ldots,x_p$ for which the correlation is evaluated.
The BFNNB comes from the estimates in~\cite{AbdesselamCG1} applied to the initial linear chains steming from the leaves
of the tree until the {\em first fusion}
(vertex of the tree with degree $\ge 3$) whose scale is given by the distance
to the {\em nearest neighbor}.
The bound is in {\em factorized} form because the correlation is computed using a Fr\'echet differential
of order $p$ (the number of points) which is continuous, i.e., bounded by the {\em product}
of norms of the inputs.

The contribution of the OPE-like factors in (\ref{efnnb}) comes from similar yet more involved
RG arguments.
The results in the present article were first derived for random fields in $\mathcal{S}'(\mathbb{Q}_p^d)$ and then later adapted to random fields in $\mathcal{S}'(\mathbb{R}^d)$.
The EFNNB over $\mathbb{Q}_p$ looks exactly the same as the one over $\mathbb{R}$ given in (\ref{efnnb}),
except for the {\em effect} of the CZ-like factors.
The proof in the $p$-adic case {\em requires} CZ-like terms, i.e., ``moving legs around'' as in BPHZ
renormalization in $x$-space~\cite[\S{II.2}]{RivasseauBook}
or as in (\ref{Qexp}) below. However, the pointwise correlations outside the diagonal are smooth
which in the $p$-adic setting means {\em locally constant}: the CZ-like elements have a vanishing contribution.
In the real case, one can use the Fundamental Theorem of Calculus
$\mathcal{O}_{C}(y)-\mathcal{O}_{C}(x)=\sum_{|\alpha|=1}\int_0^1 dt\ (y-x)^{\alpha}\partial^{\alpha}\mathcal{O}_{C}(x+t(y-x))$
``inside correlations''.
In the particular CFT case~\cite{Abdesselam3DCP},
where a derivative increases the scaling dimension by one, 
we see that such a CZ-element has the contribution written on the right-hand side of (\ref{efnnb}) with $\gamma=1$.
In order to make Theorem \ref{maintheorem} more general, we also allowed $0<\gamma<1$.
It turns out that this kind of bound is exactly what is used in harmonic analysis and the theory of Calder\'on-Zygmund operators
(see~\cite[p. 372]{DavidJ} and~\cite[p. 9]{MeyerC}), hence our choice of terminology.
Note that our motivation for the present
work is that it accomplishes one of the tasks from the program outlined in~\cite{Abdesselam3DCP}:
we showed how~\cite[Conjecture 8]{Abdesselam3DCP} implies~\cite[Conjecture 9]{Abdesselam3DCP}.
In order to help the reader better understand the EFNNB (\ref{efnnb}), a simple heuristic is as follows.
Replace $CZ_i$ by $(y_i-x_i) \partial\mathcal{O}_{B_i}(x_i)$. Replace the OPE remainder $OPE_i(y_i,x_i)$
by $\mathcal{C}_{A_i,B_i}^{E}(y_i,x_i)\mathcal{O}_{E}(x_i)$ where $E$ is next operator appearing in the expansion, i.e., the one with the smallest $[E]>\Delta_i$. Finally, apply the BFNNB. The result essentially is the EFNNB.
The upper gap $[E]-\Delta_i$ in the spectrum of scaling dimensions is what limits how large $\gamma>0$ can be. In fact,
this heuristic is turned into a rigorous proof of the EFNNB in the example provided in \S\ref{examplesec}.

\subsubsection{A conjecture and relation to other work}

We propose the following conjecture for theoretical and mathematical physicists.

\begin{Conjecture} {\bf (for physicists)}\label{physconj}
Any reasonable QFT with a gapped dimension spectrum satisfies the OPE in the strong sense expressed by the
EFNNB (\ref{efnnb}), without the $\frac{d}{2}$ restriction on scaling dimensions.
\end{Conjecture}

As a matter of convenience, we imposed that the set of labels $\mathcal{A}$ for our abstract system
of pointwise correlations be finite.
However, the above formalism also makes sense if $\mathcal{A}$ is countably infinite, provided $\forall \Delta\in\mathbb{R}$, $\{A\in\mathcal{A}\ |\ [A]\le \Delta\}$ is finite. This is what we mean by having a gapped dimension spectrum. This condition is also related to the first axiom in the definition of regularity structures~\cite[Definition 2.1]{Hairer} with the caveat that Hairer's homogeneity exponents $\alpha$ essentially
are minus our scaling dimensions $[A]$. 
 
Note that the OPE is believed by physicists to be true (in the weaker sense of (\ref{opelittleo})) for any
reasonable QFT. Because of this generality, one may wonder if the OPE would follow from the
Wightman axioms,
but this is still unclear~\cite{SchliederS,Adetunji}. See~\cite{SchroerSV,Luscher,Mack} for the better behaved case of CFT.
As for the OPE in the sense of formal power series, there has been recent progress by Hollands and his collaborators (see, e.g.,~\cite{HollandsK,FrobHH,FrobH}).
This work uses the flow equation methods developed earlier by Keller and Kopper~\cite{KellerK1,KellerK2}.
However, these are mainly Fourier or momentum space methods. For example, in~\cite{HollandsK},
the authors prove the convergence of the OPE, at every finite order of perturbation theory,
without limitation on the distance between the two points being collapsed.
This is because the spectator fields are smeared with test functions that have compact support
in Fourier space. Take the simple example from \S\ref{introtointro} and consider the partially smeared correlation
$\langle \phi(x_1)\phi(x_2)\phi(f_3)\phi(f_4)\rangle$
where $\widehat{f}_3$ and $\widehat{f}_4$ are compactly supported. Then for fixed $x_2$, this is real-analytic in $x_1$
on all of $\mathbb{R}^d\backslash\{x_2\}$. If $\widehat{f}_3$ and $\widehat{f}_4$ have constant modulus equal to one, as is the case
when localizing at points $x_3$ and $x_4$, then the singularity at $|x_1-x_2|=\min\{|x_2-x_3|,|x_2-x_4|\}$ is indeed present.
For the proof of the above conjecture, in perturbative QFT, which is formulated in $x$-space, we believe that position space methods are more
appropriate. Nevertheless the above-mentioned works manage to go quite far.
In particular,~\cite[Theorem 3]{FrobHH}
essentially proves the analogue of our conjecture above for the BFNNB. The bound~\cite[Theorem 3]{FrobHH} is {\em not} in factorized form,
but it is {\em grosso modo} equivalent to (\ref{bfnnb}).
We believe the most powerful techniques for proving our conjecture in the sense of perturbative QFT are the
multiscale $x$-space methods for BPHZ renormalization~\cite{FeldmanMRS2}
developed by the \'Ecole Polytechnique school of constructive QFT.
The power of these methods resides in the fact they can be adapted (albeit {\em with} tears) to the non-perturbative
situation~\cite{FeldmanMRS3,RivasseauBook,AbdesselamR,AbdesselamX,Unterberger}.

One could also mention
that the methods from~\cite{FeldmanMRS2}
play an important role in the recent result by Chandra and Hairer~\cite{ChandraH}.
Note that with simple modifications, our result can be applied to random space-time distributions 
$u(x,t)$, e.g., ones arising from singular SPDEs. One needs to change, in the following estimates,
Euclidean norms $|\cdot|$ by anisotropic ``norms'' (e.g., suitable for parabolic scaling) $||\cdot||_{\mathfrak{s}}$ as defined in~\cite[\S2.2]{Hairer}. The length of the multiindex $\mathfrak{s}$ should replace $d$ in the local integrability conditions, etc. One can also easily work with the more general spaces of distributions $\mathcal{D}'(U)$ instead of $\mathcal{S}'(\mathbb{R}^d)$. However, a difficulty is that one would need more than a local in time well-posedness result which is valid up to a random blow-up time. This is because our approach uses moments as input and they are required to exist at least on a tiny domain $U$. Nevertheless, a good
test example for such an adaptation of our theorem is for a solution $u(x,t)$ of the {\em linear} heat equation, following the treatment in \S\ref{examplesec}. In general, it would be of great interest to investigate the (dynamical) OPE structure in space-time for solutions of singular SPDEs.

While singular SPDEs are part of the motivation for this article, our main result was mainly designed for CFT.
From a careful reading of \S\ref{examplesec} and the heuristic given at the end of \S\ref{motivsec} it would be reasonable to expect that, in the CFT case, the EFNNB should follow from the convergence of the OPE (as a series without the
$\Delta$ cut-off on scaling dimensions) and from the BFNNB (for a larger collection of fields).
In \S\ref{examplesec} we prove the hypotheses of our theorem in the case of the FMFF. Although this is a Gaussian example, the construction of local Wick monomials given
by \S\ref{examplesec} and Theorem \ref{maintheorem} is non-Gaussian in spirit. In particular, we do not rely on hypercontractivity as is done, e.g., in~\cite{DaPratoT,EJS}.
This example is interesting since it corresponds to a generalized free field~\cite[\S6.2]{GlimmJbook}. Indeed, the FMFF satisfies Osterwalder-Schrader or reflection positivity when $[\phi]\ge\max\{0,\frac{d-2}{2}\}$ (see, e.g.,~\cite[Lemma 2.1]{FrankL}).
This example is also called a CFT of mean field type in the physics literature~\cite{KarateevKSD}.
It is nontrivial from the point of view of the conformal bootstrap because the same kinematic
(position dependent) building blocks (conformal blocks) feature in Gaussian as well as non-Gaussian CFTs.
Only different values for purely numerical quantities like scaling dimensions and OPE coefficients (in the CFT sense, not our $\mathcal{C}_{A,B}^{C}$) distinguish these two types of CFTs.
For completeness, we end \S\ref{examplesec} with a proof of global conformal invariance for the FMFF example.
This was previously sketched in~\cite{Abdesselam3DCP} to which we refer the reader for more background.
When $d=2$ and for some smooth volume cut-off function $g(x)$, our result also shows that
$\phi\longmapsto \exp\left(-\int_{\mathbb{R}^2} dx\ g(x) :\phi^4(x):\right)$
is well-defined and Borel measurable when $0<[\phi]<\frac{1}{4}$. Hence, one should be able to use Nelson's old approach for constructing the fractional $\phi_2^4$ model, at least in the massive case. The more difficult massless case is an important test example for the conformal bootstrap
approach~\cite{PaulosRRZ,BehanRRZ1,Behan}.

An interesting and beautiful example for the BFNNB and the EFNNB is the 2D Ising CFT.
Up to multiplying the (spin) field by a constant, this corresponds to the {\em very explicit}
collection of pointwise correlations given by
\[
\langle\phi(x_1)\cdots\phi(x_n)\rangle
=\sqrt{\sum_{q}\prod_{1\le i<j\le n}|x_i-x_j|^{\frac{q_i q_j}{2}}}
\]
where the sum is over ``neutral charge configurations'' $q=(q_i)_{1\le i\le n}\in\{-1,1\}^n$ such that
$\sum_{i=1}^{n}q_i=0$.
These correlations satisfy the BFNNB, as can be established with at least five different proofs.
A first proof follows from the observation that, in the notations of Eq. (\ref{explicitdef}) below,
one has
\begin{equation}
\langle\phi(x_1)\cdots\phi(x_n)\rangle\le \sum_{\mathcal{W}}
\prod_{\{a,b\}\in\mathcal{W}} \langle\phi(x_a)\phi(x_b)\rangle\ .
\label{Newmanineq}
\end{equation}
This follows from the Gaussian correlation inequality at the lattice level~\cite{Newman}, then passing to the limit using the results of~\cite{Dubedat,ChelkakHI}. The BFNNB with $\epsilon=0$ and $k=0$
then immediately follows from (\ref{Newmanineq}) as in the proof of Proposition \ref{nonlocalbd} in
\S\ref{examplesec} below.
A second proof uses~\cite[Proposition 3.5]{HairerS} which
is a bound such as (\ref{Newmanineq}) but with a possibly large constant in front.
A third proof of the BFNNB can be found in~\cite[\S3.3]{FurlanM} using properties of the random cluster model.
A fourth proof~\cite[Appendix A]{LacoinRV} uses the Gale-Shapely stable marriage theorem. Finally, 
a fifth proof~\cite[\S3.2]{JunnilaSW} can be obtained by looking at an individual charge configuration, taking logarithms and reducing the problem to an electrostatic inequality
(in the tradition of Onsager~\cite{Onsager}, Baxter~\cite{Baxter},
Lieb-Yau~\cite{LiebY}) for the two-dimensional Coulomb potential.
Although one cannot use Theorem \ref{maintheorem} to construct the square or energy field $\phi^2$ as a random distribution
(since it has scaling dimension equal to  $d/2$), it would still
be interesting to establish the EFNNB for it, i.e., for the alphabet $\mathcal{A}=\{\bbone,\phi,\phi^2\}$. The Ising CFT is the first nontrivial unitary minimal model $\mathcal{M}(4,3)$, the next one in order of complexity being the tricritical Ising CFT or $\mathcal{M}(5,4)$.
The latter has a spin field with $[\sigma]=\frac{3}{80}$ and two energy-like fields with
$[\varepsilon]=\frac{1}{10}$, $[\varepsilon']=\frac{3}{5}$ produced by the OPE which collapses two $\sigma$'s (see~\cite[\S7.4.3]{DiFrancescoMS}). Pending progress similar
to~\cite{Dubedat,ChelkakHI,CamiaGN} for the tricritical Ising model, the latter could in the future provide an example where our theorem applies and shows that
$\varepsilon,\varepsilon'$ exist as random distributions and are deterministic functionals of the spin field $\sigma$.
Another very interesting
example of a system of
pointwise correlations, though without the permutation symmetry (\ref{symmetry}),
is that of multiple SLE pure partition
functions~\cite{KytolaP,Wu,PeltolaW}. This system is now known to
satisfy the BFNNB
when $0<\kappa\le 6$. We refer to the
recent review~\cite{Peltola} for more details.
Finally, note that our setting does not cover CFTs with a continuous dimension spectrum such as Liouville theory. It would be interesting to investigate analogues of the BNNFB and the ENNFB for this situation too. For results in this direction, see the fusion estimates~\cite[Proposition 5.1]{KupiainenRV} and~\cite[Lemma 3.1]{Oikarinen}.

\subsubsection{Structure of the article}
In \S\ref{warmupsec} we show how the BFNNB implies local integrability for correlations.
It uses elementary estimates and in particular Lemma \ref{globalbetabd} as well as
Lemma \ref{localbetabd}
which are the workhorses on which the proof of Proposition \ref{combiprop} relies.
These are very special cases of~\cite[Lemma 10.14]{Hairer} which is also an important tool for
applications of the theory of regularity structures.
The easy proof in \S\ref{warmupsec} should give the gist of our method for proving Proposition \ref{combiprop} from which
Theorem \ref{maintheorem} follows.
In \S\ref{probatocomb}, we explain how to reduce the probability theory statements of Theorem \ref{maintheorem}
to the purely combinatorial estimate in Proposition \ref{combiprop}.
In order to deal with the space $\mathcal{O}_{\rm M}$, we devised an approach which we cannot resist
calling the ``multiply and conquer'' strategy: multiply by a generic function in $\mathcal{S}$, prove the needed
kernel theorem or special case of~\cite[Proposition 34]{SchwartzVect1} only using $\mathcal{S}$
and $\mathcal{S}'$, and finally undo the dammage thanks to the multiplier space characterization of $\mathcal{O}_{\rm M}$.
Section \ref{probatocomb} below gives an example of how this strategy works. 
For the main result
of this article, the core argument is in \S\ref{mainestimatesec} where Proposition \ref{combiprop} is proved.
In the rather short \S\ref{finishproof}, we finish the proof of Theorem \ref{maintheorem}.
Finally, in \S\ref{examplesec} the example of the fractional massless free field is treated in detail, including a proof of global conformal invariance.
Note that while the proofs in \S\ref{examplesec} are rigorous, we adopted a semi-formal style of presentation for this last section. This is because we felt that a more formal
handling of the combinatorics involved would have made this section harder rather than easier to understand.

\section{A warm-up example: local integrability}\label{warmupsec}

In this section we give some elementary but important lemmas for estimating analogues of Euler beta integrals. These are then immediately put to good use, in order to show how the BNNFB implies the local integrability of pointwise correlations stated in \S\ref{probareal}.
This is done with a simple pin-and-sum argument which will hopefully serve as useful pedagogical preparation for the more involved double pin-and-sum argument required for the proof of the main theorem and explained in \S\ref{doublepinsumsec} below. 

\subsection{Elementary beta integral estimates}\label{elementarybds}

Throughout the remainder of this article,
we will use the notation $\bbone\{\cdots\}$ for the sharp indicator function of the condition between braces.

\begin{Lemma}\label{globalL1bd}
$\forall \alpha\in[0,d)$, $\forall \beta\in (d,\infty)$, $\exists K>0$, $\forall x\in\mathbb{R}^d$
\[
\int_{\mathbb{R}^d\backslash\{x\}} dy
\ \langle y\rangle^{-\beta}
|x-y|^{-\alpha}\le K\ .
\]
\end{Lemma}

\noindent{\bf Proof:}
Since $\beta\ge d>0$ and $\langle y\rangle\ge 1$, we have
\[
\int_{\mathbb{R}^d\backslash\{x\}} dy\  
\frac{\bbone\{|x-y|\le 1\}}{\langle y\rangle^{\beta}|x-y|^{\alpha}}
\le \int_{\mathbb{R}^d\backslash\{x\}} dy\ 
\frac{\bbone\{|x-y|\le 1\}}{|x-y|^{-\alpha}}
=\frac{1}{d-\alpha}\times\frac{2\pi^{\frac{d}{2}}}{\Gamma\left(\frac{d}{2}\right)}\ .
\]
On the other hand, $\alpha \ge 0$ implies
\[
\int_{\mathbb{R}^d\backslash\{x\}} dy\ 
\frac{\bbone\{|x-y|> 1\}}{\langle y\rangle^{\beta}|x-y|^{\alpha}}
\le 
\int_{\mathbb{R}^d\backslash\{x\}} dy\ \bbone\{|x-y|>1\}\ \langle y\rangle^{-\beta}
\le \int_{\mathbb{R}^d\backslash\{x\}} dy\ \langle y\rangle^{-\beta}\ .
\]
Combining both pieces gives the wanted bound with
\[
K=\frac{2\pi^{\frac{d}{2}}}{\Gamma\left(\frac{d}{2}\right)}\left[
\frac{1}{d-\alpha}+\int\limits_{0}^{\infty}dr\frac{r^{d-1}}{(1+r^2)^{\frac{\beta}{2}}}\right]\ .
\]
\qed

The following lemma is the fundamental tool and resembles an estimate for the beta function.

\begin{Lemma}\label{globalbetabd}
$\forall \alpha,\beta\in[0,\frac{d}{2})$, $\forall \gamma\in (d,\infty)$,
$\exists K>0$, $\forall x,z\in\mathbb{R}^d$,
\[
\int_{\mathbb{R}^d\backslash\{x,z\}} dy\  \langle y\rangle^{-\gamma}
|x-y|^{-\alpha}|y-z|^{-\beta}\le K\ .
\]
\end{Lemma}
Note that we allow the case $x=z$ although it already follows from Lemma \ref{globalL1bd}.
  
\noindent{\bf Proof:}
If $|x-y|\le|y-z|$ then $|y-z|^{-\beta}\le|x-y|^{-\beta}$ because $\beta\ge 0$.
Therefore,
\[
\int_{\mathbb{R}^d\backslash\{x,z\}} dy\ \bbone\{|x-y|\le|y-z|\}
\ \langle y\rangle^{-\gamma}
|x-y|^{-\alpha}|y-z|^{-\beta}
\le
\]
\[
\int_{\mathbb{R}^d\backslash\{x,z\}} dy\  \bbone\{|x-y|\le|y-z|\}
\ \langle y\rangle^{-\gamma}
|x-y|^{-(\alpha+\beta)}
\le
\int_{\mathbb{R}^d\backslash\{x\}} dy 
\ \langle y\rangle^{-\gamma}
|x-y|^{-(\alpha+\beta)}\le K'
\]
where $K'$ is provided by the previous lemma with $\alpha+\beta$ in lieu of $\alpha$.
By symmetry, the other piece
\[\int_{\mathbb{R}^d\backslash\{x,z\}} dy\ \bbone\{|x-y|>|y-z|\}
\ \langle y\rangle^{-\gamma}
|x-y|^{-\alpha}|y-z|^{-\beta}
\]
is bounded by the same $K'$ since $\alpha\ge 0$.
The lemma follows with $K=2K'$.
\qed

We also need local versions of the two lemmas where the points are restricted to a closed Euclidean ball $\bar{B}(0,R)$ with $R>0$.
Note that the range of exponents is larger.

\begin{Lemma}\label{localL1bd}
$\forall \alpha\in(-\infty,d)$, $\exists K>0$, $\forall R>0$,
$\forall x\in\bar{B}(0,R)$,
\[
\int_{\bar{B}(0,R)\backslash\{x\}} dy
\ |x-y|^{-\alpha}\le K R^{d-\alpha}\ .
\]
\end{Lemma}

\noindent{\bf Proof:}
Since $|x-y|\le |x|+|y|\le 2R$, the integral is bounded by
\[
\int_{\bar{B}(x,2R)\backslash\{x\}} dy
\ |x-y|^{-\alpha}= K R^{d-\alpha}
\ \ 
{\rm with}\ \ 
K=\frac{2^{d-\alpha}}{d-\alpha}\times\frac{2\pi^{\frac{d}{2}}}{\Gamma\left(\frac{d}{2}\right)}\ .
\]
\qed

\begin{Lemma}\label{localbetabd}
$\forall \alpha,\beta\in(-\infty,\frac{d}{2})$,
$\exists K>0$, $\forall R>0$, $\forall x,z\in \bar{B}(0,R)$,
\[
\int_{\bar{B}(0,R)\backslash\{x,z\}} dy
\ |x-y|^{-\alpha}|y-z|^{-\beta}\le K R^{d-\alpha-\beta}
\]
\end{Lemma}

\noindent{\bf Proof:}
Suppose first that $\alpha<0$. Then, on the integration domain,
$|x-y|^{-\alpha}\le (2R)^{-\alpha}$ holds. Thus,
\[
\int_{\bar{B}(0,R)\backslash\{x,z\}} dy
\ |x-y|^{-\alpha}|y-z|^{-\beta}\le (2R)^{-\alpha} \int_{\bar{B}(0,R)\backslash\{x,z\}} dy\ |y-z|^{-\beta}
\]
and the previous lemma gives the desired bound. Likewise if $\beta<0$ we use $|y-z|^{-\beta}\le (2R)^{-\beta}$ 
and obtain the same conclusion. So we can assume that both $\alpha$ and $\beta$ are nonnegative.
We then have
\[
\int_{\bar{B}(0,R)\backslash\{x,z\}} dy\ 
\frac{\bbone\{|x-y|\le|y-z|\}}{|x-y|^{\alpha}|y-z|^{\beta}}
\le
\int_{\bar{B}(0,R)\backslash\{x,z\}} dy\ 
\frac{\bbone\{|x-y|\le|y-z|\}}{|x-y|^{\alpha+\beta}}
\]
and bound this with the previous lemma with $\alpha+\beta$ instead of $\alpha$.
The other piece
\[
\int_{\bar{B}(0,R)\backslash\{x,z\}} dy\ 
\frac{\bbone\{|x-y|>|y-z|\}}{|x-y|^{\alpha}|y-z|^{\beta}}
\le
\int_{\bar{B}(0,R)\backslash\{x,z\}} dy\ 
\frac{\bbone\{|x-y|>|y-z|\}}{|y-z|^{\alpha+\beta}}
\]
satisfies the same desired bound.
\qed

\subsection{The pin and sum argument with hairy cycles}\label{pinsum}

We use the notation
\[
||f||_{\alpha,k}=\sup_{x\in\mathbb{R}^d} \langle x\rangle^k |\partial^{\alpha}f(x)|
\]
for the defining seminorms of $\mathcal{S}(\mathbb{R}^d)$ indexed by the integer $k\in\mathbb{N}$ and multiindex $\alpha\in\mathbb{N}^d$.
Note that if $x\in\mathbb{R}^m$ and $y\in\mathbb{R}^n$
then $\langle x,y\rangle^2=1+|x|^2+|y|^2\le\langle x\rangle^2 \langle y\rangle^2$. More generally, we have the following nice factorization property
for concatenation of points or vectors
$\langle x_1,\ldots,x_n\rangle\le \langle x_1\rangle\cdots \langle x_n\rangle$.

Let $f(x_1,\ldots,x_p)$ be a function in $\mathcal{S}_{x_1,\ldots,x_p}(\mathbb{R}^{pd})$. Then from the BFNNB
with
\[
\epsilon=\frac{1}{2}\min_{1\le i\le p} \left(\frac{d}{2}-[B_i]\right)\ ,
\]
\begin{flalign*}
&{\rm we\ get}\ \ & &
\int_{{\rm Conf}_p}\prod_{i=1}^{p}dx_i\ |f(x_1,\ldots,x_p)|\times|\langle \mathcal{O}_{B_1}(x_1)\cdots
\mathcal{O}_{B_p}(x_p)\rangle|\le 
K ||f||_{0,k+d+1} \times
\mathcal{I} & & \ \\
&{\rm where}\ \ & &
\mathcal{I}=
\int_{{\rm Conf}_p}\prod_{i=1}^{p}dx_i\ \prod_{i=1}^{p}\langle x_i\rangle^{-(d+1)}
\prod_{i=1}^{p}\left(\min_{j\neq i}|x_i-x_j|\right)^{-([B_i]+\epsilon)}\ . & & \ 
\end{flalign*}
For $p=1$, there is nothing to prove so we assume $p\ge 2$.
Let $\mathcal{N}$ denote the set of fixed-point-free endofunctions of $[p]=\{1,\ldots,p\}$, namely, all maps $\tau:[p]\rightarrow[p]$
such that $\tau(i)\neq i$ for all $i$.
Then for any fixed configuration of points $(x_1,\ldots,x_p)\in{\rm Conf}_p$, we have
\[
1\le\sum_{\tau\in\mathcal{N}}\prod_{i=1}^{p}
\bbone\left\{
|x_i-x_{\tau(i)}|=\min_{j\neq i}|x_i-x_j|
\right\}\ .
\]
Indeed, for each point $i$ we can choose a nearest neighbor which we call $\tau(i)$.
Insert the inequality inside the integral, then
use the equalities in the indicator functions to replace the min's in the bound by $|x_i-x_{\tau(i)}|$'s. Then drop the indicator functions
and pull the sum out of the integral. Therefore,
\[
\mathcal{I}\le \sum_{\tau\in\mathcal{N}}
\int_{{\rm Conf}_p}\prod_{i=1}^{p}dx_i\ \prod_{i=1}^{p}\langle x_i\rangle^{-(d+1)}
\prod_{i=1}^{p}|x_i-x_{\tau(i)}|^{-([B_i]+\epsilon)}\ .
\]
One can draw a directed graph on the vertex set $[p]$ with edges $i\rightarrow\tau(i)$, for $i\in[p]$.
Since $\tau$ is an endofunction, the connected components are what we call ``hairy cycles''.
Namely, each component is made of a central cycle playing the role of ``root'' for a collection of trees attached to it
and oriented towards it.
Then eliminate the trees by recursively using Lemma \ref{globalL1bd}, starting from the leaves.
The remaining cycle has length at least two by the fixed-point-free condition. Use Lemma \ref{globalbetabd}
to open the cycle and erase two consecutive edges. 
Finally, eliminate the chain left over by Lemma \ref{globalL1bd} including its $\alpha=0$ case for the last point.
As a result
\[
\int_{{\rm Conf}_p}\prod_{i=1}^{p}dx_i\ |f(x_1,\ldots,x_p)|\times|\langle \mathcal{O}_{B_1}(x_1)\cdots
\mathcal{O}_{B_p}(x_p)\rangle|\le 
O(1)\ ||f||_{0,k+d+1} 
\]
and the function $\langle \mathcal{O}_{B_1}(x_1)\cdots
\mathcal{O}_{B_p}(x_p)\rangle$ canonically defines an element of $\mathcal{S}'_{x_1,\ldots,x_p}(\mathbb{R}^{pd})$.
Note that we will use $O(1)$ in order to denote constants which do not need to be made explicit.

As a picture is worth a thousand words, let us show on an example how the pin and sum
argument works.
\[
\parbox{10cm}{
\includegraphics[width=10cm]{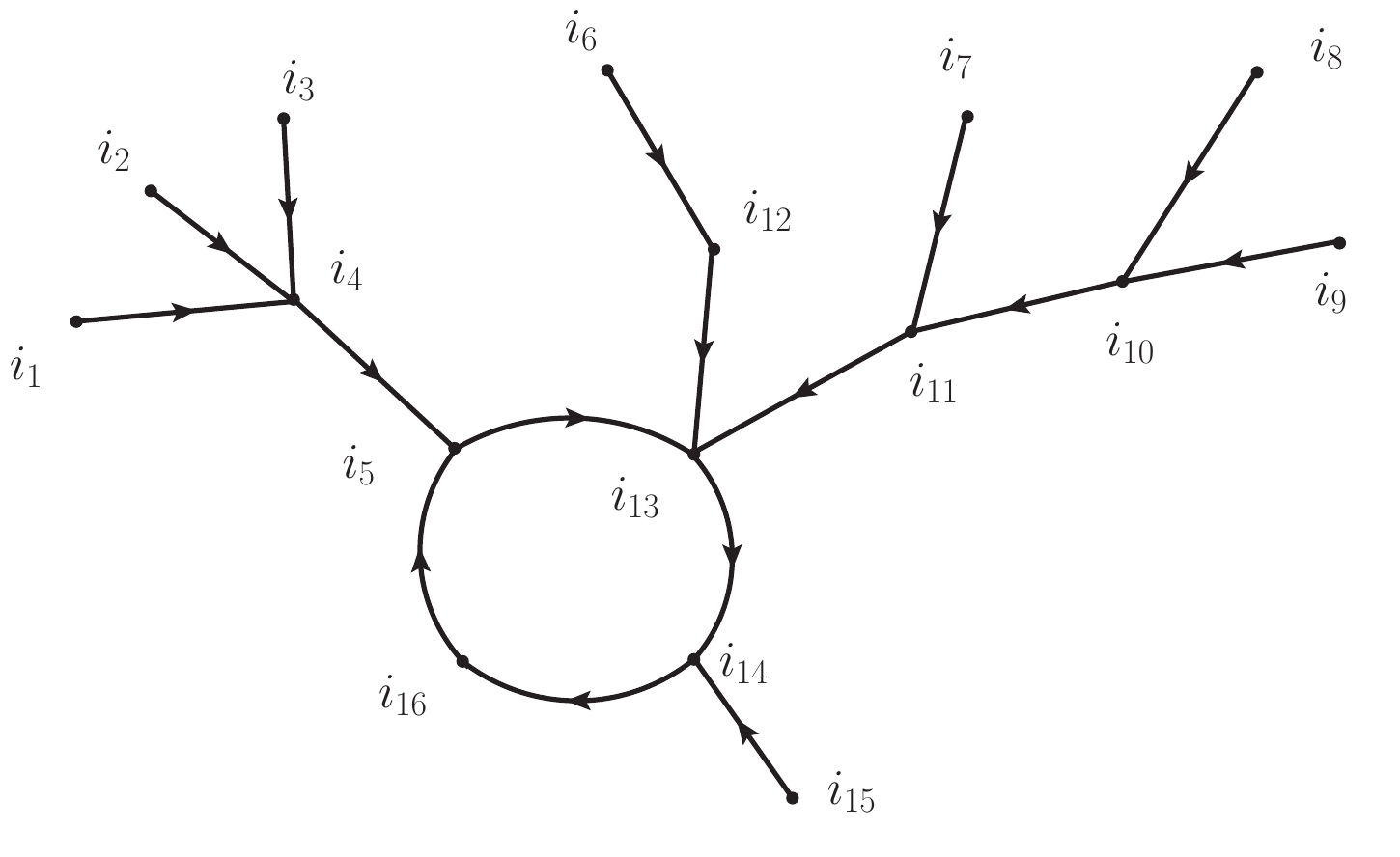}}
\]
Suppose for instance that the root is the label $i_{16}$. Namely, $x_{i_{16}}$ is the last variable to be integrated.
Then a possible order of integration is given by the succession
\[
i_1,i_2,i_3,i_4,i_6,i_{12},i_7, i_8,i_9,i_{10},i_{11},i_{15},i_5,i_{13},i_{14},i_{16}. 
\]
All these integrations are done with the
help of Lemma \ref{globalL1bd},
except when opening the cycle, i.e., integrating over $x_{i_{5}}$ which uses Lemma \ref{globalbetabd}.

\begin{Remark}\label{contrarian}
An important variant which will be used later is that one does not have to follow the arrows and pick the root in the cycle.
One can be a contrarian and
decide that the root is say $i_{10}$, i.e., that $x_{i_{10}}$ is integrated last.
The argument works just as well, if one chooses for instance the order of integration
$i_9,i_8,i_7,i_6,i_{12},i_{15},i_1,i_2,i_3,i_4,i_{16},i_5,i_{14},i_{13},i_{11},i_{10}$.
Indeed, when one is left with just the cycle attached by a path to the root, one can always open the cycle
at a vertex which is different from the one where the path touches the cycle. This is again because the latter
has length at least two.
\end{Remark}

\begin{Remark}
This kind of pin and sum argument is rather common
in constructive QFT~\cite{RivasseauBook}.
A similar use of the nearest neighbor function $\tau$ is in~\cite{GunsonP}.
\end{Remark}

\section{From probability to combinatorics}\label{probatocomb}

\subsection{Some functional-analytic technical lemmas}
We pick up the thread and notations from \S\ref{resultssec}.
The formulation of Theorem \ref{maintheorem} assumes that the random distribution $M_r$ is well defined, belongs to $\mathcal{O}_{\rm M}$ and is such that $M_r(f)\in L^p$ for all $p\ge 1$ and all test functions $f$.

\begin{Lemma}
The function $g_r$ from \S\ref{resultssec} is such that
 $\psi(x) g_r(x,z)$
belonds to $\mathcal{S}_{x,z}(\mathbb{R}^{2d})$ for all $\psi(x)$ in 
$\mathcal{S}_{x}(\mathbb{R}^d)$.
\end{Lemma}

\noindent{\bf Proof:}
Note that $g_r$
can be written as
$g_r(x,z)=\rho_r(x-z)\times\langle \mathcal{C}_{AB}^{C}(y,z),\rho_r(x-y)\rangle_y$.
More generally, we can define
$G_r(x,z,w)=\rho_r(x-z)\times\langle \mathcal{C}_{AB}^{C}(y,z),\rho_r(w-y)\rangle_y$
for $(x,z,w)\in\mathbb{R}^{3d}$. 
Let $\psi_1,\psi_2\in\mathcal{S}(\mathbb{R}^d)$, then
\[
\psi_1(x)\psi_2(w) G_r(x,z,w)=\psi_1(x)\rho_r(x-z)\times\langle \mathcal{C}_{AB}^{C}(y,z),\psi_2(w)\rho_r(w-y)\rangle_y
\]
belongs to $\mathcal{S}_{x,z,w}(\mathbb{R}^{3d})$, as results from our hypothesis (SH2) in \S\ref{softhyp}.
By restriction to the subspace $x=w$, we get that $\psi_1(x)\psi_2(x) g_r(x,z)$
belongs to $\mathcal{S}_{x,z}(\mathbb{R}^{2d})$.
However, by a lemma of Miyazaki~\cite[Lemma 1]{Miyazaki} (see also~\cite{PetzeltovaV,Voigt,Garrett1}), every $\psi(x)$ in $\mathcal{S}_{x}(\mathbb{R}^d)$
can be written as $\psi_1(x)\psi_2(x)$ with $\psi_1(x)$, $\psi_2(x)$
in $\mathcal{S}_{x}(\mathbb{R}^d)$. Thus, $\psi(x) g_r(x,z)$
belongs to $\mathcal{S}_{x,z}(\mathbb{R}^{2d})$ for all $\psi(x)$ in 
$\mathcal{S}_{x}(\mathbb{R}^d)$.
\qed

Now consider the integral used to define $Z_r(x)$, in \S\ref{resultssec}, namely,
\[
Y_r(x)=\int_{{\rm Conf}_2}dy\ dz\ \rho_r(x-y)\rho_r(x-z)\ \mathcal{C}_{AB}^{C_{\ast}}(y,z)
=\langle 1, g_r(x,z)\rangle_z
\] 
in the special case $C=C_{\ast}$.

\begin{Lemma}
The function $Y_r(x)$ is strictly positive and belongs to
$\mathcal{O}_{{\rm M},x}(\mathbb{R}^d)$. Likewise, $Z_r(x)=Y_r(x)^{-1}$ belongs to
$\mathcal{O}_{{\rm M},x}(\mathbb{R}^d)$.
\end{Lemma}

\noindent{\bf Proof:}
If $\psi\in\mathcal{S}(\mathbb{R}^d)$, then
$\psi(x)Y_r(x)=\langle 1, \psi(x)g_r(x,z)\rangle_z\ \in\ \mathcal{S}_{x}(\mathbb{R}^d)$,
by the previous lemma and by
Fubini's Theorem for distributions~\cite[Theorem IV, p. 108]{SchwartzBook1}.
From the multiplier space characterization of $\mathcal{O}_{\rm M}$ (see~\cite[Proposition 5, p. 417]{Horvath}
or~\cite[Proposition 1.6.1]{Amann}), this implies that $Y_r(x)$ is in $\mathcal{O}_{{\rm M},x}(\mathbb{R}^d)$
and therefore smooth.
For taking the inverse, the soft hypotheses from \S\ref{softhyp} are not enough. 
From (\ref{nondeg}) we get
\[
Y_{r}(x)\ge
\int_{{\rm Conf}_2}
dy\ dz\ \rho_{r}(x-y)\rho_{r}(x-z)
\frac{1}{K|y-z|^{[A]+[B]-[C_{\ast}]-\epsilon}}\langle y\rangle^{-k}\langle z\rangle^{-k}
\]
\begin{flalign*}
&{\rm or\ rather}\ \ & &
Y_{r}(x)\ge
\int_{{\rm Conf}_2}
du\ dv\ \rho_{r}(u)\rho_{r}(v)
\frac{1}{K|u-v|^{[A]+[B]-[C_{\ast}]-\epsilon}}
\langle x-u\rangle^{-k}\langle x-v\rangle^{-k}& & \ 
\end{flalign*}
after the change of variables $y=x-u$ and $z=x-v$.
Note that if $|u|\le 1$ then
\[
\langle x-u\rangle^2\le 1+(|x|+1)^2\le 1+2(|x|^2+1)\le 3\langle x\rangle^2
\]
so that $\langle x-u\rangle\le \sqrt{3}\langle x\rangle$.
From now on we assume $r\le 0$. Because of the support condition on the mollifier $\rho$, this results in
\[
Y_{r}(x)\ge K^{-1}3^{-k}\langle x\rangle^{-2k}
\int_{{\rm Conf}_2}
du\ dv\ \rho_{r}(u)\rho_{r}(v)
\frac{1}{|u-v|^{[A]+[B]-[C_{\ast}]-\epsilon}}\ .
\]
From the scaling change of variables $u=L^{r}\tilde{u}$
and $v=L^{r}\tilde{v}$
we obtain
\[
Y_{r}(x)\ge \tilde{K}\langle x\rangle^{-2k} L^{-r([A]+[B]-[C_{\ast}]-\epsilon)}
\]
\begin{flalign*}
&{\rm where}\ \ & &
\tilde{K}=K^{-1} 3^{-k}\int_{{\rm Conf}_2}
d\tilde{u}\ d\tilde{v}\ \rho(\tilde{u})\rho(\tilde{v})
|\tilde{u}-\tilde{v}|^{-[A]-[B]+[C_{\ast}]+\epsilon}>0 & & \ 
\end{flalign*}
from the assumptions on the mollifier $\rho$.
Thus $Y_r(x)>0$ and $Z_r(x)$ is a well defined smooth function of $x$ which satisfies the estimate
\[
Z_{r}(x)\le O(1)\langle x\rangle^{2k} L^{r([A]+[B]-[C_{\ast}]-\epsilon)}\ .
\]
Since derivatives of $Z_r(x)=Y_r(x)^{-1}$ are polynomials in $Z_r(x)$ and derivatives of $Y_r(r)$, we immediately conclude that 
$Z_r(x)$ belongs to $\mathcal{O}_{{\rm M},x}(\mathbb{R}^d)$.
\qed

\begin{Lemma}
The random distribution $M_r$ from \S\ref{resultssec} is well defined, belongs to $\mathcal{O}_{\rm M}$ and is such that $M_r(f)\in L^p(\Omega,\mathcal{F},\mathbb{P})$ for all $p\ge 1$ and all test function $f\in\mathcal{S}(\mathbb{R}^d)$.
\end{Lemma}

\noindent{Proof:}
For $f\in\mathcal{S}(\mathbb{R}^d)$, let
\[
M_{A,B,r}^{C}(f)=\int_{\mathbb{R}^d}dx\ 
Z_r(x) \widetilde{\mathcal{O}}_{C,r}(x) f(x)=
\langle Z_r(x),\langle \mathcal{O}_C(z), g_r(x,z)\rangle_z f(x)\rangle_x
\]
\[
=\langle Z_r(x),\langle \mathcal{O}_C(z), f(x)g_r(x,z)\rangle_z \rangle_x\ .
\]
From the above considerations, $f(x)g_r(x,z)\in \mathcal{S}_{x,z}(\mathbb{R}^{2d})$
and so by Fubini's Theorem for distributions $\langle \mathcal{O}_C(z), f(x)g_r(x,z)\rangle_z$ is in
$\mathcal{S}_x(\mathbb{R}^d)$. Moreover,
\begin{equation}
M_{A,B,r}^{C}(f)=\langle Z_r(x)\otimes \mathcal{O}_C(z), f(x) g_r(x,z)
\rangle_{x,z}
=\langle \mathcal{O}_C(z), h_{A,B,r}^{C}(z)\rangle_z
\label{ocrep}
\end{equation}
\begin{flalign*}
&{\rm where}\ \ & &
h_{A,B,r}^{C}(z)=\langle Z_r(x), f(x) g_r(x,z)\rangle_x=\int_{\mathbb{R}^d}dx\ 
Z_r(x)f(x)g_r(x,z) & & \ 
\end{flalign*}
is a {\em fixed} test function in $\mathcal{S}_z(\mathbb{R}^d)$, while $\mathcal{O}_C$ is random and depends on $\omega\in\Omega$.
From our hypotheses on probabilistic incarnations and the representation (\ref{ocrep}),
it follows that $M_{A,B,r}^{C}(f)$
is a well defined random variable on $\Omega$ with finite moments of all orders.

We now define
\[
M_{A,B,r}(f)=\int_{\mathbb{R}^d}dx\ 
Z_r(x) \mathcal{O}_{A,r}(x) \mathcal{O}_{B,r}(x) f(x)
\]
\[
=\langle Z_r(x), \langle \mathcal{O}_A(y),\rho_r(x-y)\rangle_y 
\langle \mathcal{O}_B(z),\rho_r(x-z)\rangle_z f(x)\rangle_x
\]
\[
=
\langle Z_r(x)\otimes  \mathcal{O}_A(y)\otimes  \mathcal{O}_B(z),
\rho_r(x-y)\rho_r(x-z)f(x)
\rangle_{x,y,z}
\]
by Fubini's Theorem for distributions.
This can be rewritten as
\begin{equation}
M_{A,B,r}(f)=\langle\mathcal{O}_A(y)\otimes  \mathcal{O}_B(z), h_{A,B,r}(y,z)
\rangle_{y,z}
\label{oaobrep}
\end{equation}
\begin{flalign*}
&{\rm where}\ \ & &
h_{A,B,r}(y,z)=\int_{\mathbb{R}^r}dx\ 
Z_r(x)\rho_r(x-y)\rho_r(x-z)f(x) & & \ 
\end{flalign*}
is a fixed test function in $\mathcal{S}_{y,z}(\mathbb{R}^d)$.
Since the tensor product of distributions is continuous for the strong topology, we have that $M_{A,B,r}(f)$ is $\mathcal{F}$-measurable.
By~\cite[Theorem III.7.1]{Fernique} and its corrolary, we also conclude that this random variable has moments of all orders.
\qed

\subsection{The main combinatorial proposition}
In order to proceed further, we need to generalize the setting of \S\ref{resultssec}
by jugling several renormalized product constructions at the same time.
For $1\le i\le m$, we pick $A_i$, $B_i$ in $\mathcal{B}$ and $C_{\ast i}\in\mathcal{A}\backslash\mathcal{B}$.
We set $\Delta_i=[C_{\ast i}]$ and assume as before that $\mathcal{A}(\Delta_i)\backslash\{C_{\ast i}\}\subset\mathcal{B}$.
We pick $C^{\infty}$ functions or mollifiers $\rho_{{\rm UV},i}$ or simply $\rho_{i}$ such that ${\rm supp}\ \rho_i\in\bar{B}(0,1)$,
$\rho_i\ge 0$, and $\int_{\mathbb{R}^d}dx\ \rho_i(x)=1$. We pick shifts $\Delta r_i\in\{0,1\}$ which are needed for the telescopic sum argument and bound on $||M_{r}(f)-M_{r-1}(f)||_{L^p}$ in \S\ref{finishproof}.
We also pick $n\ge m$ and $A_{m+1},\ldots,A_n\in\mathcal{A}$ for the spectator fields and we
fix a collection of test functions $f_1,\ldots,f_n$ in $\mathcal{S}(\mathbb{R}^d)$.
There will only be one varying quantity in the following discussion, namely the UV cut-off $r\in\mathbb{Z}$ which will be taken to $-\infty$.
We define $r_i=r-\Delta r_i$ and the rescaled mollifiers $\rho_{i,r_i}(x)=L^{-dr_i}\rho_{i}(L^{-r_i}x)$.
We set
\[
Z_{i,r_i}(x)=\left\{
\int_{{\rm Conf}_2} dy\ dz\ \rho_{i,r_i}(x-y)\rho_{i,r_i}(x-z)
\mathcal{C}_{A_i B_i}^{C_{\ast i}}(y,z)
\right\}^{-1}
\]
as well as $\mathcal{O}_{A_i,r_i}(x)=(\mathcal{O}_{A_i}\ast \rho_{i,r_i})(x)=\langle\mathcal{O}_{A_i}(y), \rho_{i,r_i}(x-y)\rangle_{y}$.
We likewise set
$\mathcal{O}_{B_i,r_i}(x)=
(\mathcal{O}_{B_i}\ast \rho_{i,r_i})(x)=\langle\mathcal{O}_{B_i}(z), \rho_{i,r_i}(x-z)\rangle_{z}$.
For $C\in \mathcal{A}(\Delta_i)\backslash\{C_{\ast i}\}$, we let
\[
g_{i,r_i}(x,z)=\rho_{i,r_i}(x-z)\int_{\mathbb{R}^d\backslash\{z\}}dy\ \rho_{i,r_i}(x-y)
\ \mathcal{C}_{A_i B_i}^{C}(y,z)
\]
and define
$\widetilde{\mathcal{O}}_{C,r_i}(x)=\langle\mathcal{O}_{C}(z), g_{i,r_i}(x-z)\rangle_{z}$
as before.
This gives us the candidate for the regularized product ``$\mathcal{O}_{C_{\ast i}}$'' as a function of $x$.
Namely, it is
\[
M_{i,r_i}(x)=Z_{i,r_i}(x)\left[
\mathcal{O}_{A_i,r_i}(x)\mathcal{O}_{B_i,r_i}(x)-
\sum_{C\in\mathcal{A}(\Delta_i)\backslash\{C_{\ast i}\}}
\widetilde{\mathcal{O}}_{C,r_i}(x)
\right]\ .
\]
We let
$M_{i,r_i}(f_i)=\langle M_{i,r_i}(x),f_i(x)\rangle_{x}=\int_{\mathbb{R}^d} dx\ M_{i,r_i}(x)f_i(x)$.
Our goal is to estimate the ``moment'' difference
$\Upsilon_r={\rm TM}_r-{\rm IPC}$
where
\[
{\rm TM}_r=\mathbb{E}\left[
\prod_{i=1}^{m} M_{i,r_i}(f_i)\times
\prod_{i=m+1}^{n} \mathcal{O}_{A_i}(f_i)
\right]
\]
is a true moment, while
\[
{\rm IPC}=\int_{{\rm Conf}_n}\prod_{i=1}^{n}dx_i\ \prod_{i=1}^{n} f_i(x_i)\times
\left\langle \prod_{i=1}^{m} \mathcal{O}_{C_{\ast i}}(x_i)\times
\prod_{i=m+1}^{n} \mathcal{O}_{A_i}(x_i)
\right\rangle
\]
is an integral of pointwise correlations.

We have seen that $M_{i,r_i}(f_i)$ is a random variable in $L^p(\Omega,\mathcal{F},\mathbb{P})$ for all $p\ge 1$.
By the multilinear H\"{o}lder inequality, the random variable in the expectation defining ${\rm TM}_r$ is integrable.
It is not hard to see that as a consequence of \S\ref{pinsum}, the representations (\ref{ocrep}) and (\ref{oaobrep})
for the pieces making up
$M_{i,r_i}(f_i)$ and the discussion in~\cite[\S{II.2.5}]{Fernique},
that one has the following pointwise representation
\[
{\rm TM}_r=\int_{{\rm Conf}_{m+n}}\prod_{i=1}^{m}dy_i\ \prod_{i=1}^{n}dz_i\ 
\prod_{i=1}^{m}h_{i,r_i}(y_i,z_i)\times\prod_{i=m+1}^{n}f_i(z_i) 
\]
\[
\times
\left\langle \prod_{i=1}^{m} P_i(y_i,z_i)
\prod_{i=m+1}^{n} \mathcal{O}_{A_i}(z_i)
\right\rangle
\]
where, for $1\le i\le m$, $P_i\in\mathcal{V}_2$ is given by
\[
P_i=\mathcal{O}_{A_i}\otimes\mathcal{O}_{B_i}
-\sum_{C_i\in \mathcal{A}(\Delta_i)\backslash\{C_{\ast i}\}}
\mathcal{C}_{A_i B_i}^{C_i} \mathcal{O}_{\bbone}\otimes\mathcal{O}_{C_i}
\]
and where
\begin{equation}
h_{i,r_i}(y_i,z_i)=\int_{\mathbb{R}^d}dx_i\ Z_{i,r_i}(x_i)
f_i(x_i)\rho_{i,r_i}(x_i-y_i)\rho_{i,r_i}(x_i-z_i)\ .
\label{hintegral}
\end{equation}

The proof of our main theorem from \S\ref{resultssec}
is based on the follow proposition which is a purely combinatorial estimate.

\begin{Proposition}\label{combiprop}
There exists $\nu>0$ such that$
|\Upsilon_r|\le O(1)\ L^{\nu r}$, for all $r$.
\end{Proposition}

The proof of this proposition is provided in the next section. It can be viewed as an amplification of the one given
in \S\ref{warmupsec}.

\section{The main estimate}\label{mainestimatesec}

\subsection{Preparatory steps}

We define $Q_i\in\mathcal{V}_2$ by
$Q_i=\mathcal{C}_{A_i B_i}^{C_{\ast i}} \mathcal{O}_{\bbone}\otimes\mathcal{O}_{C_{\ast i}}$
so that 
\begin{equation}
P_i=Q_i+R_{1,i}
\label{pqr}
\end{equation}
where $R_{1,i}$ is the OPE-like element
\[
R_{1,i}=\mathcal{O}_{A_i}\otimes\mathcal{O}_{B_i}
-\sum_{C_i\in\mathcal{A}(\Delta_i)} \mathcal{C}_{A_i B_i}^{C_i} \mathcal{O}_{\bbone}\otimes\mathcal{O}_{C_i}\ .
\]
By a decomposition $(I_1,\ldots,I_p)$ of a finite set $I$
we mean an ordered collection of disjoint subsets whose union is $I$.
This differs from a set partition because of the ordering and allowing the empty set.
Borrowing our notation from the theory of symmetric functions, we will write
$(I_1,\ldots,I_p)\vdash I$ in order to say that $(I_1,\ldots,I_p)$ is a decomposition of $I$.

We now expand using (\ref{pqr}) so that
\[
{\rm TM}_r=
\sum_{(I_1,I_{23})\vdash [m]}
\int_{{\rm Conf}_{m+n}}\prod_{i=1}^{m}dy_i\ \prod_{i=1}^{n}dz_i\ 
\prod_{i=1}^{m}h_{i,r_i}(y_i,z_i)\times\prod_{i=m+1}^{n}f_i(z_i) 
\]
\[
\times
\left\langle \prod_{i\in I_1} R_{1,i}(y_i,z_i)
\times \prod_{i\in I_{23}} Q_{i}(y_i,z_i)
\times \prod_{i\in I_4} \mathcal{O}_{A_i}(z_i)
\right\rangle
\]
where $I_4=[n]\backslash [m]$ is fixed.

We replace $h_{i,r_i}(y_i,z_i)$, for $i\in I_{23}$,
by the integral in (\ref{hintegral}) which introduces $|I_{23}|$ new variables of integration $x_i$.
We use the forgetful property (\ref{forgetful})
to replace $Q_i(y_i,z_i)$ by $\widetilde{Q}_i(x_i,y_i,z_i)$
where $\widetilde{Q}_i\in\mathcal{V}_3$ is given by
$\mathcal{O}_{\bbone}\otimes Q_i=(1\otimes \mathcal{C}_{A_i B_i}^{C_{\ast i}})
\ \mathcal{O}_{\bbone}\otimes\mathcal{O}_{\bbone}\otimes\mathcal{O}_{C_{\ast i}}$
where ``$1$'' simply is the function of the first argument $x_i$ which is constant and equal to one.
We define $\widehat{Q}_i=(1\otimes \mathcal{C}_{A_i B_i}^{C_{\ast i}})
\ \mathcal{O}_{C_{\ast i}}\otimes\mathcal{O}_{\bbone}\otimes\mathcal{O}_{\bbone}$ so that
\begin{equation}
\widetilde{Q}_i=\widehat{Q}_i+R_{2,i}
\label{Qexp}
\end{equation}
with 
$R_{2,i}=(1\otimes \mathcal{C}_{A_i B_i}^{C_{\ast i}})
\left[
\mathcal{O}_{\bbone}\otimes\mathcal{O}_{\bbone}\otimes\mathcal{O}_{C_{\ast i}}-
\mathcal{O}_{C_{\ast i}}\otimes\mathcal{O}_{\bbone}\otimes\mathcal{O}_{\bbone}
\right]$.
We expand using (\ref{Qexp}) and get
\[
{\rm TM}_r=
\sum_{(I_1,I_2,I_3)\vdash [m]}
\int_{{\rm Conf}_{2|I_1|+3|I_2|+3|I_3|+|I_4|}}\ 
\prod_{i\in I_2\cup I_3}dx_i\ 
\prod_{i\in [m]}dy_i\ \prod_{i\in [n]}dz_i\ 
\]
\[
\prod_{i\in I_1}h_{i,r_i}(y_i,z_i)\times\prod_{i\in I_4}f_i(z_i) 
\times
\prod_{i\in I_2\cup I_3}
\left[
Z_{i,r_i}(x_i)
f_i(x_i)\rho_{i,r_i}(x_i-y_i)\rho_{i,r_i}(x_i-z_i)
\right]
\]
\[
\times
\left\langle \prod_{i\in I_1} R_{1,i}(y_i,z_i)
\times \prod_{i\in I_2} R_{2,i}(x_i,y_i,z_i)
\times \prod_{i\in I_3} \widehat{Q}_{i}(x_i,y_i,z_i)
\times \prod_{i\in I_4} \mathcal{O}_{A_i}(z_i)
\right\rangle\ .
\]

Since $\widehat{Q}_{i}(x_i,y_i,z_i)=
\mathcal{C}_{A_i B_i}^{C_{\ast i}}(y_i,z_i)
\mathcal{O}_{C_{\ast i}}(x_i)$
``inside correlations'', we can factor $\mathcal{C}_{A_i B_i}^{C_{\ast i}}(y_i,z_i)$
out of the pointwise correlation and integrate over $y_i$, $z_i$ for $i\in I_3$.
This produces the inverse of $Z_{i,r_i}(x_i)$ by definition of the latter.
Thus,
\[
{\rm TM}_r=
\sum_{(I_1,I_2,I_3)\vdash [m]}
\int_{{\rm Conf}_{2|I_1|+3|I_2|+|I_3|+|I_4|}}\ 
\prod_{i\in I_2\cup I_3}dx_i\ 
\prod_{i\in I_1\cup I_2}dy_i\ \prod_{i\in I_1\cup I_2\cup I_4}dz_i\ 
\]
\[
\prod_{i\in I_1}h_{i,r_i}(y_i,z_i)\times\prod_{i\in I_4}f_i(z_i) 
\times\prod_{i\in I_3}f_i(x_i) 
\times\prod_{i\in I_2}
\left[
Z_{i,r_i}(x_i)
f_i(x_i)\rho_{i,r_i}(x_i-y_i)\rho_{i,r_i}(x_i-z_i)
\right]
\]
\[
\times
\left\langle \prod_{i\in I_1} R_{1,i}(y_i,z_i)
\times \prod_{i\in I_2} R_{2,i}(x_i,y_i,z_i)
\times \prod_{i\in I_3} \mathcal{O}_{C_{\ast i}}(x_i)
\times \prod_{i\in I_4} \mathcal{O}_{A_i}(z_i)
\right\rangle\ .
\]
We introduce the CZ-like elements $T_i=\mathcal{O}_{C_{\ast i}}\otimes\mathcal{O}_{\bbone}-\mathcal{O}_{\bbone}
\otimes\mathcal{O}_{C_{\ast i}}$ for $i\in I_2$.
Noting that one can write
\[
R_{2,i}(x_i,y_i,z_i)=\mathcal{C}_{A_i B_i}^{C_{\ast i}}(y_i,z_i)\left[
\mathcal{O}_{C_{\ast i}}(z_i)-\mathcal{O}_{C_{\ast i}}(x_i)
\right]
=\mathcal{C}_{A_i B_i}^{C_{\ast i}}(y_i,z_i) T_i(z_i,x_i)
\]
``inside correlations'', we factor the $\mathcal{C}_{A_i B_i}^{C_{\ast i}}(y_i,z_i)$'s
out of the pointwise correlation. We also
insert the integrals (\ref{hintegral}) and therefore create $|I_1|$ new variables of integration $x_i$ for $i\in I_1$.
Therefore
\[
{\rm TM}_r=
\sum_{(I_1,I_2,I_3)\vdash [m]}
\int_{{\rm Conf}_{3|I_1|+3|I_2|+|I_3|+|I_4|}}\ 
\prod_{i\in I_1\cup I_2\cup I_3}dx_i\ 
\prod_{i\in I_1\cup I_2}dy_i\ \prod_{i\in I_1\cup I_2\cup I_4}dz_i\ 
\]
\[
\prod_{i\in I_2}\mathcal{C}_{A_i B_i}^{C_{\ast i}}(y_i,z_i)
\times\prod_{i\in I_4}f_i(z_i) 
\times\prod_{i\in I_3}f_i(x_i) 
\times\prod_{i\in I_1\cup I_2}
\left[
Z_{i,r_i}(x_i)
f_i(x_i)\rho_{i,r_i}(x_i-y_i)\rho_{i,r_i}(x_i-z_i)
\right]
\]
\[
\times
\left\langle \prod_{i\in I_1} R_{1,i}(y_i,z_i)
\times \prod_{i\in I_2} T_i(z_i,x_i)
\times \prod_{i\in I_3} \mathcal{O}_{C_{\ast i}}(x_i)
\times \prod_{i\in I_4} \mathcal{O}_{A_i}(z_i)
\right\rangle\ .
\]
Now note that ${\rm IPC}$ is the term corresponding to the decomposition $(I_1,I_2,I_3)=(\emptyset,\emptyset,[m])$.
We now start putting
absolute values inside the integral in order to write estimates.
Hence
\[
|\Upsilon_r|\le \sum_{\substack{(I_1,I_2,I_3)\vdash [m] \\ (I_1,I_2,I_3)\neq (\emptyset,\emptyset,[m])}}
\mathcal{P}^{({\rm I})} \mathcal{I}^{({\rm I})}
\]
where $\mathcal{P}^{({\rm I})}=1$ and  
\[
\mathcal{I}^{({\rm I})}
=
\int_{{\rm Conf}_{3|I_1|+3|I_2|+|I_3|+|I_4|}}\ 
\prod_{i\in I_1\cup I_2\cup I_3}dx_i\ 
\prod_{i\in I_1\cup I_2}dy_i\ \prod_{i\in I_1\cup I_2\cup I_4}dz_i\ 
\]
\[
\prod_{i\in I_2}|\mathcal{C}_{A_i B_i}^{C_{\ast i}}(y_i,z_i)|
\times\prod_{i\in I_4}|f_i(z_i)| 
\times\prod_{i\in I_3}|f_i(x_i)| 
\times\prod_{i\in I_1\cup I_2}
\left|
Z_{i,r_i}(x_i)
f_i(x_i)\rho_{i,r_i}(x_i-y_i)\rho_{i,r_i}(x_i-z_i)
\right|
\]
\[
\times\left|
\left\langle \prod_{i\in I_1} R_{1,i}(y_i,z_i)
\times \prod_{i\in I_2} T_i(z_i,x_i)
\times \prod_{i\in I_3} \mathcal{O}_{C_{\ast i}}(x_i)
\times \prod_{i\in I_4} \mathcal{O}_{A_i}(z_i)
\right\rangle\right|
\]
and we suppressed the dependence on the decomposition from the notation for the prefactor $\mathcal{P}^{({\rm I})}$ and the integral
$\mathcal{I}^{({\rm I})}$.
There will be many more later, hence the roman numerals (I), (II), etc.
We use the bound
\[
|\rho_{i,r_i}(x_i-y_i)|\le O(1) L^{-dr}\bbone\{|x_i-y_i|\le L^r\}
\]
and from now on we will denote $r$-independent constants by $O(1)$ if convenient. These constants can depend on everything else.
For instance, the one above already ate up an $L^{d}$ factor which is needed if $\Delta r_i=1$.

In \S\ref{probatocomb}, we already established
\[
Z_{i,r_i}(x)\le O(1)\langle x\rangle^{2k} L^{r([A_i]+[B_i]-[C_{\ast i}]-\epsilon)}\ .
\]
Inserting the previous bounds for the $Z_{i,r_i}$ and
the $\rho_{i,r_i}$ as well as (\ref{hardCbd}) for the $\mathcal{C}_{A_i B_i}^{C_{\ast i}}$
we get
\[
|\Upsilon_r|\le O(1)\sum_{\substack{(I_1,I_2,I_3)\vdash [m] \\ (I_1,I_2,I_3)\neq (\emptyset,\emptyset,[m])}}
\mathcal{P}^{({\rm II})} \mathcal{I}^{({\rm II})}
\]
\begin{flalign*}
&{\rm with}\ \ & &
\mathcal{P}^{({\rm II})}=\prod_{i\in I_1\cup I_2}
L^{r([A_i]+[B_i]-\Delta_i-2d-\epsilon)} & & \ \\
&{\rm and}\ \ & &
\mathcal{I}^{({\rm II})}
=
\int_{{\rm Conf}_{3|I_1|+3|I_2|+|I_3|+|I_4|}}\ 
\prod_{i\in I_1\cup I_2\cup I_3}dx_i\ 
\prod_{i\in I_1\cup I_2}dy_i\ \prod_{i\in I_1\cup I_2\cup I_4}dz_i & & \  
\end{flalign*}
\[
\prod_{i\in I_1\cup I_2}\langle x_i\rangle^{2k}\times 
\prod_{i\in I_2}\langle y_i\rangle^{k}\times \prod_{i\in I_2}\langle z_i\rangle^{k}\times
\prod_{i\in I_1\cup I_2\cup I_3}|f_i(x_i)|\times \prod_{i\in I_4}|f_i(z_i)|
\]
\[
\times\prod_{i\in I_1\cup I_2}\left(
\ \bbone\{|x_i-y_i|\le L^r\}\ \bbone\{|x_i-z_i|\le L^r\}\ 
\right)
\times\prod_{i\in I_2}|y_i-z_i|^{-[A_i]-[B_i]+\Delta_i-\epsilon}
\]
\[
\times\left|
\left\langle \prod_{i\in I_1} R_{1,i}(y_i,z_i)
\times \prod_{i\in I_2} T_i(z_i,x_i)
\times \prod_{i\in I_3} \mathcal{O}_{C_{\ast i}}(x_i)
\times \prod_{i\in I_4} \mathcal{O}_{A_i}(z_i)
\right\rangle\right|\ .
\]
Now we make sure every point gets some long-distance decay by tapping into the nearest test function.
Since $r\le 0$ and $|u|\le 1$ implies $\langle x-u\rangle\le \sqrt{3}\langle x\rangle$ we can insert for all $i\in I_1$,
\begin{flalign*}
&{\ }\ \ & &
1=\langle y_i\rangle^{-(d+1+k)}\langle y_i\rangle^{d+1+k}\le O(1)\langle y_i\rangle^{-(d+1+k)}\langle x_i\rangle^{d+1+k} & & \ \\
&{\rm and}\ \ & & 
1=\langle z_i\rangle^{-(d+1+2k)}\langle z_i\rangle^{d+1+2k}\le O(1)\langle z_i\rangle^{-(d+1+2k)}\langle x_i\rangle^{d+1+2k}\ . & & \ 
\end{flalign*}
Likewise, for $i\in I_2$, we insert
\begin{flalign*}
&{\ }\ \ & & 
1=\langle y_i\rangle^{-k}\langle y_i\rangle^{k}\le O(1)\langle y_i\rangle^{-k}\langle x_i\rangle^{k} & & \ \\
&{\rm and}\ \ & &
1=\langle z_i\rangle^{-(d+1+2k)}\langle z_i\rangle^{d+1+2k}\le O(1)\langle z_i\rangle^{-(d+1+2k)}\langle x_i\rangle^{d+1+2k}\ . & & \ 
\end{flalign*}
Therefore
\[
\mathcal{I}^{({\rm II})}
\le O(1)
\int_{{\rm Conf}_{3|I_1|+3|I_2|+|I_3|+|I_4|}}\ 
\prod_{i\in I_1\cup I_2\cup I_3}dx_i\ 
\prod_{i\in I_1\cup I_2}dy_i\ \prod_{i\in I_1\cup I_2\cup I_4}dz_i\ 
\]
\[
\prod_{i\in I_1}\langle x_i\rangle^{2d+2+5k}\times \prod_{i\in I_2}\langle x_i\rangle^{d+1+5k}\times 
\prod_{i\in I_1}\langle y_i\rangle^{-(d+1+k)}\times \prod_{i\in I_1}\langle z_i\rangle^{-(d+1+2k)}
\]
\[
\times \prod_{i\in I_2}\langle z_i\rangle^{-(d+1+k)}\times
\prod_{i\in I_1\cup I_2\cup I_3}|f_i(x_i)|\times \prod_{i\in I_4}|f_i(z_i)|
\]
\[
\times\prod_{i\in I_1\cup I_2}\left(
\ \bbone\{|x_i-y_i|\le L^r\}\ \bbone\{|x_i-z_i|\le L^r\}\ 
\right)
\times\prod_{i\in I_2}|y_i-z_i|^{-[A_i]-[B_i]+\Delta_i-\epsilon}
\]
\[
\times\left|
\left\langle \prod_{i\in I_1} R_{1,i}(y_i,z_i)
\times \prod_{i\in I_2} T_i(z_i,x_i)
\times \prod_{i\in I_3} \mathcal{O}_{C_{\ast i}}(x_i)
\times \prod_{i\in I_4} \mathcal{O}_{A_i}(z_i)
\right\rangle\right|\ .
\]
We now bound the test functions by suitable seminorms and prepare for the integration over points which are not in the big
pointwise correlation, namely, $x_i$ for $i\in I_1$ and $y_i$ for $i\in I_2$.

For $i\in I_1$, we use the inequality
\[
\langle x_i \rangle^{2d+2+5k}|f_i(x_i)|\ \bbone\{|x_i-y_i|\le L^r\}\ \bbone\{|x_i-z_i|\le L^r\}
\le 
\]
\[
||f_i||_{0,2d+2+5k}\ \bbone\{|x_i-z_i|\le L^r\}\ \bbone\{|y_i-z_i|\le 2L^r\}\ .
\]
For $i\in I_2$, we use the inequality
\[
\langle x_i \rangle^{d+1+5k}
|f_i(x_i)|\ \bbone\{|x_i-y_i|\le L^r\}\ \bbone\{|x_i-z_i|\le L^r\}
\le 
\]
\[
\langle x_i \rangle^{-(d+1+2k)}
||f_i||_{0,2d+2+7k}\ \bbone\{|y_i-z_i|\le 2L^r\}\ \bbone\{|x_i-z_i|\le 2L^r\}\ .
\]
Note that we allowed ourselves to loose a bit on the $|x_i-z_i|$ bound in order to make this case look like the previous one
and thus ease the bookkeeping.
For $i\in I_3$, we use the inequality
\[
|f_i(x_i)|\le\langle x_i \rangle^{-(d+1+2k)}
||f_i||_{0,d+1+2k}\ .
\]
Finally, for 
$i\in I_4$, we use the inequality
\[
|f_i(z_i)|\le\langle z_i \rangle^{-(d+1+2k)}
||f_i||_{0,d+1+2k}\ .
\]
Absorbing the Schwartz seminorms into the $O(1)$ constant, we obtain
\[
\mathcal{I}^{({\rm II})}
\le O(1)
\int_{{\rm Conf}_{3|I_1|+3|I_2|+|I_3|+|I_4|}}\ 
\prod_{i\in I_1\cup I_2\cup I_3}dx_i\ 
\prod_{i\in I_1\cup I_2}dy_i\ \prod_{i\in I_1\cup I_2\cup I_4}dz_i\ 
\]
\[ 
\prod_{i\in I_1}\langle y_i\rangle^{-(d+1+k)}\times 
\prod_{i\in I_1}\langle z_i\rangle^{-(d+1+2k)}\times
\prod_{i\in I_2}\langle x_i\rangle^{-(d+1+2k)}
\times\prod_{i\in I_2}\langle z_i\rangle^{-(d+1+k)}
\]
\[
\times
\prod_{i\in I_3}\langle x_i\rangle^{-(d+1+2k)}\times 
\prod_{i\in I_4}\langle z_i\rangle^{-(d+1+2k)}
\times\prod_{i\in I_1}\left(
\ \bbone\{|x_i-z_i|\le L^r\}\ \bbone\{|y_i-z_i|\le 2L^r\}\ 
\right)
\]
\[
\times\prod_{i\in I_2}\left(
\ \bbone\{|y_i-z_i|\le 2L^r\}\ \bbone\{|x_i-z_i|\le 2L^r\}\ 
\right)
\times\prod_{i\in I_2}|y_i-z_i|^{-[A_i]-[B_i]+\Delta_i-\epsilon}
\]
\[
\times\left|
\left\langle \prod_{i\in I_1} R_{1,i}(y_i,z_i)
\times \prod_{i\in I_2} T_i(z_i,x_i)
\times \prod_{i\in I_3} \mathcal{O}_{C_{\ast i}}(x_i)
\times \prod_{i\in I_4} \mathcal{O}_{A_i}(z_i)
\right\rangle\right|\ .
\]
We now integrate over $x_i$, $i\in I_1$ and $y_i$, $i\in I_2$, respectively using the $\alpha=0$ and
$\alpha=[A_i]+[B_i]-\Delta_i+\epsilon$ cases of Lemma \ref{localL1bd}.
Thus
\[
|\Upsilon_r|\le O(1)\sum_{\substack{(I_1,I_2,I_3)\vdash [m] \\ (I_1,I_2,I_3)\neq (\emptyset,\emptyset,[m])}}
\mathcal{P}^{({\rm III})} \mathcal{I}^{({\rm III})}
\]
\begin{flalign*}
&{\rm with}\ \ & & 
\mathcal{P}^{({\rm III})}=\prod_{i\in I_1} L^{r([A_i]+[B_i]-\Delta_i-d-\epsilon)}\times
\prod_{i\in I_2} L^{-r(d+2\epsilon)} & & \ \\
&{\rm and}\ \ & &
\mathcal{I}^{({\rm III})}=
\int_{{\rm Conf}_{2|I_1|+2|I_2|+|I_3|+|I_4|}}\ 
\prod_{i\in I_2\cup I_3}dx_i\ 
\prod_{i\in I_1}dy_i\ \prod_{i\in I_1\cup I_2\cup I_4}dz_i & & \ 
\end{flalign*}
\[ 
\prod_{i\in I_1}\langle y_i\rangle^{-(d+1+k)}\times 
\prod_{i\in I_1}\langle z_i\rangle^{-(d+1+2k)}\times
\prod_{i\in I_2}\langle x_i\rangle^{-(d+1+2k)}
\times 
\prod_{i\in I_2}\langle z_i\rangle^{-(d+1+k)}
\]
\[
\times
\prod_{i\in I_3}\langle x_i\rangle^{-(d+1+2k)}\times 
\prod_{i\in I_4}\langle z_i\rangle^{-(d+1+2k)}
\times\prod_{i\in I_1}
\ \bbone\{|y_i-z_i|\le 2L^r\}\ 
\times\prod_{i\in I_2}
\ \bbone\{|x_i-z_i|\le 2L^r\}\ 
\]
\[
\times\left|
\left\langle \prod_{i\in I_1} R_{1,i}(y_i,z_i)
\times \prod_{i\in I_2} T_i(z_i,x_i)
\times \prod_{i\in I_3} \mathcal{O}_{C_{\ast i}}(x_i)
\times \prod_{i\in I_4} \mathcal{O}_{A_i}(z_i)
\right\rangle\right|\ .
\]
In order to continue, we will consolidate our notation by renaming the dummy variables of integration
as follows
\begin{itemize}
\item
For all $i\in I_1$, $y_i$ stays $y_i$.
\item
For all $i\in I_1$, $z_i$ becomes $x_i$.
\item
For all $i\in I_2$, $z_i$ becomes $y_i$.
\item
For all $i\in I_2$, $x_i$ stays $x_i$.
\item
For all $i\in I_3$, $x_i$ stays $x_i$.
\item
For all $i\in I_4$, $z_i$ becomes $x_i$.
\end{itemize}
We also introduce the notation $I_{34}=I_3\cup I_4$ together with the new field labels $D_i$, $i\in I_{34}$
defined as follows.
\begin{itemize}
\item
For all $i\in I_3$, $D_i=C_{\ast i}$ .
\item
For all $i\in I_4$, $D_i=A_i$.
\end{itemize}
These changes made, we can give a simpler
formula for the integral, i.e.,
\[
\mathcal{I}^{({\rm III})}=
\int_{{\rm Conf}_{2|I_1|+2|I_2|+|I_3|+|I_4|}}\ 
\prod_{i\in [n]}dx_i\ 
\prod_{i\in I_1\cup I_2}dy_i\ 
\]
\[ 
\prod_{i\in [n]}\langle x_i\rangle^{-(d+1+2k)}\times 
\prod_{i\in I_1\cup I_2}\langle y_i\rangle^{-(d+1+k)}
\times\prod_{i\in I_1\cup I_2}
\ \bbone\{|y_i-x_i|\le 2L^r\}\ 
\]
\[
\times\left|
\left\langle \prod_{i\in I_1} R_{1,i}(y_i,x_i)
\times \prod_{i\in I_2} T_i(y_i,x_i)
\times \prod_{i\in I_{34}} \mathcal{O}_{D_i}(x_i)
\right\rangle\right|\ .
\]

\subsection{Re-expansion}

At this point the pointwise correlation looks ready for the EFNNB.
However, we are missing the indicator functions needed for this bound.
We pick a  number $\delta>0$.
For each $i$, $1\le i\le n$, we insert in the integral $\mathcal{I}^{({\rm III})}$ the identity
\[
1=\bbone\left\{
\min_{j\neq i} |x_i-x_j|>\delta L^{r}
\right\}+\bbone\left\{
\min_{j\neq i} |x_i-x_j|\le\delta L^{r}
\right\}\ .
\]
Then we expand. This results in a sum over decompositions $(I_{\rm G},I_{\rm B})$ of $[n]$.
Let $i\in[n]$. If we pick for it the first indicator function we say that $i$ is good.
If we pick the second, we say that $i$ is bad. The set $I_{\rm G}$ is that of good labels, whereas
$I_{\rm B}$ is the set of bad ones.
The $R_{1,i}(y_i,x_i)$'s and $T_i(y_i,x_i)$ with $i\in I_{\rm G}$ are left untouched.
On the other hand, those for $i\in I_{\rm B}$ must re-expanded completely.
Recall that ``inside correlations'' we have
\begin{equation}
R_{1,i}(y_i,x_i)=\mathcal{O}_{A_i}(y_i)\mathcal{O}_{B_i}(x_i)
-\sum_{C_i\in\mathcal{A}(\Delta_i)} \mathcal{C}_{A_i B_i}^{C_i}(y_i,x_i)\mathcal{O}_{C_i}(x_i)
\label{opereexp}
\end{equation}
and
$T_i(y_i,x_i)=\mathcal{O}_{C_{\ast i}}(y_i)-\mathcal{O}_{C_{\ast i}}(x_i)$.

We need new notation for subsets of $[n]$, as follows.
\begin{itemize}
\item
We let $I_{1{\rm G}}=I_1\cap I_{\rm G}$.
\item
We let $I_{2{\rm G}}=I_2\cap I_{\rm G}$.
\item
We let $I_{1{\rm BOO}}$ be the set of $i$'s in $I_1\cap I_{\rm B}$ for which the first $\mathcal{O}\mathcal{O}$
term in (\ref{opereexp}) is chosen in the expansion.
\item
We let $I_{1{\rm BCO}}$ be the set of $i$'s in $I_1\cap I_{\rm B}$ for which a $\mathcal{C}\mathcal{O}$
term in the sum in (\ref{opereexp}) is chosen for the expansion.
\item
We let $I_{2{\rm BY}}$ be the set of $i$'s in $I_2\cap I_{\rm B}$ for which $\mathcal{O}_{C_{\ast i}}(y_i)$
is chosen.
\item
We let $I_{2{\rm BX}}$ be the set of $i$'s in $I_2\cap I_{\rm B}$ for which $\mathcal{O}_{C_{\ast i}}(x_i)$
is chosen.
\end{itemize}
After putting absolute values so that one does not have to worry about the signs produced by $i$'s in $I_{1{\rm BCO}}$ or $I_{2{\rm BX}}$,
and after factoring the $\mathcal{C}$'s out of the pointwise correlation,
the resulting estimate is
\[
\mathcal{I}^{({\rm III})}\le \sum 
\int_{{\rm Conf}_{2|I_1|+2|I_2|+|I_3|+|I_4|}}\ 
\prod_{i\in [n]}dx_i\ 
\prod_{i\in I_1\cup I_2}dy_i\ 
\]
\[ 
\prod_{i\in [n]}\langle x_i\rangle^{-(d+1+2k)}\times 
\prod_{i\in I_1\cup I_2}\langle y_i\rangle^{-(d+1+k)}
\times\prod_{i\in I_1\cup I_2}
\ \bbone\{|y_i-x_i|\le 2L^r\}\ 
\]
\[
\times
\prod_{i\in I_{\rm G}}
\bbone\left\{
\min_{j\neq i} |x_i-x_j|>\delta L^{r}
\right\}
\times
\prod_{i\in I_{\rm B}}
\bbone\left\{
\min_{j\neq i} |x_i-x_j|\le\delta L^{r}
\right\}
\]
\[
\times\prod_{i\in I_{1{\rm BCO}}} \left|\mathcal{C}_{A_i B_i}^{C_i}(y_i,x_i)\right|
\times\left|
\left\langle \prod_{i\in I_{1{\rm G}}} R_{1,i}(y_i,x_i)
\times \prod_{i\in I_{2{\rm G}}} T_i(y_i,x_i)
\times \prod_{i\in I_{1{\rm BOO}}} \left[\mathcal{O}_{A_i}(y_i)\mathcal{O}_{B_i}(x_i)\right]
\right.\right.
\]
\[
\left.\left.
\times \prod_{i\in I_{1{\rm BCO}}} \mathcal{O}_{C_i}(x_i)
\times \prod_{i\in I_{2{\rm BY}}} \mathcal{O}_{C_{\ast i}}(y_i)
\times \prod_{i\in I_{2{\rm BX}}} \mathcal{O}_{C_{\ast i}}(x_i)
\times \prod_{i\in I_{34}} \mathcal{O}_{D_i}(x_i)
\right\rangle\right|\ .
\]
For obvious reasons, we did not write the (formidable) summation index under the sum.
Indeed, the sum is now over decompositions $(I_{1{\rm G}},I_{1{\rm BOO}},I_{1{\rm BCO}})$
of $I_1$, as well as decompositions $(I_{2{\rm G}},I_{2{\rm BY}},I_{2{\rm BX}})$,
of $I_2$, followed by a summation over $C_i\in\mathcal{A}(\Delta_i)$, for each $i\in I_{1{\rm BCO}}$.

At (long) last, we are now able to use the EFNNB as well as the bound (\ref{hardCbd}) for the $\mathcal{C}$'s.
Indeed, if we pick $\delta\ge 4\eta^{-1}$ then the indicator functions present
imply the needed repulsive condition for the OPE-like elements $R_{1,i}(y_i,x_i)$, $i\in I_{1{\rm G}}$, and
CZ-like elements $T_{i}(y_i,x_i)$, $i\in I_{2{\rm G}}$.
Since the $x_i$'s and $y_i$'s are treated differently, and in order to keep the size of formulas under control,
we need to introduce two notions of nearest-neighbor distance, before writing the outcome of the EFNNB.

For all $i\in [n]\backslash I_{2{\rm BY}}$, we let
\[
{\rm NNDX}_i=\min\left\{
\min_{j\in ([n]\backslash I_{2{\rm BY}})\backslash\{i\}}|x_i-x_j|\ ,\ \min_{j\in I_{1{\rm BOO}}\cup I_{2{\rm BY}} }|x_i-y_j|
\right\}\ .
\]
For all $i\in I_{1{\rm BOO}}\cup I_{2{\rm BY}}$, we let
\[
{\rm NNDY}_i=\min\left\{
\min_{j\in [n]\backslash I_{2{\rm BY}}}|y_i-x_j|\ ,\ \min_{j\in (I_{1{\rm BOO}}\cup I_{2{\rm BY}})\backslash\{i\} }|y_i-y_j|
\right\}\ .
\]

These precautions taken, we now have
\[
|\Upsilon_r|\le O(1)\sum_{\substack{(I_1,I_2,I_3)\vdash [m] \\ (I_1,I_2,I_3)\neq (\emptyset,\emptyset,[m]) }}
\sum_{\substack{ (I_{1{\rm G}},I_{1{\rm BOO}},I_{1{\rm BCO}})\vdash I_1 \\ (I_{2{\rm G}},I_{2{\rm BY}},I_{2{\rm BX}})\vdash I_2 }}
\ \ \sum_{(C_i)_{i\in I_{1{\rm BCO}}}\in \prod_{i\in I_{1{\rm BCO}}}\mathcal{A}(\Delta_i)}
\mathcal{P}^{({\rm IV})} \mathcal{I}^{({\rm IV})}
\]
\begin{flalign*}
&{\rm with}\ \ & &
\mathcal{P}^{({\rm IV})}=\mathcal{P}^{({\rm III})}=\prod_{i\in I_1} L^{r([A_i]+[B_i]-\Delta_i-d-\epsilon)}\times
\prod_{i\in I_2} L^{-r(d+2\epsilon)} & & \ \\
&{\rm and}\ \ & &
\mathcal{I}^{({\rm IV})}=
\int_{{\rm Conf}_{2|I_1|+2|I_2|+|I_3|+|I_4|}}\ 
\prod_{i\in [n]}dx_i\ 
\prod_{i\in I_1\cup I_2}dy_i & & \  
\end{flalign*}
\[ 
\prod_{i\in I_{1{\rm BCO}}}\langle x_i\rangle^{-(d+1)}
\times\prod_{i\in I_{2{\rm BY}}}\langle x_i\rangle^{-(d+1+2k)}
\times\prod_{i\in [n]\backslash( I_{1{\rm BCO}} \cup I_{2{\rm BY}} )}\langle x_i\rangle^{-(d+1+k)}
\times\prod_{i\in I_{2{\rm BX}}}\langle y_i\rangle^{-(d+1+k)}
\]
\[
\times\prod_{i\in (I_1\cup I_2)\backslash I_{2{\rm BX}}}\langle y_i\rangle^{-(d+1)}
\times\prod_{i\in I_1\cup I_2}
\ \bbone\{|y_i-x_i|\le 2L^r\}\ 
\times
\prod_{i\in I_{\rm G}}
\bbone\left\{
\min_{j\neq i} |x_i-x_j|>\delta L^{r}
\right\}
\]
\[
\times
\prod_{i\in I_{\rm B}}
\bbone\left\{
\min_{j\neq i} |x_i-x_j|\le\delta L^{r}
\right\}
\times \prod_{i\in I_{1{\rm BCO}}} 
\frac{1}{|y_i-x_i|^{[A_i]+[B_i]-[C_i]+\epsilon}}
\]
\[
\times \prod_{i\in I_{1{\rm G}}} \frac{|y_i-x_i|^{\Delta_i+\gamma-[A_i]-[B_i]}}{{\rm NNDX}_{i}^{\Delta_i+\gamma+\epsilon}}
\times \prod_{i\in I_{2{\rm G}}} \frac{|y_i-x_i|^{\gamma}}{{\rm NNDX}_{i}^{\Delta_i+\gamma+\epsilon}}
\times \prod_{i\in I_{1{\rm BOO}}} \frac{1}{{\rm NNDY}_{i}^{[A_i]+\epsilon}\times{\rm NNDX}_{i}^{[B_i]+\epsilon}}
\]
\[
\times \prod_{i\in I_{1{\rm BCO}}} \frac{1}{{\rm NNDX}_{i}^{[C_i]+\epsilon}}
\times \prod_{i\in I_{2{\rm BY}}} \frac{1}{{\rm NNDY}_{i}^{\Delta_i+\epsilon}}
\times \prod_{i\in I_{2{\rm BX}}} \frac{1}{{\rm NNDX}_{i}^{\Delta_i+\epsilon}}
\times \prod_{i\in I_{34}} \frac{1}{{\rm NNDX}_{i}^{[D_i]+\epsilon}}\ .
\]
In order to reduce the complexity of this formula, we immediately integrate over the $y_i$'s
with $i\in I_{1{\rm G}}\cup I_{1{\rm BCO}}\cup I_{2{\rm G}}\cup I_{2{\rm BX}}$.
Indeed, the latter only couple to their corresponding $x_i$ through a simple dependence that can be dealt with using Lemma \ref{localL1bd}.

For $i\in I_{1{\rm G}}$,
we use the bound
\[
\int dy_i\ \bbone\{|y_i-x_i|\le 2L^r\}\ \langle y_i\rangle^{-(d+1)}\ |y_i-x_i|^{\Delta_i+\gamma-[A_i]-[B_i]}
\le O(1)\ L^{r(d+\Delta_i+\gamma-[A_i]-[B_i])}
\]
which follows from $\langle y_i\rangle\ge 1$ and Lemma \ref{localL1bd}.
For convenience, we did not write the domain of integration
since it is $\mathbb{R}^d$ minus a finite number of points (all the other $x_j$'s and $y_j$'s).

For $i\in I_{1{\rm BCO}}$,
we similarly use the bound
\[
\int dy_i\ \bbone\{|y_i-x_i|\le 2L^r\}\ \langle y_i\rangle^{-(d+1)}\ |y_i-x_i|^{-[A_i]-[B_i]+[C_i]-\epsilon}
\le O(1)\ L^{r(d-[A_i]-[B_i]+[C_i]-\epsilon)}\ .
\]

For $i\in I_{2{\rm G}}$,
we use the bound
\[
\int dy_i\ \bbone\{|y_i-x_i|\le 2L^r\}\ \langle y_i\rangle^{-(d+1)}\ |y_i-x_i|^{\gamma}
\le O(1)\ L^{r(d+\gamma)}\ .
\]

Finally, for $i\in I_{2{\rm BX}}$,
we use the bound
\[
\int dy_i\ \bbone\{|y_i-x_i|\le 2L^r\}\ \langle y_i\rangle^{-(d+1+k)}\ |y_i-x_i|^{0}
\le O(1)\ L^{rd}\ .
\]

As for the surviving long-distance decay factors, we bound them by the worst-case exponent, i.e., we turn them all into
$\langle\cdot\rangle^{-(d+1)}$. Therefore
\[
|\Upsilon_r|\le O(1)\sum_{\substack{(I_1,I_2,I_3)\vdash [m] \\ (I_1,I_2,I_3)\neq (\emptyset,\emptyset,[m]) }}
\sum_{\substack{ (I_{1{\rm G}},I_{1{\rm BOO}},I_{1{\rm BCO}})\vdash I_1 \\ (I_{2{\rm G}},I_{2{\rm BY}},I_{2{\rm BX}})\vdash I_2 }}
\ \ \sum_{(C_i)_{i\in I_{1{\rm BCO}}}\in \prod_{i\in I_{1{\rm BCO}}}\mathcal{A}(\Delta_i)}
\mathcal{P}^{({\rm V})} \mathcal{I}^{({\rm V})}
\]
\begin{flalign*}
&{\rm with}\ \ & &
\mathcal{P}^{({\rm V})}=
\prod_{i\in I_{1{\rm G}}} L^{r(\gamma-\epsilon)}\times
\prod_{i\in I_{1{\rm BOO}}} L^{r([A_i]+[B_i]-\Delta_i-d-\epsilon)}\times
\prod_{i\in I_{1{\rm BCO}}} L^{r([C_i]-\Delta_i-2\epsilon)} & & \ 
\end{flalign*}
\[
\times
\prod_{i\in I_{2{\rm G}}} L^{r(\gamma-2\epsilon)}\times
\prod_{i\in I_{2{\rm BY}}} L^{-r(d+2\epsilon)}\times
\prod_{i\in I_{2{\rm BX}}} L^{-2\epsilon r}
\]
and
\[
\mathcal{I}^{({\rm V})}=
\int_{{\rm Conf}_{n+|I_{1{\rm BOO}}|+|I_{2{\rm BY}}|}}\ 
\prod_{i\in [n]}dx_i\ 
\prod_{i\in I_{1{\rm BOO}}\cup I_{2{\rm BY}}}dy_i\  
\prod_{i\in [n]}\langle x_i\rangle^{-(d+1)}\times
\prod_{i\in I_{1{\rm BOO}}\cup I_{2{\rm BY}}}\langle y_i\rangle^{-(d+1)}
\]
\[
\times\prod_{I_{1{\rm BOO}}\cup I_{2{\rm BY}}}
\ \bbone\{|y_i-x_i|\le 2L^r\}\ 
\times
\prod_{i\in I_{\rm G}}
\bbone\left\{
\min_{j\neq i} |x_i-x_j|>\delta L^{r}
\right\}
\times
\prod_{i\in I_{\rm B}}
\bbone\left\{
\min_{j\neq i} |x_i-x_j|\le\delta L^{r}
\right\}
\]
\[
\times \prod_{i\in I_{1{\rm G}}} \frac{1}{{\rm NNDX}_{i}^{\Delta_i+\gamma+\epsilon}}
\times \prod_{i\in I_{2{\rm G}}} \frac{1}{{\rm NNDX}_{i}^{\Delta_i+\gamma+\epsilon}}
\times \prod_{i\in I_{1{\rm BOO}}} \frac{1}{{\rm NNDY}_{i}^{[A_i]+\epsilon}\times{\rm NNDX}_{i}^{[B_i]+\epsilon}}
\]
\[
\times \prod_{i\in I_{1{\rm BCO}}} \frac{1}{{\rm NNDX}_{i}^{[C_i]+\epsilon}}
\times \prod_{i\in I_{2{\rm BY}}} \frac{1}{{\rm NNDY}_{i}^{\Delta_i+\epsilon}}
\times \prod_{i\in I_{2{\rm BX}}} \frac{1}{{\rm NNDX}_{i}^{\Delta_i+\epsilon}}
\times \prod_{i\in I_{34}} \frac{1}{{\rm NNDX}_{i}^{[D_i]+\epsilon}}\ .
\]

Clearly, another round of consolidation of notation and relabeling is needed!
We will avoid the different treatment of the $x_i$'s versus the $y_i$'s by replacing the set of labels $[n]$ by a larger one we will
denote by $V$. The latter is a subset of $[n]\times\{{\rm X},{\rm Y}\}$
where ${\rm X}$, and ${\rm Y}$ are mere symbols. The set $V$ is nothing more nor less than what we need in order to label all
the variables in the previous integral. Thus, the set $V$ is the (disjoint) union of the following ten blocks.

\begin{itemize}
\item
Let $V_{1{\rm G}}^{\rm X}=I_{1{\rm G}}\times \{{\rm X}\}$.
\item
Let $V_{2{\rm G}}^{\rm X}=I_{2{\rm G}}\times \{{\rm X}\}$.
\item
Let $V_{1{\rm BOO}}^{\rm Y}=I_{1{\rm BOO}}\times \{{\rm Y}\}$
\item
Let $V_{1{\rm BOO}}^{\rm X}=I_{1{\rm BOO}}\times \{{\rm X}\}$
\item
Let $V_{1{\rm BCO}}^{\rm X}=I_{1{\rm BCO}}\times \{{\rm X}\}$.
\item
Let $V_{2{\rm BY}}^{\rm Y}=I_{2{\rm BY}}\times \{{\rm Y}\}$.
\item
Let $V_{2{\rm BY}}^{\rm X}=I_{2{\rm BY}}\times \{{\rm X}\}$.
\item
Let $V_{2{\rm BX}}^{\rm X}=I_{2{\rm BX}}\times \{{\rm X}\}$.
\item
Let $V_{34{\rm G}}^{\rm X}=\left(I_{34}\cap I_{\rm G}\right)\times\{{\rm X}\}$.
\item
Let $V_{34{\rm B}}^{\rm X}=\left(I_{34}\cap I_{\rm B}\right)\times\{{\rm X}\}$.
\end{itemize}

For each $a\in V$ we will define its corresponding (denominator) exponent $\beta_a$, we will also indicate its status, i.e.,
effective versus virtual from the point of view of the EFNNB we just used. This will result in a decomposition
$(V_{\rm eff},V_{\rm virt})\vdash V$. Elements $a\in V$ which were not present in the correlation estimated by the EFNNB, namely, those of
$V_{2{\rm BY}}^{X}$,
are declared virtual.
The point is that only effective elements can be somebody else's neighbor for the calculation of the ${\rm NNDX}_i$ and
${\rm NNDY}_i$ functions.
All of this is summarized in the following table which should help the reader follow the rest of the proof.

\bigskip
\begin{center}
\begin{tabular}{|c|c|c|} \hline
Block & Status & $\beta_{(i,s)}$ \\ \hline
$V_{1{\rm G}}^{\rm X}$ & eff & $\Delta_i+\gamma+\epsilon$ \\ \hline
$V_{2{\rm G}}^{\rm X}$ & eff & $\Delta_i+\gamma+\epsilon$ \\ \hline
$V_{1{\rm BOO}}^{\rm Y}$ & eff & $[A_i]+\epsilon$ \\ \hline
$V_{1{\rm BOO}}^{\rm X}$ & eff & $[B_i]+\epsilon$ \\ \hline
$V_{1{\rm BCO}}^{\rm X}$ & eff & $[C_i]+\epsilon$ \\ \hline
$V_{2{\rm BY}}^{\rm Y}$ & eff & $\Delta_i+\epsilon$ \\ \hline
$V_{2{\rm BY}}^{\rm X}$ & virt & $0$ \\ \hline
$V_{2{\rm BX}}^{\rm X}$ & eff & $\Delta_i+\epsilon$ \\ \hline
$V_{34{\rm G}}^{\rm X}$ & eff & $[D_i]+\epsilon$ \\ \hline
$V_{34{\rm B}}^{\rm X}$ & eff & $[D_i]+\epsilon$ \\ \hline
\end{tabular}
\end{center}

We let $V^{\rm X}$ denote the set of all $(i,s)\in V$ with $s={\rm X}$. We likewise define
$V^{\rm Y}$ denote the set of all $(i,s)\in V$ with $s={\rm Y}$.
We let $V_{\rm G}$ denote the set of all $(i,s)\in V$ with $i\in I_{\rm G}$.
We likewise define
$V_{\rm B}$ as the set of all $(i,s)\in V$ with $i\in I_{\rm B}$.

We also introduce notation for the following sets.
\begin{itemize}
\item
We let $V_{\rm G}^{\rm X}=V_{\rm G}\cap V^{\rm X}=V_{1{\rm G}}^{\rm X}\cup V_{2{\rm G}}^{\rm X}\cup V_{34{\rm G}}^{\rm X}$.
\item
We let $V_{\rm B}^{\rm X}=V_{\rm B}\cap V^{\rm X}=V_{1{\rm BOO}}^{\rm X}\cup V_{2{\rm BY}}^{\rm X}
\cup V_{2{\rm BX}}^{\rm X}\cup V_{34{\rm B}}^{\rm X}$.
\item
We let $V_{\rm B}^{\rm Y}=V_{\rm B}\cap V^{\rm Y}=V_{1{\rm BOO}}^{\rm Y}\cup V_{2{\rm BY}}^{\rm Y}$.
\end{itemize}
Of course, there is no need for a fourth set $V_{\rm G}^{\rm Y}$ since it would be empty.

We define an involution $\iota:V\rightarrow V$ as follows.
\begin{itemize}
\item
For all $i\in I_{1{\rm BOO}}\cup I_{2{\rm BY}}$, we let $\iota(i,{\rm Y})=(i,{\rm X})$ and $\iota(i,{\rm X})=(i,{\rm Y})$.
\item
For all $i\in [n]\backslash\left( I_{1{\rm BOO}}\cup I_{2{\rm BY}} \right) $, we let $\iota(i,{\rm X})=(i,{\rm X})$
\end{itemize}

We can now rewrite the last integral as
\[
\mathcal{I}^{({\rm V})}=
\int_{{\rm Conf}_{|V|}}\ 
\prod_{a\in V}du_a
\ \prod_{a\in V}\langle u_a\rangle^{-(d+1)}
\times\prod_{a\in V_{1{\rm BOO}}^{\rm Y}\cup V_{2{\rm BY}}^{\rm Y}}
\ \bbone\{|u_a-u_{\iota(a)}|\le 2L^r\}
\]
\[
\times
\prod_{a\in V_{\rm G}^{\rm X}}
\bbone\left\{
\min_{b\in V^{\rm X}\backslash\{a\}} |u_a-u_b|>\delta L^{r}
\right\}
\times
\prod_{a\in V_{\rm B}^{\rm X}}
\bbone\left\{
\min_{b\in V^{\rm X}\backslash\{a\}} |u_a-u_b|\le \delta L^{r}
\right\}
\]
\[
\times \prod_{a\in V_{\rm eff}}
\left(\min_{b\in V_{\rm eff}\backslash\{ a\}} |u_a-u_b|\right)^{-\beta_a}\ .
\]
The next step which we call preemptive rerouting is necessary in order to keep the complexity of the graphs arising
in the next section under control.
We will overestimate the last integral by replacing $V_{\rm eff}$ simply by $V$.
Indeed, $V\backslash V_{\rm eff}=V_{2{\rm BY}}^{\rm X}$.
There is no harm in adding the factors
\[
\left(\min_{b\in V\backslash\{ a\}} |u_a-u_b|\right)^{-\beta_a}
\]
for $a\in V_{2{\rm BY}}^{\rm X}$ since the corresponding exponent is $\beta_a=0$.
As for the elements $a\in V_{\rm eff}$, we have
\[
\left(\min_{b\in V_{\rm eff}\backslash\{ a\}} |u_a-u_b|\right)^{-\beta_a}
\le
\left(\min_{b\in V\backslash\{ a\}} |u_a-u_b|\right)^{-\beta_a}
\]
since all exponents $\beta_a$ in the table are nonnegative.

As a result,
\[
|\Upsilon_r|\le O(1)\sum_{\substack{(I_1,I_2,I_3)\vdash [m] \\ (I_1,I_2,I_3)\neq (\emptyset,\emptyset,[m]) }}
\sum_{\substack{ (I_{1{\rm G}},I_{1{\rm BOO}},I_{1{\rm BCO}})\vdash I_1 \\ (I_{2{\rm G}},I_{2{\rm BY}},I_{2{\rm BX}})\vdash I_2 }}
\ \ \sum_{(C_i)_{i\in I_{1{\rm BCO}}}\in \prod_{i\in I_{1{\rm BCO}}}\mathcal{A}(\Delta_i)}
\mathcal{P}^{({\rm VI})} \mathcal{I}^{({\rm VI})}
\]
with
\[
\mathcal{P}^{({\rm VI})}=\mathcal{P}^{({\rm V})}=
\prod_{i\in I_{1{\rm G}}} L^{r(\gamma-\epsilon)}\times
\prod_{i\in I_{1{\rm BOO}}} L^{r([A_i]+[B_i]-\Delta_i-d-\epsilon)}\times
\prod_{i\in I_{1{\rm BCO}}} L^{r([C_i]-\Delta_i-2\epsilon)}
\]
\[
\times
\prod_{i\in I_{2{\rm G}}} L^{r(\gamma-2\epsilon)}\times
\prod_{i\in I_{2{\rm BY}}} L^{-r(d+2\epsilon)}\times
\prod_{i\in I_{2{\rm BX}}} L^{-2\epsilon r}
\]
and
\[
\mathcal{I}^{({\rm VI})}=
\int_{{\rm Conf}_{|V|}}\ 
\prod_{a\in V}du_a 
\ \prod_{a\in V}\langle u_a\rangle^{-(d+1)}
\times\prod_{a\in V_{1{\rm BOO}}^{\rm Y}\cup V_{2{\rm BY}}^{\rm Y}}
\ \bbone\{|u_a-u_{\iota(a)}|\le 2L^r\}
\]
\[
\times
\prod_{a\in V_{\rm G}^{\rm X}}
\bbone\left\{
\min_{b\in V^{\rm X}\backslash\{a\}} |u_a-u_b|>\delta L^{r}
\right\}
\times
\prod_{a\in V_{\rm B}^{\rm X}}
\bbone\left\{
\min_{b\in V^{\rm X}\backslash\{a\}} |u_a-u_b|\le \delta L^{r}
\right\}
\]
\begin{equation}
\times \prod_{a\in V} \left(\min_{b\in V\backslash\{ a\}} |u_a-u_b|\right)^{-\beta_a}\ .
\label{Isixeq}
\end{equation}

\subsection{Graph construction}

We are now ready for a more involved version of the argument in \S\ref{pinsum}.
We let $\mathcal{G}$ be the set of all pairs $(\tau,\sigma)$ made of fixed-point-free endofunctions
$\tau:V\rightarrow V$ and $\sigma:V_{\rm B}\rightarrow V_{\rm B}$ which satisfy the crucial requirement
$\tau(V_{\rm B})\subset V_{B}$.

{\bf Claim:} For any point configuration $(u_a)_{a\in V}\in {\rm Conf}_{|V|}$ which satisfies all the conditions in the indicator functions appearing in (\ref{Isixeq}),
the following inequality holds
\[
1\le \sum_{(\tau,\sigma)\in\mathcal{G}}
\prod_{a\in V}\bbone\left\{|u_a-u_{\tau(a)}|=\min_{b\in V\backslash\{a\}}|u_a-u_b|\right\}
\]
\[
\times
\prod_{a\in V_{\rm B}}\bbone\left\{|u_a-u_{\sigma(a)}|\le \delta L^{r} \right\}
\times
\prod_{a\in V_{\rm B}}\bbone\left\{|u_a-u_{\tau(a)}|\le \delta L^{r} \right\}\ .
\]

To prove the claim, we need to construct a pair $(\tau,\sigma)$ with the desired properties.
For $\tau$ we proceed exactly as in \S\ref{pinsum},
namely, we choose a nearest neighbor. This already takes care of the
first set of indicator functions.
As for the construction of $\sigma$, it is done as follows. For $a\in V_{\rm B}^{\rm Y}$, we let $\sigma(a)=\iota(a)\in V_{\rm B}^{\rm X}$.
For $a\in V_{\rm B}^{\rm X}$ we let $\sigma(a)$ be some choice of element $b\in V^{X}\backslash\{a\}$ such that
$|u_a-u_b|\le \delta L^r$. Since we can arrange for $\delta\ge 2$,
this definition clearly satisfies all conditions in the second group of indicator functions.
Note that such a $b$ cannot belong to $V_{\rm G}^{\rm X}$ by the repulsive condition in the definition of the latter.
So $\sigma$ is indeed an endofunction of $V_{\rm B}$.

We only need to show that, for any $a\in V_{\rm B}$, we have $\tau(a)\in V_{\rm B}$ and $|u_a-u_{\tau(a)}|\le \delta L^r$.
For $a\in V_{\rm B}^{\rm Y}$, $\sigma(a)\in V_{\rm B}^{\rm X}\subset V_{\rm B}$ by construction, while the inequalities
\[
2L^{r}\ge |u_a-u_{\iota(a)}|\ge \min_{b\in V\backslash\{a\}}|u_a-u_b|=|u_a-u_{\tau(a)}|
\]
and $\delta\ge 2$ imply $|u_a-u_{\tau(a)}|\le \delta L^r$.
Now let $a\in V_{\rm B}^{\rm X}$. Then $\sigma(a)\in V_{\rm B}^{\rm X}\backslash\{a\}$ satisfies
\[
\delta L^r\ge |u_a-u_{\sigma(a)}|\ge \min_{b\in V\backslash\{a\}}|u_a-u_b|=|u_a-u_{\tau(a)}|\ .
\]
Finally, we argue by contradiction an assume that $\tau(a)\notin V_{\rm B}$.
Then $\tau(a)\in V_{\rm G}^{\rm X}$ which by definition implies
$|u_{\tau(a)}-u_{\iota(a)}|>\delta L^r$
because $\iota(a)\in V_{\rm B}^{\rm X}$ must be different from $\tau(a)$.
On the other hand $\iota(a)\neq a$ satisfies
\[
2L^r\ge |u_a-u_{\iota(a)}|\ge \min_{b\in V\backslash\{a\}}|u_a-u_b|=|u_a-u_{\tau(a)}|
\]
\begin{flalign*}
&{\rm while}\ \ & &
|u_a-u_{\tau(a)}|\ge |u_{\tau(a)}-u_{\iota(a)}|-|u_{a}-u_{\iota(a)}|>(\delta-2)L^r & & \ 
\end{flalign*}
which is a contradiction because we can arrange for $\delta\ge 4$.

Now that the claim is proved, we do as in \S\ref{pinsum}. We pull the sum over $(\tau,\sigma)$ out of the integral, use the first set
of indicator functions to replace the $\min$'s by $|u_a-u_{\tau(a)}|$'s. Finally, we discard the old indicator functions and
the ones from the first group we just used, and we keep the new ones from the second and third group.
As a result, we have 
\[
|\Upsilon_r|\le O(1)\sum_{\substack{(I_1,I_2,I_3)\vdash [m] \\ (I_1,I_2,I_3)\neq (\emptyset,\emptyset,[m]) }}
\sum_{\substack{ (I_{1{\rm G}},I_{1{\rm BOO}},I_{1{\rm BCO}})\vdash I_1 \\ (I_{2{\rm G}},I_{2{\rm BY}},I_{2{\rm BX}})\vdash I_2 }}
\ \ \sum_{(C_i)_{i\in I_{1{\rm BCO}}}\in \prod_{i\in I_{1{\rm BCO}}}\mathcal{A}(\Delta_i)}
\sum_{(\tau,\sigma)\in\mathcal{G}}
\mathcal{P}^{({\rm VII})} \mathcal{I}^{({\rm VII})}
\]
\begin{flalign*}
&{\rm with}\ \ & &
\mathcal{P}^{({\rm VII})}=\mathcal{P}^{({\rm VI})}=
\prod_{i\in I_{1{\rm G}}} L^{r(\gamma-\epsilon)}\times
\prod_{i\in I_{1{\rm BOO}}} L^{r([A_i]+[B_i]-\Delta_i-d-\epsilon)}\times
\prod_{i\in I_{1{\rm BCO}}} L^{r([C_i]-\Delta_i-2\epsilon)} & & \ 
\end{flalign*}
\[
\times
\prod_{i\in I_{2{\rm G}}} L^{r(\gamma-2\epsilon)}\times
\prod_{i\in I_{2{\rm BY}}} L^{-r(d+2\epsilon)}\times
\prod_{i\in I_{2{\rm BX}}} L^{-2\epsilon r}
\]
\begin{flalign*}
&{\rm and}\ \ & &
\mathcal{I}^{({\rm VII})}=
\int_{{\rm Conf}_{|V|}}\ 
\prod_{a\in V}du_a
\ \prod_{a\in V}\langle u_a\rangle^{-(d+1)}
\times
\prod_{a\in V_{\rm B}}\bbone\left\{|u_a-u_{\sigma(a)}|\le \delta L^{r} \right\} & & \ 
\end{flalign*}
\begin{equation}
\times
\prod_{a\in V_{\rm B}}\bbone\left\{|u_a-u_{\tau(a)}|\le \delta L^{r} \right\}
\times \prod_{a\in V} |u_a-u_{\tau(a)}|^{-\beta_a}\ .
\label{integralVII}
\end{equation}

\subsection{Two-scale pin and sum argument}\label{doublepinsumsec}

The workhorses in \S\ref{elementarybds} having now been provided with proper steering, we can sit back and watch them
take care of the remaining integral. Indeed, one has a directed graph on $V$ with two types of edges: $a\rightarrow \tau(a)$ and
$a\rightarrow\sigma(a)$. Moreover, we have a decomposition $(V_{\rm G},V_{\rm B})\vdash V$.
Clearly, the restriction of the graph to $V_{\rm G}$ is such that there can only be two kinds of connected components:
1) isolated pure $\tau$ hairy cycles or, 2) directed $\tau$ trees attached to some vertex in $V_{\rm B}$.
There is no $\sigma$ edge incident to a vertex in $V_{\rm G}$. Furthermore, the only possible communication between $V_{\rm G}$ and
$V_{\rm B}$ is through an edge going from a vertex in $V_{\rm G}$ to a vertex in $V_{\rm B}$,
and this can only happen at most once per connected component
of $V_{\rm G}$.
We therefore repeat the pin and sum argument from \S\ref{pinsum}
and integrate over the vertices of $V_{\rm G}$ which produces an $O(1)$ factor.

Let $(W_1,\ldots,W_q)\vdash V_{\rm B}$ be a decomposition into connected components for the remaining graph made of both $\tau$
and $\sigma$ edges. For $i$ with $1\le i\le q$, further decompose $W_{i}$ into connected components for the subgraph only made of $\tau$ edges. 
This gives a decomposition $(W_{i,1},\ldots,W_{i,l_{q}})\vdash W_i$. Then delete all $\sigma$ edges internal to the $W_{i,j}$.
Also, keep enough $\sigma$ edges, i.e., $l_{q}-1$ of them
which connect $W_{i,1},\ldots,W_{i,l_q}$ together and delete the rest. When a $\sigma$ edge is deleted
we discard the corresponding indicator function in (\ref{integralVII}).
Also pick some $b_i\in W_{i,1}$ which will serve as root for all of $W_i$.
Discard the $\langle u_a\rangle^{-(d+1)}$ factors except those of $b_1,\ldots,b_q$.
As a result, each $\tau$ edge $a\rightarrow \tau(a)$
contributes an $L^{-\beta_a r}$ and each vertex contributes $L^{dr}$ except the $b_i$'s which have to be integrated with Lemma
\ref{globalL1bd} instead of Lemma \ref{localL1bd}.

Again, an explicit example with a picture should help the reader follow the argument.
\[
\parbox{12cm}{
\includegraphics[width=12cm]{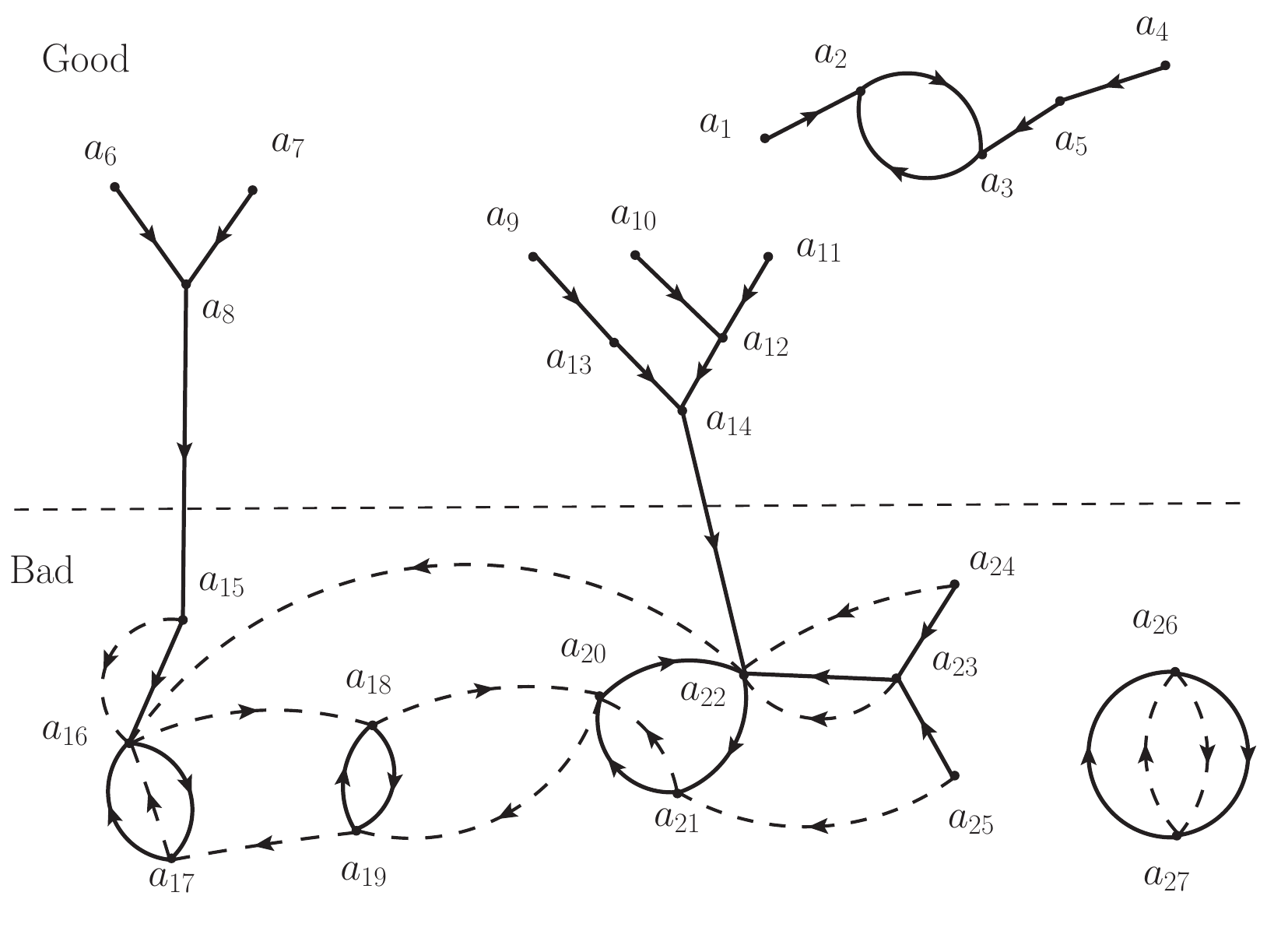}}
\]

In the above graph, the solid lines represent $\tau$ edges and the dashed ones represent the $\sigma$ edges.
The elements in $V_{\rm G}$ are at the top and the ones in $V_{\rm B}$ are at the bottom.
We first perform the succession of integrations corresponding to the sequence
\[
a_1,a_4,a_5,a_3,a_2,a_6,a_7,a_8,a_9,a_{10},a_{11},a_{12},a_{13},a_{14}
\]
in order to get rid of the good points.
These operations all use Lemma \ref{globalL1bd}, except the treatment of $a_3$ which uses Lemma \ref{globalbetabd}.
We are then left with the following picture.
\[
\parbox{12cm}{
\includegraphics[width=12cm]{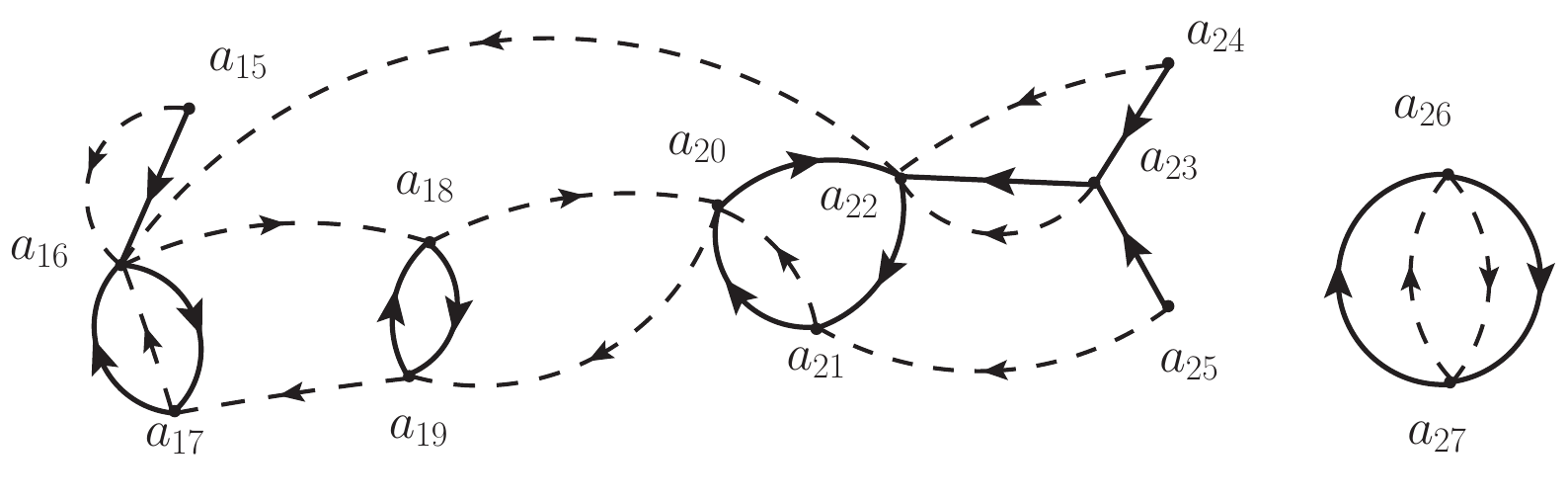}}
\]
In this example we have
$
W_1=\{a_{15},a_{16},a_{17},a_{18},a_{19},a_{20},a_{21},a_{22},a_{23},a_{24},a_{25}\}$
and
$W_2=\{a_{26},a_{27}\}$, while
the list of sub-components is given by
$W_{1,1}=\{a_{15},a_{16},a_{17}\}$,
$W_{1,2}=\{ a_{18},a_{19} \}$,
$W_{1,3}=\{ a_{20},a_{21},a_{22},a_{23},a_{24},a_{25}\}$
and
$W_{2,1}=\{a_{26},a_{27}\}$.
Then we remove $\sigma$ edges which do not improve connectivity, until we have a spanning tree
between $W_{i,j}$'s in each $W_i$.
Of course there are many ways to do so. For example, one could obtain the following graph.
\[
\parbox{12cm}{
\includegraphics[width=12cm]{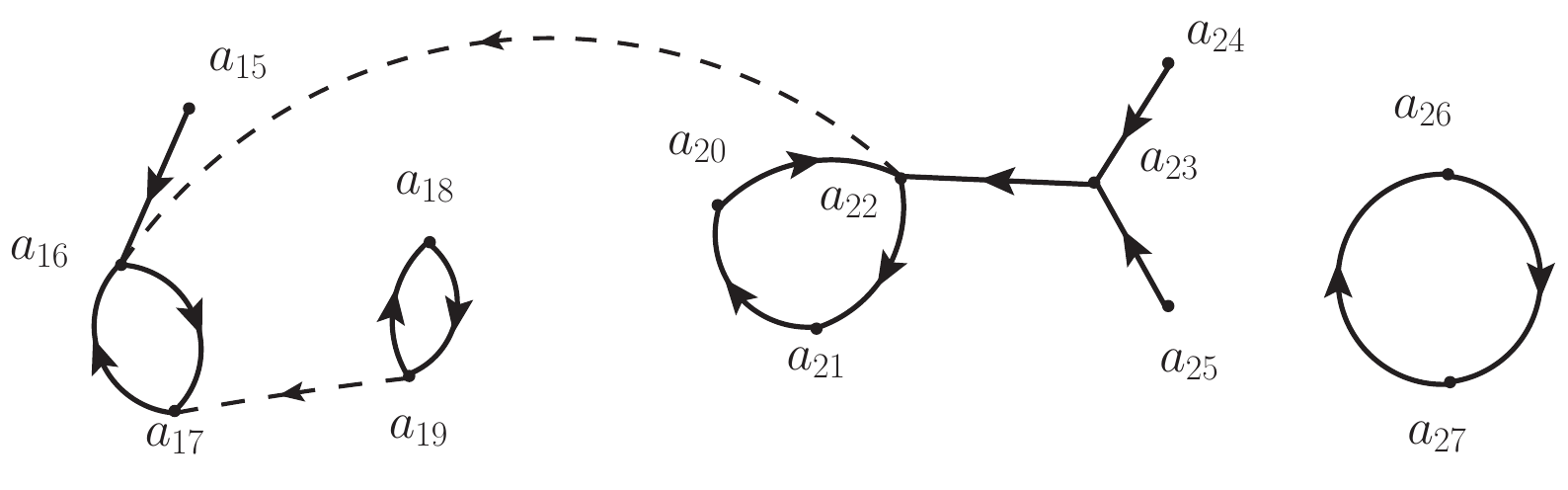}}
\]
For the choices of roots, we could for instance take $b_1=a_{15}$ and $b_2=a_{26}$.
A possible order of integration is then 
$a_{27},a_{26},a_{24},a_{25},a_{23},a_{20},a_{21},a_{22},a_{18},a_{19},a_{17},a_{16},a_{15}$.
Note that we used the idea from Remark \ref{contrarian} in the last three integrations.
For $a_{15}$ and $a_{26}$, we used Lemma \ref{globalL1bd}. For $a_{27}$, $a_{20}$, $a_{18}$ and $a_{17}$, we used
Lemma \ref{localbetabd}. Finally, for the remaining points, we used Lemma \ref{localL1bd}.

\subsection{Power-counting}

At this point we are left with an accounting problem which is to make sure the overall coefficient of $r$ in the exponent of $L$
is positive.
Indeed,
\[
|\Upsilon_r|\le O(1)\sum_{\substack{(I_1,I_2,I_3)\vdash [m] \\ (I_1,I_2,I_3)\neq (\emptyset,\emptyset,[m]) }}
\sum_{\substack{ (I_{1{\rm G}},I_{1{\rm BOO}},I_{1{\rm BCO}})\vdash I_1 \\ (I_{2{\rm G}},I_{2{\rm BY}},I_{2{\rm BX}})\vdash I_2 }}
\ \ \sum_{(C_i)_{i\in I_{1{\rm BCO}}}\in \prod_{i\in I_{1{\rm BCO}}}\mathcal{A}(\Delta_i)}
\sum_{(\tau,\sigma)\in\mathcal{G}}
\mathcal{P}^{({\rm VIII})} \mathcal{I}^{({\rm VIII})}
\]
with $\mathcal{I}^{({\rm VIII})}=1$
but also
\[
\mathcal{P}^{({\rm VIII})}=
\prod_{i\in I_{1{\rm G}}} L^{r(\gamma-\epsilon)}\times
\prod_{i\in I_{1{\rm BOO}}} L^{r([A_i]+[B_i]-\Delta_i-d-\epsilon)}\times
\prod_{i\in I_{1{\rm BCO}}} L^{r([C_i]-\Delta_i-2\epsilon)}
\]
\[
\times
\prod_{i\in I_{2{\rm G}}} L^{r(\gamma-2\epsilon)}\times
\prod_{i\in I_{2{\rm BY}}} L^{-r(d+2\epsilon)}\times
\prod_{i\in I_{2{\rm BX}}} L^{-2\epsilon r}
\times L^{-rdq}\times \prod_{a\in V_{\rm B}} L^{(d-\beta_a) r}\ .
\]
Using the table, this can be reorganized according to subsets of $I$ as follows.
We have
$\mathcal{P}^{({\rm VIII})}=L^{r \alpha_{\rm Total}}$
with
\[
\alpha_{\rm Total}=-dq+
(\gamma-\epsilon)|I_{1{\rm G}}|+(\gamma-2\epsilon)|I_{2{\rm G}}|
\]
\[
+\sum_{i\in I_{1{\rm BOO}}\cup I_{1{\rm BCO}}\cup I_{2{\rm BY}}\cup I_{2{\rm BX}}}
(d-\Delta_i-3\epsilon)
+\sum_{i\in I_{34}\cap I_{\rm B}}
(d-[D_i]-\epsilon)\ .
\]
Now notice that each connected component $W_i$ must contain at least two elements of $V_{\rm B}^{\rm X}$.
Indeed, it is not empty by definition. If it contains $a\in V_{\rm B}^{\rm Y}$ then by construction it also
contains $\sigma(a)\in V_{\rm B}^{\rm X}$.
But if we caught one element $b$ of  $V_{\rm B}^{\rm X}$, then we can also catch another one, namely, $\sigma(b)$.
Thus
\[
2q\le |V_{\rm B}^{\rm X}|=|I_{1{\rm BOO}}|+ |I_{1{\rm BCO}}|+ |I_{2{\rm BY}}|+ |I_{2{\rm BX}}|+
|I_{34}\cap I_{\rm B}|=|I_{\rm B}|
\]
and therefore
\[
\alpha_{\rm Total}\ge -\frac{d}{2}|I_{\rm B}|
+(\gamma-2\epsilon)\left(|I_{1{\rm G}}|+|I_{2{\rm G}}|\right)
\]
\[
+|I_{\rm B}|\times
\min\left\{
\min_{i\in I_1\cup I_2}(d-\Delta_i-3\epsilon),
\min_{i\in I_3}(d-\Delta_i-\epsilon),
\min_{i\in I_4}(d-[A_i]-\epsilon)
\right\}
\]
\[
\ge (\gamma-2\epsilon)\left(|I_{1{\rm G}}|+|I_{2{\rm G}}|\right)
\]
\[
+\left(|I_{1{\rm B}}|+|I_{2{\rm B}}|+|I_{34}\cap I_{\rm B}|\right)
\times
\min\left\{
\min_{i\in I_1\cup I_2\cup I_3}\left(\frac{d}{2}-\Delta_i-3\epsilon\right),
\min_{i\in I_4}\left(\frac{d}{2}-[A_i]-\epsilon\right)
\right\}\ .
\]
We toss $|I_{34}\cap I_{\rm B}|$ away and use $|I_1|=|I_{1{\rm G}}|+|I_{1{\rm B}}|$ as well as $|I_2|=|I_{2{\rm G}}|+|I_{2{\rm B}}|$ in
order to write
$\alpha_{\rm Total}\ge
\nu\left(|I_1|+|I_2|\right)$
where
\[
\nu=\min\left\{\gamma-2\epsilon,
\min_{1\le i\le m}\left(\frac{d}{2}-\Delta_i-3\epsilon\right),
\min_{m+1\le i\le n}\left(\frac{d}{2}-[A_i]-\epsilon\right)
\right\}>0\ .
\]
Recalling that $(I_1,I_2,I_3)\neq(\emptyset,\emptyset,[m])$ and $r\le 0$, and all the sums being finite, we finally obtain
$|\Upsilon_r|\le O(1) L^{\nu r}$,
namely, the statement of Proposition \ref{combiprop}.
\qed

\section{Proof of the main theorem}\label{finishproof}

We will be brief since now the rest of the proof of Theorem \ref{maintheorem} is standard.
For Part (1), it is enough to take $p\ge 2$ to be an even integrer.
We write
\[
||M_{r}(f)-M_{r-1}(f)||_{L^p}^p=
\mathbb{E}\left[M_{r}(f)-M_{r-1}(f)
\right]^p
=\sum_{q=0}^p \left(\begin{array}{c} p \\ q \end{array}\right)
(-1)^{p-q}\  \mathbb{E}\left[ M_{r}(f)^{p-q} M_{r-1}(f)^{q}\right] 
\]
\[
=\sum_{q=0}^p \left(\begin{array}{c} p \\ q \end{array}\right)
(-1)^{p-q} \left(\mathbb{E}\left[ M_{r}(f)^{p-q} M_{r-1}(f)^{q}\right]-{\rm IPC}\right) 
\]
\begin{flalign*}
&{\rm with}\ \ & &
{\rm IPC}=\int_{{\rm Conf}_p}\prod_{i=1}^{p} dx_i\ \prod_{i=1}^{p} f(x_i)
\left\langle
\prod_{i=1}^{p}\mathcal{O}_{C_{\ast i}}(x_i)
\right\rangle\ . & & \ 
\end{flalign*}
We apply Proposition \ref{combiprop} to each term and deduce convergence of the telescopic series in $L^p$, and
thus in $L^1$ from which the almost sure convergence follows too.
We then prove part (2).
If one used two different mollifiers $\rho_1$ and $\rho_2$
and constructed two versions $\mathcal{O}_{C_{\ast},1}(f)$ and $\mathcal{O}_{C_{\ast},2}(f)$ of the
smeared renormalized product corresponding to the label $C_{\ast}$, then
\[
||\mathcal{O}_{C_{\ast},1}(f)-\mathcal{O}_{C_{\ast},1}(f)||_{L^2}^{2}=
\mathbb{E}\ \mathcal{O}_{C_{\ast},1}(f)^2
-2\ \mathbb{E}\ \mathcal{O}_{C_{\ast},1}(f)\mathcal{O}_{C_{\ast},2}(f)
+\mathbb{E}\ \mathcal{O}_{C_{\ast},2}(f)^2=0
\]
by a similar $r\rightarrow -\infty$ limit and application of Proposition \ref{combiprop}.

It is also easy to similarly show that by adjoining the random variables $\mathcal{O}_{C_{\ast}}(f)$, $f\in \mathcal{S}(\mathbb{R}^d)$,
the pointwise representation
from \S\ref{probareal} still holds.
Hence $f\rightarrow \mathcal{O}_{C_{\ast}}(f)$ is continuous
in say $L^2$ and thus in probability.
It is just a matter of collating this generalized random field into a random distribution (see~\cite[Proposition III.4.2 (a)]{Fernique}),
in order to finish the proof of the theorem.
\qed 

\section{A detailed example: the fractional massless free field}\label{examplesec}

\subsection{The explicit abstract system of pointwise correlations and the BNNFB}

We fix throughout this section the parameter $[\phi]>0$.
For any integer $r\in\mathbb{N}$, let $\mathcal{A}_r=(\mathbb{N}^d)^r/\mathfrak{S}_r$, namely the set of
$r$-tuples $(\nu(1),\ldots,\nu(r))$ of multiindices $\nu(i)\in\mathbb{N}^d$ modulo the action of the symmetric group $\mathfrak{S}_r$ by permutation of these multiindices. We will denote the equivalence class or $\mathfrak{S}_r$-orbit of the tuple $(\nu(1),\ldots,\nu(r))$ by $[\nu(1),\ldots,\nu(r)]$.
More generally, for any finite set $F$ of cardinality $r$ and any function $\nu:F\rightarrow\mathbb{N}^d$ we define a corresponding element of $\mathcal{A}_r$ given by $\{F;\nu\}=[\nu(a_1),\ldots,\nu(a_r)]$
where $a_1,\ldots,a_r$ is some enumeration of the elements of $F$.
We define the disjoint union $\mathcal{A}^{\infty}=\cup_{r\ge 0}\mathcal{A}_r$.
The unique element of $\mathcal{A}_{0}$ is denoted by $\bbone$.
Using the language of multilinear algebra, $\mathcal{A}^{\infty}$ can also be seen as a set of labels for a ``monomials of monomials'' linear basis for the plethysm ${\rm Sym}({\rm Sym}(\mathbb{R}^{d}))$
where ${\rm Sym}(V)=\oplus_{r\ge 0}{\rm Sym}^r (V)$ is the symmetric algebra obtained from a real vector space
$V$. 
For any $A\in\mathcal{A}^{\infty}$ we define a notion of scaling dimension 
\[
[A]=r[\phi]+|\nu(1)|+\ldots+|\nu(r)|
\]
and degree ${\rm deg}(A)=r$, if $A=[\nu(1),\ldots,\nu(r)]$. We recall that for a multiindex $\alpha=(\alpha_1,\ldots,\alpha_d)\in\mathbb{N}^d$ we will use the standard notations for the length $|\alpha|=\alpha_1+\cdots+\alpha_d$, the factorial $\alpha!=\alpha_1!\cdots\alpha_d!$, the differential operator $\partial_x^{\alpha}=
\frac{\partial^{|\alpha|}}{\partial x_1^{\alpha_1}\cdots\partial x_d^{\alpha_d}}$ and the monomial $x^\alpha=x_1^{\alpha_1}\cdots x_d^{\alpha_d}$.
Such an element $A$ will serve as a label for the composite field formally given by
\[
\mathcal{O}_{A}(x)=:\partial^{\nu(1)}\phi(x)\cdots\partial^{\nu(r)}\phi(x):\ .
\]
Because of the hypothesis $[\phi]>0$, for any $\Delta\in\mathbb{R}$, the set
$\{A\in\mathcal{A}^{\infty}\ |\  [A]\le\Delta\}$ is finite. We also have that the quantity
${\rm next}(\Delta)=\min\{[A]\ |\ A\in\mathcal{A}^{\infty}, [A]>\Delta\}$ is well defined and satisfies
${\rm next}(\Delta)>\Delta$.
The present example thus has gapped dimension spectrum as in the statement of Conjecture \ref{physconj}.
The application of our main theorem will concern the finite alphabet $\mathcal{A}=\{A\in\mathcal{A}^{\infty}\ |\ [A]<\frac{d}{2}\}$, but we will first define an abstract system of pointwise correlations
$\langle\mathcal{O}_{A_1}(x_1)\cdots\mathcal{O}_{A_n}(x_n)\rangle$, more generally for any labels $A_1,\ldots,A_n$ in $\mathcal{A}^{\infty}$. 
For such a collection of labels, pick a finite set $F$ of cardinality $|F|={\rm deg}(A_1)+\cdots+{\rm deg}(A_n)$.
Pick a decomposition (ordered collection of disjoint subsets whose union is $F$) $F_1,\ldots,F_n$ such that
$|F_i|={\rm deg}(A_i)$. Also pick a function $\nu:F\rightarrow\mathbb{N}^d$ such that, for all $i$, $\{F_i;\nu|_{F_i}\}=A_i$.
The previous data determine a function $\iota:F\rightarrow[n]=\{1,\ldots,n\}$ which maps $a\in F$ to the unique $i$ such that $a\in F_i$. 

Recall from \S\ref{introtointro}
that for $x\neq y$ in $\mathbb{R}^d$, $C(x,y)=\frac{\varkappa}{|x-y|^{2[\phi]}}$. We introduce the more general notation
$C^{\alpha,\beta}(x,y)=\partial_x^{\alpha}\partial_y^{\beta}C(x,y)$
for any $\alpha,\beta\in\mathbb{N}^d$. As a consequence of the symmetry $C(y,x)=C(x,y)$, we have
\begin{equation}
C^{\beta,\alpha}(y,x)=C^{\alpha,\beta}(x,y)\ .
\label{Csymmetry}
\end{equation}
By {\em definition}, we let for $(x_1,\ldots,x_n)\in{\rm Conf}_n$,
\begin{equation}
\langle\mathcal{O}_{A_1}(x_1)\cdots\mathcal{O}_{A_n}(x_n)\rangle=
\sum_{\mathcal{W}\ {\rm off-diag}}
\prod_{\{a,b\}\in\mathcal{W}} C^{\nu(a),\nu(b)}(x_{\iota(a)},x_{\iota(b)})\ .
\label{explicitdef}
\end{equation}
Here $\mathcal{W}$ is a set partition of $F$ with pairs only, i.e., a perfect matching.
The sum is over all such $\mathcal{W}$'s which are off-diagonal. This means that no (unordered) pair $\{a,b\}\in\mathcal{W}$ is allowed to be contained in one of the blocks $F_i$ of the given decomposition of $F$.
Thus $\{a,b\}\in\mathcal{W}$ implies $\iota(a)\neq\iota(b)$ and therefore $x_{\iota(a)}\neq x_{\iota(b)}$.
Also note that by the symmetry (\ref{Csymmetry}), a factor $C^{\nu(a),\nu(b)}(x_{\iota(a)},x_{\iota(b)})$ does not depend on the choice of ordered pair $(a,b)$ or $(b,a)$ corresponding to a given unordered pair $\{a,b\}$.

If $A_1=[\nu(1,1),\ldots,\nu(1,r_1)],\ldots,A_n=[\nu(n,1),\ldots,\nu(n,r_n)]$, then
in better agreement with physics literature notation, the previous pointwise correlation could be written
more explicitly as
\[
\langle\mathcal{O}_{A_1}(x_1)\cdots\mathcal{O}_{A_n}(x_n)\rangle=
\]
\begin{equation}
\langle
:\partial^{\nu(1,1)}\phi(x_1)\cdots\partial^{\nu(1,r_1)}\phi(x_1):\times\cdots\times
:\partial^{\nu(n,1)}\phi(x_n)\cdots\partial^{\nu(n,r_n)}\phi(x_n):
\rangle
\label{defcorrel}
\end{equation}
which formally represents a correlation of several Wick-ordered local monomials in the elementary scalar field $\phi$ and its derivatives.
In this case, we take $F=\{(i,j)\ |\ 1\le i\le n, 1\le j\le r_i\}\subset\mathbb{N}^2$ and
$F_i=\{(i,j)\ | 1\le j\le r_i\}$ for $1\le i\le n$. The function $\nu:F\rightarrow \mathbb{N}^d$ is the already defined $(i,j)\mapsto\nu(i,j)$ and $\iota(i,j)=i$ simply is the projection on the first coordinate.

In order to prove our wanted bounds we will need more general objects involving nonlocal Wick monomials.
Let 
\[
z=(z_{i,j})_{\substack{1\le i\le n\\ 1\le j\le r_i}}\in U\subset (\mathbb{R}^d)^{r_1+\cdots+r_n}
\]
where $U$ is the open set defined by the conditions $z_{i,j}\neq z_{i',j'}$ whenever $i\neq i'$.
We define the function $z\in U$ given by 
\[
\langle
:\partial^{\nu(1,1)}\phi(z_{1,1})\cdots\partial^{\nu(1,r_1)}\phi(z_{1,r_1}):\times\cdots\times
:\partial^{\nu(n,1)}\phi(z_{n,1})\cdots\partial^{\nu(n,r_n)}\phi(z_{n,r_n}):
\rangle=
\]
\begin{equation}
\sum_{\mathcal{W}\ {\rm off-diag}}
\prod_{\{a,b\}\in\mathcal{W}} C^{\nu(a),\nu(b)}(z_a,z_b)
\label{nonlocaldef}
\end{equation}
where we used the same choices for $F$, $(F_i)_{1\le i\le n}$ and $\nu$ as before.
Of course, the correlation $\langle\mathcal{O}_{A_1}(x_1)\cdots\mathcal{O}_{A_n}(x_n)\rangle$ corresponds to the specialization
\begin{equation}
\left\{
\begin{array}{ccccccc}
z_{1,1} & = & \cdots & = & z_{1,r_1} & = & x_1 \\
\vdots & & & & & & \\
z_{n,1} & = & \cdots & = & z_{n,r_n} & = & x_n\ .
\end{array}
\right.
\label{Zspecial}
\end{equation}

\begin{Lemma}
For all $\alpha\in\mathbb{N}^d$, $\delta\in\mathbb{R}$, $u\in\mathbb{R}^d\backslash\{0\}$,
one has the inequality
\[
|\partial_u^{\alpha}|u|^{-\delta}|\le O(1)\ |u|^{-\delta-\alpha}
\]
where the constant $O(1)$ can depend on $d,\delta$ and $\alpha$ but not on the position $u$.
\end{Lemma}

\noindent{\bf Proof:}
By a trivial induction on $|\alpha|$
one can write 
$\partial_u^{\alpha}(u^2)^{-\frac{\delta}{2}}=P_\alpha(u)\times (u^2)^{-\frac{\delta}{2}-\alpha}$
for some polynomial $P_{\alpha}(u)$ which is homogeneous of degree $|\alpha|$. The lemma then follows from
the obvious inequality $|P_{\alpha}(u)|\le O(1)|u|^{|\alpha|}$.
\qed

As a consequence of the lemma, one has for any points $x\neq y$, and any multiindices $\alpha,\beta$,
\begin{equation}
|C^{\alpha,\beta}(x,y)|\le O(1)\ |x-y|^{-2[\phi]-|\alpha|-|\beta|}\ .
\label{propagest}
\end{equation}
From now on we will use $O(1)$ to denote unspecified constants which are independent of the locations of the points concerned by the inequality at hand.

\begin{Proposition}\label{nonlocalbd}
For all $z\in U$ one has
\[
\left|\left\langle
\prod_{i=1}^{n}:\prod_{j=1}^{r_i}\partial^{\nu(i,j)}\phi(z_{i,j}):
\right\rangle\right|\le O(1)\times
\prod_{a\in F}\left(\min_{b\in F\backslash F_{\iota(a)}}|z_a-z_b|\right)^{-[\phi]-|\nu(a)|}\ .
\]
\end{Proposition}

\noindent{\bf Proof:}
It is enough to bound each term in the finite sum over off-diagonal $\mathcal{W}$'s.
By (\ref{propagest}), such a term is bounded by a constant times
\[
\prod_{\{a,b\}\in\mathcal{W}}\left(
|z_a-z_b|^{-[\phi]-|\nu(a)|}\times |z_b-z_a|^{-[\phi]-|\nu(b)|}
\right)
\le
\]
\[
\prod_{\{a,b\}\in\mathcal{W}}\left[
\left(\min_{c\in F\backslash F_{\iota(a)}}|z_a-z_c|\right)^{-[\phi]-|\nu(a)|}
\times 
\left(\min_{c\in F\backslash F_{\iota(b)}}|z_b-z_c|\right)^{-[\phi]-|\nu(b)|}
\right]
\]
since $\mathcal{W}$ is off-diagonal and the exponents are negative.
Finally, collecting the resulting factors according to elements $a$ in $F$ instead of pairs $\{a,b\}$
in the perfect matching $\mathcal{W}$, we obtain the desired bound.
\qed

\begin{Corollary}
One has the bound
\[
\left|\langle\mathcal{O}_{A_1}(x_1)\cdots\mathcal{O}_{A_n}(x_n)\rangle\right|\le O(1)
\prod_{i=1}^{n}\left(\min_{j\neq i} |x_i-x_j|\right)^{-[A_i]}
\]
and the BNNFB holds for the abstract system of correlations under consideration.
\end{Corollary}

\noindent{\bf Proof:}
Given that, for all $i$, $\sum_{a\in F_i}([\phi]+|\nu(a)|)=[A_i]$, the first inequality follows immediately from Proposition \ref{nonlocalbd} and the specialization (\ref{Zspecial}).
This gives a clean and sharp form of the BNNFB with $\epsilon=0$ and $k=0$ as expected for a conformal field theory. In order to derive the BNNFB as formulated in \S\ref{hardhyp} with some $\epsilon>0$ one can of course divide and multiply by factors of the form $|x-y|^\epsilon$ which are bounded by taking $k$ large enough, thanks to the elementary inequality $|x-y|\le \sqrt{2}\langle x\rangle\langle y\rangle$.
\qed

\begin{Corollary}\label{ballcoro}
Suppose $\mathbb{B}_1,\ldots,\mathbb{B}_n$ are disjoint closed Euclidean balls in $\mathbb{R}^d$.
Suppose $z=(z_{i,j})$ satisfies $z_{i,j}\in\mathbb{B}_i$ for all $i$, $1\le i\le n$ and for all $j$,
$1\le j\le r_i$. Then
\[
\left|\left\langle
\prod_{i=1}^{n}:\prod_{j=1}^{r_i}\partial^{\nu(i,j)}\phi(z_{i,j}):
\right\rangle\right|\le O(1)\times
\prod_{i=1}^{n}\left(\min_{l\neq i} d(\mathbb{B}_i,\mathbb{B}_l)
\right)^{-r_i[\phi]-\sum_{j=1}^{r_i}|\nu(i,j)|}
\]
where we used the standard notation $d(\mathbb{A},\mathbb{B})=\inf\{|x-y|\ |\ x\in\mathbb{A},\ y\in\mathbb{B}\}$
for the distance between two subsets of $\mathbb{R}^d$.
\end{Corollary}

\noindent{\bf Proof:}
For an element $a\in F_{i}$, we have $\iota(a)=i$ and $z_a\in \mathbb{B}_i$. If $b\in F\backslash F_{\iota(a)}$, then $b\in F_{l}$ for some $l\neq i$ and thus $z_b\in
\mathbb{B}_l$. Hence,
$|z_a-z_b|\ge \min_{l\neq i} d(\mathbb{B}_i,\mathbb{B}_l)$ and the corollary follows by a simple rearrangement of the bound from Proposition \ref{nonlocalbd}.
\qed

\subsection{The OPE structure}\label{fgfOPE}

We now discuss the OPE structure of the abstract system of pointwise correlations defined in (\ref{explicitdef}).
For given labels $D_1,\ldots,D_n\in\mathcal{A}^{\infty}$ and given distinct points $u_1,\ldots,u_n$, consider the correlator $\langle\mathcal{O}_{D_1}(u_1)\cdots\mathcal{O}_{D_n}(u_n)\rangle$
where some operators are grouped in pairs say $\{1,2\}$, $\{3,4\}$, etc. inside $[n]=\{1,\ldots,n\}$. The OPE arises when {\em conditioning} the sum over $\mathcal{W}$ according to the contractions $\{a,b\}\in\mathcal{W}$ which are contained in $F_1\cup F_2$, $F_3\cup F_4$, etc.
Let us focus on one pair say $\{1,2\}\subset [n]$.
In this case we condition the sum over $\mathcal{W}$ by fixing
$\mathcal{V}_{\mathcal{W}}=\{\{a,b\}\in\mathcal{W}\ |\ \{a,b\}\subset F_1\cup F_2\}$.
``Inside correlations'', one has the identity
\[
\mathcal{O}_{D_1}(u_1)\mathcal{O}_{D_2}(u_2)=\ 
:\partial^{\nu(1,1)}\phi(u_{1})\cdots\partial^{\nu(1,r_1)}\phi(u_{1}):\times
:\partial^{\nu(2,1)}\phi(u_{2})\cdots\partial^{\nu(2,r_2)}\phi(u_{2}):
\]
\[
=\sum_{\mathcal{V}}\prod_{\substack{(a,b)\in F_1\times F_2\\ \{a,b\}\in\mathcal{V}}}
C^{\nu(a),\nu(b)}(u_1,u_2)
:\prod_{a\in F_1\backslash(\cup\mathcal{V})}\partial^{\nu(a)}\phi(u_1)\times\prod_{b\in F_2\backslash(\cup\mathcal{V})}\partial^{\nu(b)}\phi(u_2):\ .
\]
The sum over $\mathcal{V}$ is over sets of disjoint (unordered) pairs contained in $F_1\cup F_2$ but not in $F_1$ nor $F_2$. Namely, $\mathcal{V}$ is a (not necessarily perfect) matching of $F_1\cup F_2$ where one never matches two elements of $F_1$ or two elements of $F_2$.
Unmatched elements of $F_1$ form the set $F_1\backslash(\cup\mathcal{V})$ whereas
unmatched elements of $F_2$ form the set $F_2\backslash(\cup\mathcal{V})$.
One can use the nonlocal generalization (\ref{nonlocaldef}) in order to write the more explicit identity
\[
\langle\mathcal{O}_{D_1}(u_1)\cdots\mathcal{O}_{D_n}(u_n)\rangle
=\sum_{\mathcal{V}}\left[\prod_{\substack{(a,b)\in F_1\times F_2\\ \{a,b\}\in\mathcal{V}}}
C^{\nu(a),\nu(b)}(u_1,u_2)\right]
\]
\[
\times\left\langle:\prod_{a\in F_1\backslash(\cup\mathcal{V})}\partial^{\nu(a)}\phi(u_1)\times\prod_{b\in F_2\backslash(\cup\mathcal{V})}\partial^{\nu(b)}\phi(u_2):
\times\mathcal{O}_{D_3}(u_3)\cdots\mathcal{O}_{D_n}(u_n)\right\rangle\ .
\]
The next step is to perform a Taylor expansion to sufficiently high order in the position $u_1$ with respect to the reference point $u_2$.
Still in shorthand notation, this amounts to setting
\[
f(t)=\ \ :\prod_{a\in F_1\backslash(\cup\mathcal{V})}\partial^{\nu(a)}\phi(u_2+t(u_1-u_2))\times\prod_{b\in F_2\backslash(\cup\mathcal{V})}\partial^{\nu(b)}\phi(u_2):
\]
and writing
\[
f(1)=\sum_{q=0}^{p_{\mathcal{V}}}\frac{f^{(q)}(0)}{q!}
+\int\limits_{0}^{1}dt\ \frac{(1-t)^{p_{\mathcal{V}}}}{p_{\mathcal{V}}!} f^{(p_{\mathcal{V}}+1)}(t)
\]
for some well chosen $p_{\mathcal{V}}$ to be defined shortly.
By the multivariate chain rule,
\[
f^{(q)}(t)=\sum_{\substack{\alpha\in\mathbb{N}^d\\ |\alpha|=q}}
\frac{q!}{\alpha!}(u_1-u_2)^{\alpha}
:\partial^{\alpha}\left[\prod_{a\in F_1\backslash(\cup\mathcal{V})}\partial^{\nu(a)}\phi\right]
(u_2+t(u_1-u_2))\times\prod_{b\in F_2\backslash(\cup\mathcal{V})}\partial^{\nu(b)}\phi(u_2):
\]
which by the multivariate Leibnitz rule gives
\[
f^{(q)}(t)=q!\sum_{\beta,||\beta||=q}
\frac{(u_1-u_2)^{|\beta|}}{\beta!}
:\prod_{a\in F_1\backslash(\cup\mathcal{V})}\partial^{\nu(a)+\beta(a)}\phi(u_2+t(u_1-u_2))\times\prod_{b\in F_2\backslash(\cup\mathcal{V})}\partial^{\nu(b)}\phi(u_2):
\]
where the following notation was used.
The sum is over maps $\beta:F_1\backslash(\cup\mathcal{V})\rightarrow\mathbb{N}^{d}$. Therefore, $\beta$ is not a multiindex but a collection of multiindices $\beta(a)$, one for each unmatched element $a$ in
$F_1\backslash(\cup\mathcal{V})$. We thus introduced the new notation
\[
\beta!=\prod_{F_1\backslash(\cup\mathcal{V})}\beta(a)!\ \in \mathbb{N}\ ,
\ \ |\beta|=\sum_{F_1\backslash(\cup\mathcal{V})}\beta(a)\ \in \mathbb{N}^d
\ \ {\rm and}\ \ 
||\beta||=\sum_{F_1\backslash(\cup\mathcal{V})}|\beta(a)|\ \in \mathbb{N}\ .
\]
As a result, we have ``inside correlations'' the identity
\[
:\partial^{\nu(1,1)}\phi(u_{1})\cdots\partial^{\nu(1,r_1)}\phi(u_{1}):\times
:\partial^{\nu(2,1)}\phi(u_{2})\cdots\partial^{\nu(2,r_2)}\phi(u_{2}):
\]
\[
=\sum_{\mathcal{V}}\prod_{\substack{(a,b)\in F_1\times F_2\\ \{a,b\}\in\mathcal{V}}}
C^{\nu(a),\nu(b)}(u_1,u_2)\times
\sum_{\beta,||\beta||\le p_{\mathcal{V}}}\frac{(u_1-u_2)^{|\beta|}}{\beta!}
\]
\[
\times
:\prod_{a\in F_1\backslash(\cup\mathcal{V})}\partial^{\nu(a)+\beta(a)}\phi(u_2)\times\prod_{b\in F_2\backslash(\cup\mathcal{V})}\partial^{\nu(b)}\phi(u_2):
\]
\[
+\sum_{\mathcal{V}}\prod_{\substack{(a,b)\in F_1\times F_2\\ \{a,b\}\in\mathcal{V}}}
C^{\nu(a),\nu(b)}(u_1,u_2)\times
\sum_{\beta,||\beta||=p_{\mathcal{V}}+1}\frac{(u_1-u_2)^{|\beta|}}{\beta!}
\]
\begin{equation}
\times
(p_{\mathcal{V}}+1)
\int\limits_{0}^{1}dt\ (1-t)^{p_{\mathcal{V}}}
:\prod_{a\in F_1\backslash(\cup\mathcal{V})}\partial^{\nu(a)+\beta(a)}\phi(u_2+t(u_1-u_2))\times\prod_{b\in F_2\backslash(\cup\mathcal{V})}\partial^{\nu(b)}\phi(u_2):\ .
\label{fullTaylor}
\end{equation}
Suppose one is performing the OPE up to a scaling dimension cut-off $\Theta$. Then we make a minimal choice for $p_{\mathcal{V}}$, namely, we let
\[
p_{\mathcal{V}}=\max\left\{
0\ ,\ \left\lfloor\ \Theta-|(F_1\cup F_2)\backslash(\cup\mathcal{V})|
\times [\phi]-\sum_{a\in (F_1\cup F_2)\backslash(\cup\mathcal{V})}|\nu(a)|\ \right\rfloor
\right\}\ .
\]
We are now able to introduce the needed OPE coefficients as suggested by the terms
appearing in the first sum of (\ref{fullTaylor}). We let, {\em by definition},
\begin{equation}
\mathcal{C}_{D_1,D_2}^{E}(u_1,u_2)=
\sum_{\mathcal{V}}\prod_{\substack{(a,b)\in F_1\times F_2\\ \{a,b\}\in\mathcal{V}}}
C^{\nu(a),\nu(b)}(u_1,u_2)\times
\sum_{\beta\in (\mathbb{N}^d)^{F_1\backslash(\cup\mathcal{V})}}\frac{(u_1-u_2)^{|\beta|}}{\beta!}
\times\bbone\{(\mathcal{V},\beta)\rightarrow E\}
\label{OPEcoeffdef}
\end{equation}
where the condition $(\mathcal{V},\beta)\rightarrow E$ means the symbolic equality
\[
:\prod_{a\in F_1\backslash(\cup\mathcal{V})}\partial^{\nu(a)+\beta(a)}\phi\times\prod_{b\in F_2\backslash(\cup\mathcal{V})}\partial^{\nu(b)}\phi:\ =\mathcal{O}_{E}\ .
\]
More precisely, construct the finite set $G=(F_1\backslash(\cup\mathcal{V}))\cup(F_2\backslash(\cup\mathcal{V}))$
and the map $\omega:G\rightarrow\mathbb{N}^d$ given by
\[
\omega(a)=
\left\{
\begin{array}{ll}
\nu(a)+\beta(a) &  {\rm if}\ a\in F_1\backslash(\cup\mathcal{V})\ , \\
\nu(a) &  {\rm if}\ a\in F_2\backslash(\cup\mathcal{V})\ .
\end{array}
\right.
\]
The condition $(\mathcal{V},\beta)\rightarrow E$, by definition, means that
$\{G;\omega\}=E$.
Notice that the previous constraint $||\beta||\le p_{\mathcal{V}}$ which featured in (\ref{fullTaylor})
is not included in (\ref{OPEcoeffdef}).

In fact, one can organize the first sum in (\ref{fullTaylor}) according to the unique label $E\in\mathcal{A}$ 
which satisfies $(\mathcal{V},\beta)\rightarrow E$. One then splits this sum according to whether $[E]\le \Theta$
or $[E]>\Theta$.
Hence
\begin{equation}
\mathcal{O}_{D_1}(u_1)\mathcal{O}_{D_2}(u_2)=\sum_{E,[E]\le \Theta}
\mathcal{C}_{D_1,D_2}^{E}(u_1,u_2)\mathcal{O}_{E}(u_2)\ + {\rm Rem}_1 +{\rm Rem}_2
\label{rem12defeq}
\end{equation}
where the first remainder is
\[
{\rm Rem}_1=\sum_{E,[E]>\Theta}
\sum_{\mathcal{V}}\prod_{\substack{(a,b)\in F_1\times F_2\\ \{a,b\}\in\mathcal{V}}}
C^{\nu(a),\nu(b)}(u_1,u_2)\times
\sum_{\beta,||\beta||\le p_{\mathcal{V}}}\frac{(u_1-u_2)^{|\beta|}}{\beta!}
\]
\begin{equation}
\times
:\prod_{a\in F_1\backslash(\cup\mathcal{V})}\partial^{\nu(a)+\beta(a)}\phi(u_2)\times\prod_{b\in F_2\backslash(\cup\mathcal{V})}\partial^{\nu(b)}\phi(u_2):\ \times
\bbone\{(\mathcal{V},\beta)\rightarrow E\}
\label{rem1defeq}
\end{equation}
and the second one is
\[
{\rm Rem}_2=\sum_{\mathcal{V}}\prod_{\substack{(a,b)\in F_1\times F_2\\ \{a,b\}\in\mathcal{V}}}
C^{\nu(a),\nu(b)}(u_1,u_2)\times
\sum_{\beta,||\beta||=p_{\mathcal{V}}+1}\frac{(u_1-u_2)^{|\beta|}}{\beta!}
\]
\begin{equation}
\times
(p_{\mathcal{V}}+1)
\int\limits_{0}^{1}dt\ (1-t)^{p_{\mathcal{V}}}
:\prod_{a\in F_1\backslash(\cup\mathcal{V})}\partial^{\nu(a)+\beta(a)}\phi(u_2+t(u_1-u_2))\times\prod_{b\in F_2\backslash(\cup\mathcal{V})}\partial^{\nu(b)}\phi(u_2):\ .
\label{rem2defeq}
\end{equation}
Note that the first sum in (\ref{rem12defeq}) should, in principle,
include the constraint $||\beta||\le p_{\mathcal{V}}$. However, the latter is redundant and can be removed which leads to the neater expression $\sum_{E,[E]\le \Theta}
\mathcal{C}_{D_1,D_2}^{E}(u_1,u_2)\mathcal{O}_{E}(u_2)$.
Indeed, the condition $(\mathcal{V},\beta)\rightarrow E$ implies 
$[E]=\delta+||\beta||$
where
\[
\delta=|(F_1\cup F_2)\backslash(\cup\mathcal{V})|\times[\phi]
+\sum_{a\in (F_1\cup F_2)\backslash(\cup\mathcal{V})}|\nu(a)|\ .
\]
Since the condition $[E]\le\Theta$ is in force,
the nonnegative integer $||\beta||$ must satisfy $||\beta||\le \Theta-\delta$
and therefore
$||\beta||\le \lfloor\Theta-\delta\rfloor
\le\max\{0,\lfloor\Theta-\delta\rfloor\}=p_{\mathcal{V}}$.
As a result, the OPE element for the product $\mathcal{O}_{D_1}(u_1)\mathcal{O}_{D_2}(u_2)$,
with $u_2$ as point of reference and subtracting all field labels $E$ with scaling dimension $[E]\le \Theta$,
is given by
\begin{equation}
{\rm OPE}(u_1,u_2)={\rm Rem}_1+{\rm Rem}_2\ .
\label{OPEremains}
\end{equation}

\subsection{A proof of the ENNFB}

We will now show that the abstract system of pointwise correlations corresponding to the fractional Gaussian field with scaling dimension $[\phi]>0$ satisfies the ENNFB as stated in \S\ref{hardhyp}.
Reverting to the notations of \S\ref{hardhyp}, we now pick some $\eta\in(0,\frac{1}{3}]$ and first establish the bound with $\epsilon=0$ and $k=0$.
For $1\le i\le m+n+p$, we let $R_i=\eta\min_{l\neq i}|x_i-x_l|$, where the minimum is over all $l$'s such that $1\le l\le m+n+p$ and $l\neq i$.
We let $\mathbb{B}_i$ denote the closed Euclidean ball centered at $x_i$ and of radius
$R_i$. We also let $\lambda_i=\min_{l\neq i}d(\mathbb{B}_i,\mathbb{B}_l)$.
If $i\neq j$, then clearly
$R_i\le \eta|x_i-x_j|$
and
$R_j\le \eta|x_i-x_j|$.
From this simple observation and the fact
$\eta<\frac{1}{2}$, we see that the balls $\mathbb{B}_i$ are all disjoint.
Also note that from elementary geometry we have
\[
d(\mathbb{B}_i,\mathbb{B}_j)=|x_i-x_j|-R_i-R_j\ge (1-2\eta)|x_i-x_j|\ .
\]
As a result, we have
\begin{equation}
\lambda_i\ge (1-2\eta)\min_{l\neq i}|x_i-x_l|=\frac{(1-2\eta)}{\eta} R_i
\label{lambdaineq}
\end{equation}
and hence
\begin{equation}
\frac{R_i}{\lambda_i}\le\frac{\eta}{1-2\eta}\le 1
\label{geometric}
\end{equation}
because we picked $\eta\le\frac{1}{3}$.

Regarding the EFNNB, 
there is nothing to prove unless $y_i\in \mathbb{B}_i$ for all $i$, $1\le i\le m+n$, which we now assume.
In order to bound 
\[
{\rm LHS}=
\left\langle
\prod_{i=1}^{m} {\rm OPE}_i(y_i,x_i)
\prod_{i=m+1}^{m+n} {\rm CZ}_i(y_i,x_i)
\prod_{i=m+n+1}^{m+n+p} \mathcal{O}_{B_i}(x_i)
\right\rangle\ ,
\]
we decompose each ${\rm OPE}$ factor as ${\rm Rem}_1+{\rm Rem}_2$ following (\ref{OPEremains}).
For any CZ term $\mathcal{O}_{B_i}(y_i)-\mathcal{O}_{B_i}(x_i)$
with say $B_i=[\nu(1),\ldots,\nu(r)]$, we write
\[
\mathcal{O}_{B_i}(y_i)-\mathcal{O}_{B_i}(x_i)=
:\partial^{\nu(1)}\phi(y_i)\cdots\partial^{\nu(r)}\phi(y_i):
-:\partial^{\nu(1)}\phi(x_i)\cdots\partial^{\nu(r)}\phi(x_i):
\]
\[
=\int\limits_{0}^{1}dt\ \frac{d}{dt}
:\partial^{\nu(1)}\phi(x_i+t(y_i-x_i))\cdots\partial^{\nu(r)}\phi(x_i+t(y_i-x_i)):
\]
\begin{equation}
=\sum_{\beta,||\beta||=1}(y_i-x_i)^{|\beta|}\int\limits_{0}^{1}dt
\ :\prod_{a=1}^{r}\partial^{\nu(a)+\beta(a)}\phi(x_i+t(y_i-x_i)):
\label{CZgradeq}
\end{equation}
where the sum is over collections of multiindices $\beta=(\beta(a))_{1\le a\le r}\in(\mathbb{N}^d)^r$ and we used the notations of the previous section for $|\beta|\in\mathbb{N}^d$ and $||\beta||\in\mathbb{N}$.

Implementing the expansion described above 
produces for LHS a rather formidable yet finite sum which we will not write explicitly.
Each term is then bounded using Corollary \ref{ballcoro}. The latter gives a constant times a product of factors which can be collected by successively examining each OPE, CZ and spectator factor. The reader might find it useful to go over the following
list of cases twice: once to check the applicability of Corollary \ref{ballcoro} by making sure all the points belong to the proper balls $\mathbb{B}$, and once to review the collection of factors produced by the application of Corollary \ref{ballcoro}.

\smallskip
\noindent{\bf 1st case:} OPE factor of ${\rm Rem}_1$ type.

First let us consider a ${\rm Rem}_1$ contribution for an OPE term ${\rm OPE}_i(y_i,x_i)$, $1\le i\le m$, which has further been decomposed using a sum over $E,\mathcal{V},\beta$ with in fact $E$ determined by $\mathcal{V},\beta$ via the condition $(\mathcal{V},\beta)\rightarrow E$.

As a matter of managing notations, we first need to connect the globally defined variables from \S\ref{hardhyp}
with the locally defined variables used in \S\ref{fgfOPE}.
We thus set $u_1=y_i$, $u_2=x_i$, $D_1=A_i$, $D_2=B_i$ as well as $\Theta=\Delta_i$.
From (\ref{rem1defeq}), we see that we have a first expression
\[
\prod_{\substack{(a,b)\in F_1\times F_2\\ \{a,b\}\in\mathcal{V}}}
C^{\nu(a),\nu(b)}(u_1,u_2)\times
(u_1-u_2)^{|\beta|}\ ,
\]
outside the big correlator, which is
readily bounded using (\ref{propagest}) by
\begin{equation}
O(1)\times |u_1-u_2|^{-2|\mathcal{V}|[\phi]-\sum_{a\in \cup\mathcal{V}}|\nu(a)|+||\beta||}\ .
\label{fgftypeI}
\end{equation}
We also have a second expression
\[
:\prod_{a\in F_1\backslash(\cup\mathcal{V})}\partial^{\nu(a)+\beta(a)}\phi(u_2)\times\prod_{b\in F_2\backslash(\cup\mathcal{V})}\partial^{\nu(b)}\phi(u_2):\ ,
\]
inside the big correlator, which is bounded using
Corollary \ref{ballcoro}. The corresponding contribution to the global bound is
\begin{equation}
O(1)\times \lambda_{i}^{-|(F_1\cup F_2)\backslash(\cup\mathcal{V})|[\phi]
-\sum_{a\in(F_1\cup F_2)\backslash(\cup\mathcal{V})}|\nu(a)|-||\beta||}\ .
\label{fgftypeII}
\end{equation}
Note that the Wick monomial only involves $u_2=x_i\in\mathbb{B}_i$.

Recall that the condition $(\mathcal{V},\beta)\rightarrow E$
imposes
\[
[E]=|(F_1\cup F_2)\backslash(\cup\mathcal{V})|\times [\phi]
+\sum_{a\in(F_1\cup F_2)\backslash(\cup\mathcal{V})}|\nu(a)|+||\beta||
\]
while we also have
\[
[D_1]=|F_1|[\phi]+\sum_{a\in F_1}|\nu(a)|\ \ 
{\rm and}
\ \ 
[D_2]=|F_2|[\phi]+\sum_{a\in F_2}|\nu(a)|\ .
\]
Therefore
\begin{equation}
[D_1]+[D_2]-[E]=2|\mathcal{V}|[\phi]+\sum_{a\in \cup\mathcal{V}}|\nu(a)|-||\beta||
\label{dimsumeq}
\end{equation}
and the product of (\ref{fgftypeI}) and (\ref{fgftypeII})
gives an overall factor
\[
O(1) \times |u_1-u_2|^{[E]-[D_1]-[D_2]}\times \lambda_{i}^{-[E]}\ .
\]
By (\ref{geometric}), and the hypothesis $y_i\in\mathbb{B}_i$, we have
\[
\frac{|u_1-u_2|}{\lambda_i}\le \frac{R_i}{\lambda_i}\le 1\ .
\]
Also recall the condition $[E]>\Theta$ which implies $[E]\ge {\rm next}(\Theta)$.
Hence
\[
\left(\frac{|u_1-u_2|}{\lambda_i}\right)^{[E]}\le
\left(\frac{|u_1-u_2|}{\lambda_i}\right)^{{\rm next}(\Theta)}\ .
\]
We use the inequality in (\ref{lambdaineq}) and the fact
${\rm next}(\Theta)\ge 0$ (even if $\Theta<0$) to bound
$\lambda_{i}^{-{\rm next(\Theta)}}$ by
$(\min_{l\neq i}|x_i-x_l|)^{-{\rm next}(\Theta)}$
times
a harmless factor
$(1-2\eta)^{-{\rm next}(\Theta)}$ which gets absorbed
into $O(1)$.
Reverting back to global variables,
we see that the present case contributes
an overall factor
\[
O(1)\times |y_i-x_i|^{{\rm next}(\Delta_i)-[A_i]-[B_i]}
\times \left(\min_{j\neq i} |x_i-x_j|\right)^{-{\rm next}(\Delta_i)}
\]
to our final bound.

\smallskip
\noindent{\bf 2nd case:} OPE factor of ${\rm Rem}_2$ type.

Now let us consider a ${\rm Rem}_2$ contribution for an OPE term ${\rm OPE}_i(y_i,x_i)$, $1\le i\le m$, which has further been decomposed using a sum over $\mathcal{V},\beta$. While there is no $E$ explicitly appearing
in the formula (\ref{rem2defeq}), there is no harm in reintroducing it by defining $E$ using the condition
$(\mathcal{V},\beta)\rightarrow E$.
We will also take a supremum over $t\in [0,1]$ in our estimates.
Note that $y_i,x_i$ or rather $u_1,u_2$ are both in $\mathbb{B}_i$ and so is $u_2+t(u_1-u_2)$ by convexity. 
Namely, all points involved in the Wick monomial (which now is nonlocal) are contained in $\mathbb{B}_i$.
Note that $[E]=\delta+||\beta||$
with
\[
\delta=|(F_1\cup F_2)\backslash(\cup\mathcal{V})|\times [\phi]
+\sum_{a\in(F_1\cup F_2)\backslash(\cup\mathcal{V})}|\nu(a)|\ .
\]
However, now
\[
||\beta||=p_{\mathcal{V}}+1=\max\{0,\lfloor\Theta-\delta\rfloor\}+1\ge
\lfloor\Theta-\delta\rfloor+1>\Theta-\delta
\]
which implies $[E]>\Theta$ and therefore $[E]\ge{\rm next}(\Theta)$.
Now the same reasoning as in the previous case gives the same overall contribution to the final bound.

\smallskip
\noindent{\bf 3rd case:} CZ factor.

Let us now consider a CZ factor ${\rm CZ}_i(y_i,x_i)$, $m+1\le i\le m+n$.
Using the notation of (\ref{CZgradeq}), we now have a finite sum to bound.
Each term gives a factor outside the big correlator which is trivially bounded by
\begin{equation}
|(y_i-x_i)^{|\beta|}|\le |y_i-x_i|^{||\beta||}=|y_i-x_i|
\label{fgfCZtypeI}
\end{equation}
since $||\beta||=1$.
We also have a Wick monomial
\[
:\prod_{a=1}^{r}\partial^{\nu(a)+\beta(a)}\phi(x_i+t(y_i-x_i)):
\]
inside the big correlator.
We again make the easy yet important remark that $x_i+t(y_i-x_i)\in\mathbb{B}_i$, by convexity.
Via the bound provided by Corollary \ref{ballcoro}, this Wick monomial contributes
a factor
\begin{equation}
O(1)\times \lambda_{i}^{-r[\phi]-\sum_{a=1}^{r}|\nu(a)|-||\beta||}=\lambda_{i}^{-[B_i]-1}\ .
\label{fgfCZtypeII}
\end{equation}
Combining (\ref{fgfCZtypeI}) and (\ref{fgfCZtypeII}) and by similar geometric reasoning as in the 1st case
we thus obtain an overall factor
\[
O(1)\times
|y_i-x_i|\times \left(\min_{j\neq i}|x_i-x_j|\right)^{-[B_i]-1}
\]
for our final bound.

\smallskip
\noindent{\bf 4th case:} Spectator factor.

We finally consider a spectator
factor $\mathcal{O}_{B_i}(x_i)$, $m+n\le i\le m+n+p$.
Trivially $x_i\in\mathbb{B}_i$ and Corollary \ref{ballcoro} immediately gives the desired
factor
\[
O(1)\times
\left(\min_{j\neq i}|x_i-x_j|\right)^{-[B_i]}\ .
\]

\medskip
This concludes the case discussion. Collecting all the factors considered above, we obtain:
\[
\prod_{i=1}^{m+n}\bbone\left\{
|y_i-x_i|\le\eta\min_{j\neq i}|x_i-x_j|
\right\}\times
\left|
\left\langle
\prod_{i=1}^{m} {\rm OPE}_i(y_i,x_i)
\prod_{i=m+1}^{m+n} {\rm CZ}_i(y_i,x_i)
\prod_{i=m+n+1}^{m+n+p} \mathcal{O}_{B_i}(x_i)
\right\rangle
\right|\le
\]
\[
O(1)\times
\prod_{i=1}^{m}\left\{|y_i-x_i|^{\Delta_i+\gamma_i-[A_i]-[B_i]}\times\left(\min_{j\neq i}|x_i-x_j|\right)^{-\Delta_i-\gamma_i}\right\}
\]
\[
\times
\prod_{i=m+1}^{m+n}\left\{|y_i-x_i|^{\gamma_i}\times\left(\min_{j\neq i}|x_i-x_j|\right)^{-[B_i]-\gamma_i}\right\}
\times
\prod_{i=m+n+1}^{m+n+p}\left(\min_{j\neq i}|x_i-x_j|\right)^{-[B_i]}
\]
where
\[
\gamma_i=\left\{
\begin{array}{lll}
{\rm next}(\Delta_i)-\Delta_i & {\rm if} & 1\le i\le m\ , \\
1 & {\rm if} & m+1\le i\le m+n\ .
\end{array}
\right.
\]
This is a clean and sharp form of the ENNFB which is especially suitable for a conformal field theory.
Using the worst case scenario planning from \S\ref{worstsec}
one can set
\[
\gamma=\min\left\{1,\min_{1\le i\le m+n} \gamma_i\right\}>0
\]
and one can also turn on $\epsilon>0$, at the expense of a suitable power $k$ of inhomogeneous norm factors
$\langle\cdot\rangle$,
in order to establish the ENNFB exactly as it is stated in \S\ref{hardhyp}. \qed

\subsection{A proof of the soft and hard hypotheses on the OPE coefficients}\label{softhardsec}

In order to use Theorem \ref{maintheorem}, we need to check that the present system of pointwise correlations
with its OPE structure satisfies the soft and hard hypotheses of \S\ref{softhyp} and \S\ref{hardhyp}.

For any $D_1,D_2,E\in\mathcal{A}^{\infty}$ and using the definition (\ref{OPEcoeffdef})
of an OPE coefficient $\mathcal{C}_{D_1,D_2}^{E}(u_1,u_2)$, one immediately sees that this is indeed a $C^{\infty}$ function
on ${\rm Conf}_2$. Moreover, each term in the sum over $\mathcal{V},\beta$, is easily bounded using (\ref{propagest}) by
\[
O(1)\times |u_1-u_2|^{-2|\mathcal{V}|[\phi]-\sum_{a\in\cup\mathcal{V}}|\nu(a)|+||\beta||}\ .
\]
Since the relation (\ref{dimsumeq}) also holds here, we obtain
\begin{equation}
|\mathcal{C}_{D_1,D_2}^{E}(u_1,u_2)|\le O(1)\times |u_1-u_2|^{[E]-[D_1]-[D_2]}
\label{OPEcoeffsharp}
\end{equation}
which is the version of (\ref{hardCbd}) with $\epsilon=0$ and $k=0$.
On can again turn on $\epsilon>0$ and adjust $k$ accordingly, so that the hypothesis (\ref{hardCbd}) holds.

We now restrict ourselves to $D_1,D_2,E$ in $\mathcal{A}=\{A\in\mathcal{A}^{\infty}\ |\ [A]<\frac{d}{2}\}$ and examine the requirements of \S\ref{softhyp}.
Since $[D_1]+[D_2]-[E]\le [D_1]+[D_2]<d$, the bound (\ref{OPEcoeffsharp}) guarantees the local integrability
of $\mathcal{C}_{D_1,D_2}^{E}(u_1,u_2)$ on the diagonal $u_1=u_2$.
Moreover, the condition (\ref{softbound}) here is trivial.

We now consider the condition (SH1) from \S\ref{softhyp}
and take Schwartz functions $f(u_1)$ and $g(u_2)$. Using the translation invariance
$\mathcal{C}_{D_1,D_2}^{E}(u_1,u_2)=\mathcal{C}_{D_1,D_2}^{E}(u_1-u_2,0)$
and a trivial change of variables, we have 
\[
g(u_2)\langle\mathcal{C}_{D_1,D_2}^{E}(u_1,u_2),f(u_1)\rangle_{u_1}
=g(u_2)\int_{\mathbb{R}^d\backslash\{0\}}dv
\ \mathcal{C}_{D_1,D_2}^{E}(v,0)\ f(u_2+v)\ .
\]
Let $k\in\mathbb{N}$ and $\alpha\in\mathbb{N}^d$.
By the Leibniz rule and the theorem of derivation under the integral sign, the estimation of 
$\langle u_2\rangle^k\ \partial_{u_2}^{\alpha}\left[
g(u_2)\langle\mathcal{C}_{D_1,D_2}^{E}(u_1,u_2),f(u_1)\rangle_{u_1}
\right]$
amounts to that of finitely many terms of the form
\[
\langle u_2\rangle^k\ \partial^{\beta}g(u_2)
\int_{\mathbb{R}^d\backslash\{0\}}dv
\ \mathcal{C}_{D_1,D_2}^{E}(v,0)\ \partial^{\alpha-\beta}f(u_2+v)
\]
where $\beta$ is a multiindex which is bounded component-wise by $\alpha$.
Let $p,q\in\mathbb{N}$ and introduce corresponding Schwartz seminorms for $f$ and $g$.
Putting absolute values and using (\ref{OPEcoeffsharp}), the previous expression is bounded by
\[
O(1)\times \langle u_2\rangle^{k-p}\ ||g||_{\beta,p}
\int_{\mathbb{R}^d\backslash\{0\}}dv\ 
\frac{||f||_{\alpha-\beta,q}}{|v|^{[D_1]+[D_2]-[E]}\times \langle u_2+v\rangle^q}\ .
\]
Form the elementary inequality $\langle x+y\rangle\le\sqrt{2}
\langle x\rangle \langle y\rangle$
we get
$\frac{1}{\langle u_2+v\rangle}\le \sqrt{2}\times
\frac{\langle u_2\rangle}{\langle v\rangle}$
which results in the new bound
\[
O(1)\times \langle u_2\rangle^{k+q-p}\ ||g||_{\beta,p}||f||_{\alpha-\beta,q}
\int_{\mathbb{R}^d\backslash\{0\}}dv\ 
\frac{1}{|v|^{[D_1]+[D_2]-[E]}\times \langle v\rangle^q}\ .
\]
Taking a supremum over $u_2$ while choosing
$q=d+1$, $p=k+d+1$, and the earlier inequality $[D_1]+[D_2]-[E]<d$
show that the function
$u_2\mapsto g(u_2)\langle\mathcal{C}_{D_1,D_2}^{E}(u_1,u_2),f(u_1)\rangle_{u_1}$
is in Schwartz space and depends continuously on $f$. Hence (SH1) is established. 
The proof of (SH2) is similar and left to the reader who might find it convenient to use the concatenation notation from the beginning of \S\ref{pinsum} and elementary inequalities between $\langle x_1,\ldots,x_n\rangle$ and $\langle x_1\rangle\cdots\langle x_n\rangle$.

\subsection{Application of the main theorem to the construction of local Wick monomials as random distributions}

Suppose we have an $\mathcal{S}'(\mathbb{R}^d)$-valued random variable $\phi$ on a probability space $(\Omega,\mathcal{F},\mathbb{P})$ and with law given by the FMFF
with $[\phi]>0$.  
We add the field $\mathcal{O}_{\bbone}$ which is constant in space and with respect to the random sample and equal to one. Provided $[\phi]<\frac{d}{2}$, these two fields form an incarnation of the system of pointwise correlations corresponding to the 
label subset $\{\bbone,\phi\}\subset\mathcal{A}$.
Since derivatives in the sense of distributions $\partial^{\alpha}:\mathcal{S}'(\mathbb{R}^d)\rightarrow
\mathcal{S}'(\mathbb{R}^d)$ are continuous and thus Borel measurable, one immediately gets $\partial^{\alpha}\phi$, for all $\alpha\in\mathbb{N}^d$,
as honest random Schwartz distributions on the same probability space. However, only for the subset
corresponding to $[\phi]+|\alpha|<\frac{d}{2}$ are we guaranteed to have an incarnation for the corresponding pointwise correlations.
Indeed, the needed identity for arbitrary test functions $f_1,\ldots,f_n$ and sequence of multiindices $\alpha(1),\ldots,\alpha(n)$
is
\[
\mathbb{E} \left[
\prod_{i=1}^{n} (-1)^{|\alpha(i)|}\phi(\partial^{\alpha(i)}f_i)
\right]
=\int_{{\rm Conf}_n}dx_1\ldots dx_n\ 
\langle \partial^{\alpha(1)}\phi(x_1)\cdots\partial^{\alpha(n)}\phi(x_n) 
\rangle\ f_1(x_1)\cdots f_n(x_n)\ .
\]
By the definition (\ref{defcorrel})
of $\langle \partial^{\alpha(1)}\phi(x_1)\cdots\partial^{\alpha(n)}\phi(x_n) \rangle$, this identity reduces to the $n=2$ case which amounts to the integration by parts
\[
(-1)^{|\alpha(1)|+|\alpha(2)|}
\int_{{\rm Conf}_2}dx_1 dx_2\ \frac{1}{|x_1-x_2|^{2[\phi]}}\ \partial^{\alpha(1)}f_1(x_1)\ \partial^{\alpha(1)}f_2(x_2) 
= 
\]
\[
\int_{{\rm Conf}_2}dx_1 dx_2\ 
\left[\partial_{x_1}^{\alpha(1)}\partial_{x_2}^{\alpha(2)}
\frac{1}{|x_1-x_2|^{2[\phi]}}\right]
f_1(x_1) f_2(x_2) \ .
\]
It is an easy exercise to check that this equality holds without boundary terms at infinity or the origin because our restriction on labels imposes the local integrability condition $2[\phi]+|\alpha(1)|+
|\alpha(2)|<d$.

The construction of random distributions giving an incarnation of Wick monomials indexed by $\mathcal{A}$ is done recursively with respect to the degree of the monomial.
Consider a new field to be constructed 
$\mathcal{O}_{C_{\ast}}=:\partial^{\alpha(1)}\phi
\cdots\partial^{\alpha(n)}\phi:$ with $[C_{\ast}]=n[\phi]+|\alpha(1)|+\cdots+|\alpha(n)|<\frac{d}{2}$.
It is obtained by Theorem \ref{maintheorem} applied to the OPE given by the product $\mathcal{O}_A\mathcal{O}_{B}$
where $A=[\alpha(1),\cdots,\alpha(n-1)]$ and $B=[\alpha(n)]$. Using the definition (\ref{OPEcoeffdef}), it is easy to check that in this particular case $\mathcal{C}_{A,B}^{C_{\ast}}(x_1,x_2)=1$
and the non-degeneracy condition  (\ref{nondeg}) holds trivially.

\subsection{Global conformal invariance} 

Global conformal invariance for the elementary scalar field $\phi$
means that we have the pointwise correlation identity
\[
|J_f(x_1)|^{\frac{[\phi]}{d}}\cdots |J_f(x_n)|^{\frac{[\phi]}{d}}
\langle\phi(f(x_1))\cdots\phi(f(x_n))\rangle=
\langle\phi(x_1)\cdots\phi(x_n)\rangle
\]
for all global conformal transformations $f$ of $\mathbb{R}^d$, and all collections of distinct points which do not get mapped to the point at infinity.
We used $J_f(x)$ to denote the Jacobian determinant of $f$ at $x$.
It is clearly sufficient to show the identity for generators of the global conformal group, e.g., Euclidean isometries, homotheties and the unit sphere inversion $I(x)=\frac{1}{|x|^2}x$. From the definition (\ref{defcorrel}), it is also clear that the case of general $n$ reduces to $n=2$, and only the case of the transformation $I$ needs explaining.
In components, the Jacobian matrix of $I$ is given by
$|x|^{-2}(
\delta_{ij}-2x_i x_j|x|^{-2})$.
A simple linear algebra computation shows that the determinant is $J_I(x)=-|x|^{-2d}$ and the needed
identity follows from the elementary equality
$|I(x_1)-I(x_2)|=|x_1-x_2|\times|x_1|^{-1}\times|x_2|^{-1}$.

\medskip
\noindent{\bf Acknowledgements:}
{\small
This work would have been impossible without a semester leave from teaching and administrative duties in the form
of a Sesquicentennial Associate Award. Therefore, the support of the Mathematics Department and the College
and Graduate School of Arts and Sciences at the
University of Virginia is very gratefully acknowledged. For useful discussions or correspondence we thank A. Chandra, J. Dubedat,
J. Fageot, M. Furlan, C. Garban,
M. Gubinelli, M. Hairer, S. Hollands, C. Hongler, C. Kopper, A. Kupiainen,
P. Mitter,
J.-C. Mourrat, J. Oikarinen,
E. Peltola, V. Rychkov, E. Seiler, D. Simmons-Duffin, V. Vargas and E. Witten.
Finally, we thank the anonymous referees for suggestions which helped improve this article.}

\end{document}